\documentclass[preprint,11pt]{elsarticle}

\usepackage{ifpdf}
\usepackage{amsmath}
\usepackage{amsfonts}
\usepackage{amsthm}
\usepackage{mathrsfs}
\usepackage{color}
\usepackage{graphicx}
\usepackage{subfigure}
\usepackage{float}
\usepackage{overpic}
\usepackage{threeparttable}
\usepackage{epstopdf}
\usepackage{dcolumn}
\usepackage{multirow}
\usepackage{booktabs}
\biboptions{numbers,sort&compress}

\usepackage{graphicx,natbib,amssymb,lineno}
\ifpdf
\usepackage[%
  pdftitle={Instructions for use of the document class
    elsart},%
  pdfauthor={},%
  pdfsubject={The preprint document class elsart},%
  pdfkeywords={instructions for use, elsart, document class},%
  pdfstartview=FitH,%
  bookmarks=true,%
  bookmarksopen=true,%
  breaklinks=true,%
  colorlinks=true,%
  linkcolor=blue,anchorcolor=blue,%
  citecolor=blue,filecolor=blue,%
  menucolor=blue,pagecolor=blue,%
  urlcolor=blue]{hyperref}
\else
\usepackage[%
  breaklinks=true,%
  colorlinks=true,%
  linkcolor=blue,anchorcolor=blue,%
  citecolor=blue,filecolor=blue,%
  menucolor=blue,pagecolor=blue,%
  urlcolor=blue]{hyperref}
\fi

\makeatletter
\def\elsartstyle{%
    \def\normalsize{\@setfontsize\normalsize\@xiipt{14.5}}
    \def\small{\@setfontsize\small\@xipt{13.6}}
    \let\footnotesize=\small
    \def\large{\@setfontsize\large\@xivpt{18}}
    \def\Large{\@setfontsize\Large\@xviipt{22}}
    \skip\@mpfootins = 18\p@ \@plus 2\p@
    \normalsize
} \@ifundefined{square}{}{} \makeatother

\newtheorem{theorem}{Theorem}[section]
\newtheorem{lemma}[theorem]{Lemma}
\newtheorem{proposition}[theorem]{Proposition}

\theoremstyle{definition}

\newtheorem{example}[theorem]{Example}

\newtheorem{remark}[theorem]{Remark}

\makeatletter
\def\ps@pprintTitle{%
  \let\@oddhead\@empty
  \let\@evenhead\@empty
  \let\@oddfoot\@empty
  \let\@evenfoot\@oddfoot
}
\makeatother

\def\ps@pprintTitle{%
  \let\@oddhead\@empty
  \let\@evenhead\@empty
  \def\@oddfoot{\reset@font\hfil\thepage\hfil}
  \let\@evenfoot\@oddfoot
}

\begin{document}

\begin{frontmatter}

\title{Bi-center conditions and bifurcation of limit cycles
in a class of $Z_2$-equivariant cubic switching
systems with two nilpotent points}

\author[mymainaddress,myaddress]
{Ting Chen}
\ead{chenting0715@126.com}
\author[mysecondaryaddress]
{Feng Li}
\ead{lf0539@126.com}
\author[mythirdaddress]
{Yun Tian}
\ead{ytian22@shnu.edu.cn}
\author[myforthaddress]
{Pei Yu\corref{mycorrespondingauthor}}
\cortext[mycorrespondingauthor]{Corresponding author}
\ead{pyu@uwo.ca}
\address[mymainaddress]{School of Statistics and Mathematics, Guangdong
University of Finance and Economics, \\
 Guangzhou, 510320, China \vspace*{0.05in}}
\address[myaddress]{College of Science, National University of Defense Technology, Changsha, 410073,  China  \vspace*{0.05in}}
\address[mysecondaryaddress]{School of Mathematics and Statistics,
Linyi University, Linyi, 276005, China \vspace*{0.05in}}
\address[mythirdaddress]{Department of Mathematics,
Shanghai Normal University, Shanghai, 200234, China \vspace*{0.05in}}
\address[myforthaddress]
{Department of Mathematics,
Western University, London, Ontario, N6A~5B7, Canada}

\begin{abstract}
In this paper, we generalize the Poincar\'e-Lyapunov method for
systems with linear
type centers to study nilpotent centers in switching polynomial systems
and use it to investigate the bi-center problem of planar
$Z_2$-equivariant cubic switching systems associated
with two symmetric nilpotent singular points.
With a properly designed perturbation, 6
explicit bi-center conditions for such
polynomial systems are derived.
Then, based on the $6$ center conditions,
by using Bogdanov-Takens bifurcation theory with
general perturbations, we prove that there exist at least $20$
small-amplitude limit cycles
around the nilpotent bi-center for a class of $Z_2$-equivariant
cubic switching systems. This is a new lower bound of cyclicity
for such cubic polynomial systems, increased from $12$ to $20$.
\end{abstract}


\end{frontmatter}

\section{Introduction}

In recent years, there have been extensive studies on the qualitative
behaviours of nonsmooth dynamical systems, associated with the vector fields
which are either discontinuous or nondifferentiable.
The main reason of the increasing interest in these studies is that
many practical systems contain nonlinearity and nonsmoothness
in their equations and display rich complex dynamical phenomena,
as observed in biology \cite{Lv},
nonlinear oscillation \cite{ADRIANA,Kunze}, dry friction of mechanical
engineering \cite{Bernardo,Leine}, power electronics \cite{Banerjee}
and so on. In fact, fundamental mathematical theory for such
nonsmooth systems was established several decades ago,
see \cite{Andronov,Minorsky}.

Switching system is one important class of nonsmooth systems,
which is divided into several regions with different
definitions of smooth vector fields. In this paper, we consider the
switching systems described in the form of
\begin{equation}\label{Eqn1}
(\dot{x},\dot{y})
=\left\{\begin{aligned}\left(f^+(x,y), \,
g^+(x,y)\right), \quad \text{if} \ \ S(x,y)\in\Sigma^+,\\
\left(f^-(x,y), \, g^-(x,y)\right), \quad \text{if} \ \ S(x,y)\in\Sigma^-,
\end{aligned}\right.
\end{equation}
where $f^{\pm}(x,y)$, $g^{\pm}(x,y)$ are analytic functions,
and the boundary $\Sigma$ between
these two regions $\Sigma^{\pm}$ defines the switching manifold.
Although the vector field could be neither continuous nor differentiable
on the boundary $\Sigma$, we recall the Filippov convention \cite{Filippov}
that the vector field may be crossing, sliding on or escaping from
the boundary. Before presenting our main results, we give an overview on
two main problems in the qualitative theory of planar differential vector
fields, namely, the center problem and the cyclicity problem.

We recall that an isolated singular point
is monodromy in a planar vector field if all
nearby orbits rotate about this point. It is well known \cite{Algaba}
that a monodromic singular point can only be a center or a focus.
The center problem is to determine whether a monodromic critical
point is a center or a focus, see \cite{Garcia}. But this problem becomes
extremely complicated in the switching systems described by \eqref{Eqn1}.
To overcome the difficulty, the authors of \cite{Gasull} presented
a useful approach for computing the Lyapunov constants of switching systems,
which can be used to derive the linear type center conditions for
an elementary singular point characterized by a pair of
purely imaginary eigenvalues. This type of centers can be determined
by vanishing all the Lyapunov constants.
In \cite{X.W1,X.W3,L.Guo2018} a complete classification with conditions
is given to determine if a singular point is a linear type center
for several classes of switching systems.
Nevertheless, for a given particular class of polynomial systems,
finding a complete solution of the center problem is still
extremely difficult.

For a given family of differential systems, the Lyapunov constants
needed to solve the center problem are also related to the so called
cyclicity problem, i.e., finding the maximum number of limit cycles around
a singular point. The investigation on bifurcation of limit cycles
in switching systems started a half century ago. In \cite{Han1},
Han and Zhang presented an interesting example to show that limit cycles
can bifurcate in switching linear systems from a pseudo-focus, i.e.,
a focus of either focus-focus, parabolic-focus or parabolic-parabolic type,
which cannot occur in smooth linear systems.
Later, Freire {\it et al.} \cite{Freire3}, Huan and Yang \cite{Huan}
independently proved that switching linear
systems can have at least $3$ limit cycles. Bastos {\it et al.}
\cite{Bastos} proved that the number of crossing limit cycles in
switching linear systems is at least $7$, which has two regions separated
by a cubic curve. Tian and Yu \cite{TIAN} constructed a special class of
switching Bautin systems to show the existence of $10$ limit cycles
in the neighborhood of a center.
Recently, by using the averaging method, Braun {\it et al.} \cite{Braun}
obtained a new lower bound with 12 limit cycles
bifurcating in quadratic switching systems. On the other hand, by the so-called pseudo-Hopf bifurcation analysis,
one more limit cycle can be generated by a proper perturbation.
However, not many works have been focused on the studies of
cubic switching systems.
Yu {\it et al.} \cite{Yu2021} proved that cubic switching systems can
have at least 18 limit cycles around two elementary foci. Chen {\it et al.}
\cite{Chen3} presented a method to study the cyclicity problem in the
neighborhood of the origin and infinity for a class of cubic switching
systems. The maximal number of limit cycles obtained so far for
cubic switching systems by perturbing an elementary center
is $26$, see \cite{Gouveia,tian2022}.

The aim of this paper is to study the center and cyclicity problems
associated with nilpotent singular points in switching systems.
Let us recall some known results
about nilpotent singular points. An isolated singular
point of a planar polynomial system is said to be a nilpotent singular
point if both eigenvalues of the Jacobian matrix of the system
evaluated at this singular point are zero
but the Jacobian matrix is not identically null.
The authors of \cite{Colak2,Giacomini,Li,Strozyna,Liu2014,Yang} developed
some computationally efficient methods for studying the nilpotent
center problem of smooth systems.
Garc\'ia \cite{Garcia} presented a technique for bounding the maximal
number of limit cycles bifurcating from a family of nilpotent singular
points in some symmetric smooth polynomial systems.
Recently, Li {\it et al.} \cite{Li2}
obtained the conditions on two nilpotent singular points to be centers
in a class of cubic smooth systems, and proved that $4$ limit cycles
exist around each of the two nilpotent singular points, with $4$ more
limit cycles bifurcating from two elementary second-order foci which are
generated by breaking the first-order nilpotent foci, leading to a total
of $12$ limit cycles, which is up to now the best result for cubic
systems with nilpotent singular points.

However, the center problem and bifurcation of
small-amplitude limit cycles for nilpotent singular
points in switching systems become more challenging than that in
the case of elementary centers. Facing the challenging, we generalize the
Poincar\'e-Lyapunov method to compute what will be called
generalized Lyapunov constants for switching nilpotent systems, and apply the method to study the bi-center problem and the cyclicity problem for a class of $Z_2$-equivariant cubic switching systems with
two symmetric nilpotent singular points. Here, we say
that a $Z_q$-equivariant differential system, whose global phase portraits
are invariant under a rotation of $2\pi/q$ ($q\in\mathrm{Z^+}$)
radians around the origin, has a bi-center at the singular
points $e_1$ and $e_2$ if both the two singular points are centers,
which will be called bi-center problem. The study of $Z_q$-equivariant
polynomial systems is closely related to the well-known Hilbert's 16th problem,
for more details see \cite{T.Chen2,Li3,Li4,J.Li2003,J.Li2010}.

Without loss of generality, the $Z_2$-equivariant cubic switching systems,
based on \eqref{Eqn1} which satisfies (with the switching
manifold $\Sigma$ defined as the $x$-axis)
$$
f^+(-x,-y)=-f^-(x,y), \quad g^+(-x,-y)=-g^-(x,y),
$$
can be written in the form of the differential equations,
\begin{equation}\label{Eqn2}
\left(\begin{array}{cc}\dot{x}\\
\dot{y}
\end{array}\right)
=\left\{\begin{aligned}&\left(\begin{aligned}
&a_{00}+a_{10}x+a_{01}y+a_{20}x^2+a_{11}xy+a_{02}y^2\\
&+a_{30}x^3+a_{21}x^2y+a_{12}xy^2+a_{03}y^3,\\
&b_{00}+b_{10}x+b_{01}y+b_{20}x^2+b_{11}xy+b_{02}y^2\\
&+b_{30}x^3+b_{21}x^2y+b_{12}xy^2+b_{03}y^3,
\end{aligned}\right), &\text{if} \ \ y>0,\\
&\left(\begin{aligned}&-a_{00}+a_{10}x+a_{01}y-a_{20}x^2-a_{11}xy
-a_{02}y^2\\
&+a_{30}x^3+a_{21}x^2y+a_{12}xy^2+a_{03}y^3,\\
&-b_{00}+b_{10}x+b_{01}y-b_{20}x^2-b_{11}xy-b_{02}y^2\\
&+b_{30}x^3+b_{21}x^2y+b_{12}xy^2+b_{03}y^3,
\end{aligned}\right), &\text{if} \ \ y<0,
\end{aligned}
\right.
\end{equation}
where the $x$-axis is the switching manifold and
all parameters are real. As is known, the type of nilpotent singular
points can generate much more rich dynamics than that of the elementary one,
such that the center of system \eqref{Eqn2} in the switching manifold
can be made up of two center-focus (a center or a focus), or one cusp
and one center-focus, or two cusps. In this paper,
we assume that system \eqref{Eqn2}
has two nilpotent singular points located at $(\pm 1, 0)$, and derive
the conditions which ensure
the nilpotent singular points $(\pm1,0)$ of \eqref{Eqn2}
to be bi-center for two cases of critical points,
namely, the third-order critical
point and the second-order critical point.
Moreover, we apply general perturbations to prove the existence
of at least $20$ small-amplitude limit cycles bifurcating from
the nilpotent bi-center $(\pm 1,0)$, two of them are obtained from
a symmetric pseudo-Hopf bifurcation by adding a suitable additional
perturbation term. This is a new lower bound (increased from
$12$ to $20$) on the number of limit cycles
bifurcating in such cubic switching systems associated with
nilpotent singular points.

The rest of the paper is organized as follows.
In Section 2, we will simplify system \eqref{Eqn2}
for the convenience of analysis and present our main results.
In Section 3, we present some formulas which are needed in Sections 4 and 5
for proving the two main theorems.
Section 4 is devoted to derive $6$ conditions for
$(\pm 1, 0)$ of system \eqref{Eqn2} to be bi-center.
Then, in Section 5 we use the $6$ center conditions to
construct perturbed systems of \eqref{Eqn2}
to show the bifurcation of $20$ limit cycles
from the nilpotent singular points $(\pm 1,0)$. Finally,
conclusion with discussion on future work is given in Section 6.

\section{Simplification of system \eqref{Eqn2} and the main results}

Assuming that system \eqref{Eqn2} has two singular points at $(\pm1,0)$
yields that
\begin{equation}\label{Eqn3}
a_{00}=-a_{20}, \quad b_{00}=-b_{20}, \quad a_{10}=-a_{30}, \quad
b_{10}=-b_{30}.
\end{equation}
The Jacobian matrices of \eqref{Eqn2} evaluated at
$(\pm 1,0)$ are given by
\begin{equation}\label{Eqn4}
\mathcal{J}^{\pm}=\left(\begin{array}{cc}
\pm2 a_{20}+2a_{30}& a_{01}\pm a_{11}+a_{21} \\
\pm2 b_{20}+2b_{30}&  b_{01}\pm b_{11}+b_{21}
\end{array}\right).
\end{equation}
It follows from $\mathcal{J}^+=\mathcal{J}^-$ that
\begin{equation}\label{Eqn5}
a_{20}=a_{11}=b_{20}=b_{11}=0.
\end{equation}

To have $(\pm 1,0)$ being isolated nilpotent singular
points, the necessary and sufficient conditions are
${\rm Tr}(\mathcal{J}^{\pm})=\det (\mathcal{J}^{\pm})=0$,
but $J^{\pm}$ is not identically zero.
It is easy to obtain that $a_{30}=b_{01}+b_{21}=0$, and
$b_{30} (a_{01}+a_{21}) = 0$ (but not $b_{30}=a_{01}+a_{21} = 0$),
which gives either $b_{30}\ne 0, \, a_{01}+a_{21} = 0$ or
$b_{30}= 0, \, a_{01}+a_{21} \ne 0$, that is
\begin{equation}\label{Eqn6}
\mathcal{J}^{\pm}=\left(\begin{array}{cc}
0& 0 \\
2b_{30}& 0
\end{array}\right)
\qquad \textrm{or} \quad
\mathcal{J}^{\pm}=\left(\begin{array}{cc}
0& a_{01}+a_{21} \\
0& 0
\end{array}\right).
\end{equation}

However, when $b_{30}=0$ and $a_{01}+a_{21}\ne0$,
system \eqref{Eqn2} becomes
\begin{equation}\label{Eqn7}
\left(\begin{array}{cc}\dot{x}\\
\dot{y}
\end{array}\right)
=\left\{\begin{aligned}&\left(\begin{aligned}
&(a_{01}+a_{02}y+a_{21}x^2+a_{12}xy+a_{03}y^2)y, \\
&(-b_{21}+b_{02}y+b_{21}x^2+b_{12}xy+b_{03}y^2)y,
\end{aligned}\right), &\text{if} \ \ y>0,\\
&\left(\begin{aligned}&(a_{01}-a_{02}y+a_{21}x^2+a_{12}xy+a_{03}y^2)y, \\
&(-b_{21}-b_{02}y+b_{21}x^2+b_{12}xy+b_{03}y^2)y,
\end{aligned}\right), &\text{if} \ \ y<0.
\end{aligned}
\right.
\end{equation}
It is easy to note that the polynomial equations in
\eqref{Eqn7} have a common factor $y$,
and so $(\pm 1,0)$ are not isolated singular points.
Thus, $b_{30}\ne 0$.

Further, to make $(\pm1,0)$ be isolated nilpotent singular points of
system \eqref{Eqn2}, we set
\begin{equation}\label{Eqn8}
\mathcal{J}^{\pm}=\left(\begin{array}{cc}
0& 0 \\
1& 0
\end{array}\right),
\end{equation}
which leads to $b_{30}=\frac{1}{2}$.
With the above results, system \eqref{Eqn7} is reduced to
\begin{equation}\label{Eqn9}
\left(\begin{array}{cc}\dot{x}\\
\dot{y}
\end{array}\right)
=\left\{\begin{aligned}&\left(\begin{aligned}
&-a_{21}y+a_{02}y^2+a_{21}x^2y+a_{12}xy^2+a_{03} y^3, \\
&-\frac{x}{2}-b_{21}y+b_{02}y^2+\frac{x^3}{2}+b_{21}x^2y+b_{12}xy^2
+b_{03}y^3, \end{aligned}\right), &\text{if} \ \ y>0,\\[1.0ex]
&\left(\begin{aligned}&-a_{21}y-a_{02}y^2+a_{21}x^2y+a_{12}xy^2
+a_{03} y^3, \\
&-\frac{x}{2}-b_{21}y-b_{02}y^2+\frac{x^3}{2}+b_{21}x^2y+b_{12}xy^2
+b_{03}y^3, \end{aligned}\right), &\text{if} \ \ y<0.
\end{aligned}
\right.
\end{equation}
For the sake of convenience, we call the system
in \eqref{Eqn9} for $y >0$ ``the  upper system'' and that
for $y < 0$ ``the lower system''.

Now, we need to discuss the multiplicity of the nilpotent singular
points $(\pm1,0)$ for the upper and the lower systems in \eqref{Eqn9}.
In fact, by the symmetry of system \eqref{Eqn9}, we only need to analyze
the singular point $(1,0)$. Introducing the transformation
$x\rightarrow x+1$ into system \eqref{Eqn9} results in
\begin{equation}\label{Eqn10}
\left(\begin{array}{cc}\dot{x}\\
\dot{y}
\end{array}\right)
=\left\{\begin{aligned}&\left(\begin{aligned}
&2a_{21}xy+(a_{02}+a_{12})y^2+a_{21}x^2y+a_{12}xy^2\\
&+a_{03}y^3 \stackrel{\triangle}{=}\Psi^+(x,y),\\
&x+\frac{3x^2}{2}+2b_{21}xy+(b_{02}+b_{12})y^2+\frac{x^3}{2}+b_{21}x^2y\\
&+b_{12}xy^2+b_{03}y^3 \stackrel{\triangle}{=} \Phi^+(x,y)+x,
\end{aligned}\right), &\text{if} \ \ y>0,\\[1.0ex]
&\left(\begin{aligned}&2a_{21}xy+(a_{12}-a_{02})y^2+a_{21}x^2y+a_{12}xy^2\\
&+a_{03}y^3 \stackrel{\triangle}{=} \Psi^-(x,y),\\
&x+\frac{3x^2}{2}+2b_{21}xy+(b_{12}-b_{02})y^2+\frac{x^3}{2}+b_{21}x^2y\\
&+b_{12}xy^2 +b_{03}y^3 \stackrel{\triangle}{=} \Phi^-(x,y)+x,
\end{aligned}\right), &\text{if} \ \ y<0,
\end{aligned}
\right.
\end{equation}
and thus the singular point $(1,0)$ of system \eqref{Eqn9} is
shifted to the origin of system \eqref{Eqn10}. Assume that
$$
x^{\pm}=\mathfrak{f}^{\pm}(y)=\sum^\infty_{k=2}c_k^{\pm}y^k
$$
are the unique solutions of the implicit function equations
$\Phi^{\pm}(x,y)+x=0$. Further, we denote that
\begin{equation}\label{Eqn11}
\begin{aligned}
&F^{\pm}(y)=\Psi^{\pm}(\mathfrak{f}^{\pm}(y),y)
=\sum^\infty_{k=2}\mathfrak{f}_k^{\pm}y^k,\\
&G^{\pm}(y)=\bigg[\frac{\partial \Psi^{\pm}}{\partial x}
+\frac{\partial \Phi^{\pm}}{\partial y}\bigg]_{(\mathfrak{f}^{\pm}(y),y)}
=\sum^\infty_{k=1}\mathfrak{g}_k^{\pm}y^k,
\end{aligned}
\end{equation}
where
\begin{equation}\label{Eqn12}
\begin{aligned}
&\mathfrak{g}_1^{\pm}=2(a_{21}\pm b_{02}+b_{12}),\\
&\mathfrak{f}_2^{\pm}=\pm a_{02}+a_{12},\\
&\mathfrak{f}_3^{\pm}=a_{03}\mp 2a_{21}b_{02}-2a_{21}b_{12}.
\end{aligned}
\end{equation}

For planar smooth systems, if $k$ is the smallest integer
satisfying $\mathfrak{f}_k\neq0$,
then the multiplicity of the nilpotent singular point
is exactly $k$, for more details see \cite{Li2}. Thus, we can use
Theorem 3.5 in \cite{Dumortier2006} to determine the
type of the nilpotent singular point
$(0,0)$ of \eqref{Eqn10}
for the smooth polynomial equations either in the upper system or
in the lower system.

More precisely, when $\mathfrak{g}_n=0$ and $\mathfrak{f}_m\neq0$,
we have that
\begin{equation}\label{Eqn13}
\left\{\begin{aligned}&m=2k+1,~\left\{\begin{aligned} &\mathfrak{f}_m<0,\
\textrm{$(0,0)$ is a center or focus},\\
&\mathfrak{f}_m>0,\ \textrm{$(0,0)$ is a saddle},\end{aligned}\right.\\
&m=2k,\ \textrm{$(0,0)$ is a cusp}.\\
\end{aligned}
\right.
\end{equation}
When $\mathfrak{g}_n\neq0$, $\mathfrak{f}_m\neq0$ and
$\Delta=4(n+1)\mathfrak{f}_{m}+\mathfrak{g}_n^2$, we have that
\begin{equation}\label{Eqn14}
\left\{\begin{aligned}
&m=2k+1,\
\left\{\begin{aligned} &\mathfrak{f}_m>0, \ \textrm{$(0,0)$ is a saddle},\\
&\mathfrak{f}_m<0, \
\left\{\begin{aligned}&k>n,\ \textrm{or} \ k=n \ \textrm{and} \ \Delta\geq0,
\ \left\{\begin{aligned} &n\ \textrm{odd, $(0,0)$ is a H-E},\\
&n \ \textrm{even, $(0,0)$ is a node},\end{aligned}\right.\\
 &k<n,\ \textrm{or} \
k=n \ \textrm{and} \ \Delta<0, \ \textrm{$(0,0)$ is a center or focus},\\
\end{aligned}\right.\end{aligned}\right.\\
&m=2k,\ \left\{\begin{aligned}
&k>n,\ \textrm{$(0,0)$ is a saddle-node},\\
&k\leq n, \
\textrm{$(0,0)$ is a cusp},
\end{aligned}\right.
\end{aligned}\right.
\end{equation}
where H-E denotes the local phase portrait with one singular
point consisting of one hyperbolic sector and one elliptic sector.

\begin{proposition}\label{Prop2.1}
The multiplicity of the nilpotent singular point $(1,0)$ of the upper
system $($or $(-1,0)$ of the lower system$)$ of \eqref{Eqn9} is at most 6.
\end{proposition}

\begin{proof}
Using \eqref{Eqn12} with $\mathfrak{f}_2^+=\mathfrak{f}_3^+=0$, we have
\begin{equation*}
a_{02}=-a_{12}, \quad a_{03}=2a_{21}(b_{02}+b_{12}).
\end{equation*}
First, assume that $a_{21} \!=\!0$. Then, we obtain
$\mathfrak{f}_4^+ \!=\! -a_{12}(b_{02} + b_{12})$.
If $a_{12} \!=\! 0$, we have $\Psi^+(x,y)=0$, and so $(1,0)$ is not an isolated
singular point, implying that $\mathfrak{f}_4^+ \ne 0$.
If $b_{02}=-b_{12}$, we have $\mathfrak{f}_5^+=-a_{12}b_{03}$.
Letting $b_{03}=0$ yields $\Psi^+(x,y)=a_{12}xy^2$, and then $\Psi^+(x,y)$
and $x+\Phi^+(x,y)$ have a common factor $x$. Hence,
$\mathfrak{f}_5^+ \ne 0$ if $(1,0)$
is an isolated singular point.

Next, assume that $a_{21}\neq0$. Then, we have $\mathfrak{f}_4^+=(b_{02}
+b_{12})(4a_{21}b_{21}-a_{12})-2a_{21}b_{03}$. Setting
$$
b_{03}=\frac{(b_{02}+b_{12})(4a_{21}b_{21}-a_{12})}{2 a_{21}},
$$
we have $\mathfrak{f}_4^+=0$ and
$$
\mathfrak{f}_5^+ =-\frac{(b_{02}+b_{12})(-a_{12}^2+4a_{21}^2b_{02}
+4a_{12}a_{21}b_{21})}{2a_{21}}.
$$
If $b_{02}=-b_{12}$, we have $b_{03}=0$.
Further, we obtain that $\Psi^+(x,y)$
and $\Phi^+(x,y)+x$ have a common factor $x$, leading to
$\mathfrak{f}_5^+ \ne 0$. Otherwise,
we assume that $b_{02}=\frac{a_{12}^2-4a_{12}a_{21}b_{21}}{4a_{21}^2}$,
under which $\mathfrak{f}_5^+ = 0$. Then, we have
$$
\mathfrak{f}_6^+=\frac{a_{12}(a_{12}^2+4a_{21}^2b_{12}
-4a_{12}a_{21}b_{21})^2}{32 a_{21}^4}.
$$
If $a_{12}=0$, we obtain that $\Psi^+(x,y)$ and $\Phi^+(x,y)+x$
have a common factor $2x+x^2+2b_{12}y^2$.
If $b_{12}=\frac{-a_{12}^2+4a_{12}a_{21}b_{21}}{4a_{21}^2}$,
we have $\Psi^+(x,y)=xy[a_{21}(2+x)+a_{12} y]$.
Then $\Psi^+(x,y)$ and $\Phi^+(x,y)+x$ have a common factor $x$.
Hence, $(1,0)$ is not an isolated singular point,
implying that $\mathfrak{f}_6^+ \ne 0$.
\end{proof}

Next, we derive the conditions for the nilpotent singular points ($\pm1$,0)
of the $Z_2$-equivariant cubic switching system \eqref{Eqn9}
to be bi-center, yielding the following result.

\begin{theorem}\label{Thm2.2}
The nilpotent singular points {\rm($\pm1$,0)} of
system \eqref{Eqn9} become bi-center if one of
the following conditions holds:
\begin{equation}\label{Eqn15}
\begin{aligned}
\mathrm{I}\!:& \ \, a_{02}+a_{12}=a_{21}+b_{12}=a_{12}+3b_{03}=b_{02}=b_{21}=0,
\ a_{12}>0,\ \, a_{03}+2a_{21}^2<0;\\
\mathrm{II}\!:& \ \, a_{02}=a_{12}=b_{02}=b_{03}=b_{21}=0,\
2 a_{03} + (b_{12}-a_{21})^2 < 0; \\
\mathrm{III}\!:& \ \, a_{02}=a_{12}=b_{02}=a_{21}+b_{12}=3b_{03}+2a_{21}b_{21}
\\
&\hspace*{0.29in}
= 9 (a_{03}+2 a_{21}^2)^2 + 8 a_{21}b_{21}^2 (3 a_{03} + 2 a_{21}^2)=0,
\ \, a_{03}+2a_{21}^2<0;\\
\mathrm{IV}\!:& \ \,  a_{02}=a_{12}=b_{02}=8a_{21}+3b_{21}^2
=16a_{03}-3b_{21}^2(4b_{12}+b_{21}^2)\\
& \hspace*{0.29in} =8b_{03}-b_{21}(8b_{12}-b_{21}^2)=0,\ \,
-\tfrac{9+4 \sqrt{3}}{8}\, b_{21}^2 < b_{12}
< -\tfrac{9-4 \sqrt{3}}{8}\, b_{21}^2; \\
\mathrm{V}\!:& \ \, a_{12}=b_{02}=b_{03}=b_{21}=0,\ a_{02}<0;\\
\mathrm{VI}\!:& \ \, a_{21}+ b_{12}=a_{12}+ 3b_{03}=b_{02}=b_{21}=0,
\ \, \text{either}\ \, a_{02} + |a_{12}|<0 \ \, \text{or}\ \, a_{02}=a_{12}<0.
\end{aligned}
\end{equation}
\end{theorem}

Moreover, to find the maximal number of limit cycles around the
nilpotent singular points ($\pm1$,0) of system \eqref{Eqn9},
we construct perturbed $Z_2$-equivariant switching systems
using the $6$ bi-center conditions and prove that
system \eqref{Eqn9} has at least $20$ limit cycles bifurcating
from the nilpotent singular points ($\pm1$,0).
More precisely, we have the following theorem.

\begin{theorem}\label{Thm2.3}
Under each of the $6$ nilpotent bi-center conditions in Theorem
{\rm \ref{Thm2.2}}, the $Z_2$-equivariant cubic switching system
\eqref{Eqn9} has at least $18$ limit cycles bifurcating from
the two symmetric nilpotent singular points $(\pm1,0)$
by using all $\varepsilon^k$-order cubic perturbation,
and at least $20$ limit cycles by, in addition, applying
a symmetric pseudo-Hopf bifurcation.
\end{theorem}

\section{The generalized Poincar\'e-Lyapunov method}

In this section, for the convenience of reading, we first briefly
describe the Poincar\'e-Lyapunov method for determining
the linear type center in switching polynomial systems divided
by the $x$-axis, and then generate this method to study switching systems
associated with isolated nilpotent singular points.
As described in \cite{X.W3}, we consider the following switching system,
\begin{equation}\label{Eqn16}
(\dot{x},\, \dot{y})=\left\{\begin{aligned}
\Big(\delta x-\lambda^+y+\sum\limits_{i+j=2}^{n}a_{ij}^+x^iy^j, \
\lambda^+x+\delta y+\sum\limits_{i+j=2}^{n}b_{ij}^+x^iy^j\Big),
\ \ \text{if} \ \
y>0,\\[0.0ex]
\Big(\delta x-\lambda^-y+\sum\limits_{i+j=2}^{n}a_{ij}^-x^iy^j, \
\lambda^-x+\delta y+\sum\limits_{i+j=2}^{n}b_{ij}^-x^iy^j\Big),
\ \ \text{if} \ \
y<0,
\end{aligned}\right.
\end{equation}
where $\delta, a_{ij}^{\pm}, b_{ij}^{\pm}\in R$, $\lambda^{\pm}>0$.
Let $\Lambda^{\pm}=(\delta,a_{ij}^{\pm},b_{ij}^{\pm},\lambda^{\pm})$
represent two parameter vectors.
The origin is a common singular point in both upper and lower systems
of \eqref{Eqn16}. With the polar coordinate transformation:
$x=r \cos\theta$ and $y=r \sin\theta$, system \eqref{Eqn16}
can be rewritten as
\begin{equation}\label{Eqn17}
\frac{{\rm d}r}{{\rm d}\theta}=\left\{\begin{aligned}
&\frac{\delta r+\sum\limits_{i+j=3}^{n+1}c_{ij}^+(a_{ij}^+,b_{ij}^+)
\cos\theta^i\sin\theta^jr^{i+j-1}}{\lambda^+
+\sum\limits_{i+j=3}^{n+1}d_{ij}^+(a_{ij}^+,b_{ij}^+)\cos\theta^i
\sin\theta^jr^{i+j-2}}, \quad \text{if} \ \ \theta\in[0,\pi],\\[1.0ex]
&\frac{\delta r+\sum\limits_{i+j=3}^{n+1}c_{ij}^-(a_{ij}^-,b_{ij}^-)
\cos\theta^i\sin\theta^jr^{i+j-1}}{\lambda^-
+\sum\limits_{i+j=3}^{n+1}d_{ij}^-(a_{ij}^-,b_{ij}^-)\cos\theta^i
\sin\theta^jr^{i+j-2}}, \quad \text{if} \ \ \theta\in[\pi,2\pi],
\end{aligned} \right.
\end{equation}
where $c_{ij}^\pm(a_{ij}^\pm,b_{ij}^\pm)$ and
$d_{ij}^\pm(a_{ij}^\pm,b_{ij}^\pm)$ are polynomials in the
parameters $a_{ij}^\pm$ and $b_{ij}^\pm$.
We suppose that the solutions for the upper and the lower systems of
\eqref{Eqn17} are respectively obtained as
$r^+(\varrho,\Lambda^+,\theta)=\sum\limits_{k\geq1}
v_{k}^+(\Lambda^+,\theta)\varrho^k$ and
$r^-(\varrho,\Lambda^-,\theta)=\sum\limits_{k\geq1}
v_{k}^-(\Lambda^+,\theta)\varrho^k$ satisfying the initial
condition $r^+(\varrho,\Lambda^+,0)=r^-(\varrho,\Lambda^-,\pi)=\varrho$.
Then, we denote by
$$
\Upsilon^+(\varrho)=r^+(\varrho,\Lambda^+,\pi)
=e^{\pi \frac{\delta}{\lambda^+}}\varrho
+\sum\limits_{k\geq2}v_{k}^+\varrho^k
$$
and
$$
\Upsilon^-(\varrho)=r^-(\varrho,\Lambda^-,2\pi)
=e^{\pi \frac{\delta}{\lambda^-}}\varrho+\sum\limits_{k\geq2}v_{k}^-\varrho^k,
$$
the upper half-return map $\Upsilon^+(\varrho)$ and the lower half-return
map $\Upsilon^-(\varrho)$, respectively, where $v_k^{\pm}$'s
are the coefficients
in Taylor expansions. Following the procedure in \cite{Gasull},
we obtain the displacement function,
\begin{equation}\label{Eqn18}
\begin{aligned}
d(\varrho)=\Upsilon^+(\varrho)-(\Upsilon^-)^{-1}(\varrho)
=\sum\limits_{k\geq1}V_k\varrho^k,
\end{aligned}
\end{equation}
where $V_{k}$ is called the $k$th-order Lyapunov constant of system
\eqref{Eqn16}. The origin of system \eqref{Eqn16} is a center if and only
if all the Lyapunov constants in \eqref{Eqn18} vanish,
i.e., $d(\varrho)\equiv0$ for $0<\varrho\ll1$.
If there exists $\chi_{*}\in(\Lambda^+,\Lambda^-)$ such that $
V_{1}(\chi_{*})=V_{2}(\chi_{*})=\cdots=V_{k}(\chi_{*})=0$
and $V_{k+1}(\chi_{*})\neq0$, then any perturbations to system \eqref{Eqn16}
can yield at most $k$ limit cycles bifurcating from the origin.
Based on Lemma 4 in \cite{TIAN}, we give the sufficient conditions
for proving the existence
of limit cycles bifurcating from the origin of system \eqref{Eqn16}.

\begin{lemma}[\cite{TIAN}]\label{Lem3.1}
If there exists a  critical point
$\chi_*=(a_{1c},a_{2c},\cdots,a_{kc})$ such that
$V_{i_1}(\chi_*)=V_{i_2}(\chi_*)=\cdots=V_{i_k}(\chi_*)=0$,
$V_{i_{k+1}}(\chi_*)\neq0$,  with $1= i_1<i_2<\cdots<i_k$, and
\begin{equation}\label{Eqn19}
\begin{array}{l}
{\rm{det}}\bigg[\displaystyle\frac{\partial(V_{i_1},V_{i_2},\cdots,
V_{i_k})}{\partial(a_{1c},a_{2c},\cdots,a_{kc})} (\chi_*)\bigg]\neq0,
\end{array}
\end{equation}
then appropriate small perturbations about $\chi=\chi_*$
lead to that system \eqref{Eqn16} has exactly $k$ limit cycles bifurcating
from the origin.
\end{lemma}

\begin{figure}[!htb]
\vspace*{0.20in}
\hspace*{-0.70in}
\centering
\begin{overpic}[width=1.3\textwidth,height=0.17\textheight]{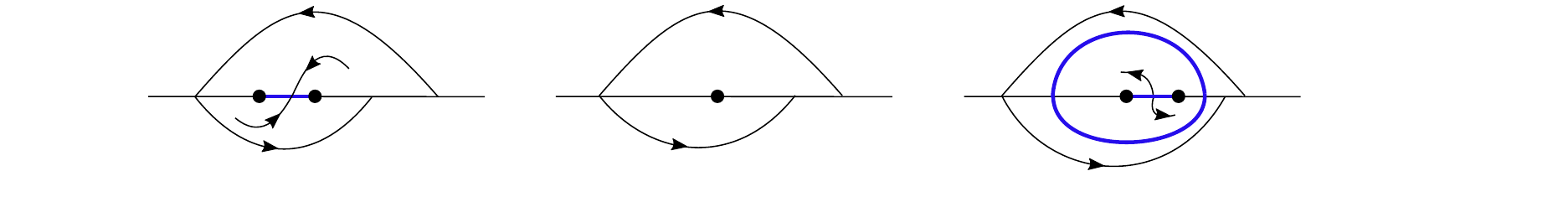}
\put(18,2){$b<0$}
\put(44,2){$b=0$}
\put(69,2){$b>0$}
\end{overpic}

\vspace*{-0.15in}
\caption{Illustration of the idea of pseudo-Hopf bifurcation.}\label{Fig1}
\vspace*{0.10in}
\end{figure}

It is worth mentioning that a switching polynomial system can have
one more small-amplitude limit cycle when the sliding segment changes
its stability by adding constant terms, which is called pesudo-Hopf
bifurcation, see \cite{Cruz,Filippov}. With $\delta>0$ and
$b\in \mathbb{R}$, we assume that $f^{\pm}(x,y)$ and $g^+(x,y)$ are
given in \eqref{Eqn16}, but $g^-(x,y)$ has the following form
\begin{equation}\label{Eqn20}
g^-(x,y)
=b+\lambda^-x+\delta y+\sum\limits_{i+j=2}^{n}b_{ij}^-x^iy^j.
\end{equation}
Then, the upper system of \eqref{Eqn16} still has a singular point
at the origin, while the lower system of \eqref{Eqn16} has a singular
point near the origin. System \eqref{Eqn16} would have a sliding
segment on the switching manifold $y=0$ (the $x$-axis).
For small enough $b$,
the new switching system exhibits a pseudo-Hopf bifurcation at
$b=0$, as shown in Figure \ref{Fig1}, which can produce one more
small-amplitude limit cycle, see more details in \cite{Castillo}.

In the above discussions, we have briefly described the Poincar\'e-Lyapunov
method and the pesudo-Hopf bifurcation for switching polynomial systems
with linear type centers. Now, we consider the switching polynomial
system with a nilpotent critical point at the origin,
\begin{equation}\label{Eqn21}
(\dot{x},\, \dot{y})=\left\{\begin{aligned}\left(F_1^+(x,y),\
x+F_2^+(x,y)\right), \ \ \text{if} \ \ y>0,\\
\left(F_1^-(x,y),\ x+F_2^-(x,y)\right), \ \ \text{if} \ \ y<0,
\end{aligned}\right.
\end{equation}
where $F_1^{\pm}(x,y)$ and $F_2^{\pm}(x,y)$ are real analytic functions
without constant and linear terms. We will show that the
Poincar\'e-Lyapunov method can be generalized to
analyze center problem and bifurcation of limit cycles
for the nilpotent singular point.

We give the following example to illustrate our idea, as the normal form
of the nilpotent differential system can be simplified into
the system (for example, see \cite{Han2}),
\begin{equation}\label{Eqn22}
\begin{aligned}
\dot{x}=&\ y,\\
\dot{y}=&\ c_ix^i[1+g(x)]+d_ix^{i-1}y[1+h(x)]+y^2q(x,y),
\end{aligned}
\end{equation}
where $i\geq2$, $g(x)$, $h(x)$ and $q(x,y)$ are analytic functions satisfying
$g(0)=h(0)=0$. For the sake of simplicity, we consider a codimension-2
Bogdanov-Takens (B-T) bifurcation of the $Z_2$ cubic normal form
\eqref{Eqn22}. Then, we obtain the following normal form with unfolding:
\begin{equation}\label{Eqn23}
\begin{aligned}
\dot{x}=&\ y,\\
\dot{y}=&\ \alpha x+\beta y+c_3x^3+d_3x^2y,
\end{aligned}
\end{equation}
where the term $\alpha x+\beta y$ is called unfolding with sufficiently
small parameters $\alpha$ and $\beta$. Note that the origin of system
\eqref{Eqn23} is a linear type singular point when $\alpha<0$ and $\beta=0$.
When $\beta=0$,
$\alpha\rightarrow0^{-}$, the linear type origin becomes
a nilpotent singular point.

We remark that the linear type origin of system \eqref{Eqn23} is a center
if and only if $d_3=0$, and that the nilpotent monodromic origin
(when $c_3<0$) of system \eqref{Eqn23} is also a center when $d_3=0$.
Hence, the main idea of our method is to perturb system \eqref{Eqn21}
and to establish a relation between the unperturbed system and the
perturbed system based on the B-T bifurcation theory, which will
generate perturbed Lyapunov constants.

To achieve this, we consider the following perturbed system of \eqref{Eqn21},
\begin{equation}\label{Eqn24}
(\dot{x},\, \dot{y})=\left\{\begin{aligned}
\left(-\varepsilon^2 y+F_1^+(x,y)+\varepsilon G_1^+(x,y,\varepsilon), \
x+F_2^+(x,y)+\varepsilon G_2^+(x,y,\varepsilon)\right),\ \ \text{if} \ \ y>0,\\
\left(-\varepsilon^2 y+F_1^-(x,y)+\varepsilon G_1^-(x,y,\varepsilon), \
x+F_2^-(x,y)+\varepsilon G_2^-(x,y,\varepsilon)\right),\ \ \text{if} \ \ y<0,
\end{aligned}\right.
\end{equation}
where $-\,\varepsilon^2 y$ is called unfolding with sufficiently
small $|\varepsilon| \! \ll \! 1$, $G_1^{\pm}(x,y,\varepsilon)$ and
$G_2^{\pm}(x,y,\varepsilon)$ are real polynomials.

It should be noted that the
$\varepsilon$-perturbation terms in \eqref{Eqn24} are applied to
the whole system, not restricted to the local behavior.
We use the $\varepsilon^2$ rather than $\varepsilon$ in the
linear perturbation term to avoid the $\varepsilon^{-i}$ ($i>0$)
and $\sqrt{|\varepsilon|}$ terms generalted in the later
transformed systems. Based on B-T bifurcation theory and the
relation established above for the two systems \eqref{Eqn21}
and \eqref{Eqn24}, we know that a nilpotent center can be
$\varepsilon$-approximated by a linear type center.
More detailed discussions on this subject can be found in \cite{Chen4}.

For this type system of \eqref{Eqn24},
we can apply the Poincar\'e-Lyapunov method
and compute the generalized Lyapunov constants $V_k(\varepsilon)$.
As a matter of fact, we have the
displacement function of system \eqref{Eqn24}, given by
\begin{equation}\label{Eqn25}
d(\varrho)=\sum_{k\geq1} V_k(\varepsilon)\varrho^k,
\quad \textrm{where} \quad
V_i(\varepsilon)=\sum_{j=1}^\infty \varepsilon^j V_{ij}, \quad i=1,2,\cdots,
\end{equation}
in which $V_{ij}$ denotes the $i$th $\varepsilon^j$-order Lyapunov constant.
We can determine the center conditions for system \eqref{Eqn24}
by vanishing the $\varepsilon$ terms in these generalized Lyapunov constants,
thus leading to a set of algebraic conditions
which characterize the existence of a nilpotent center of
system \eqref{Eqn21}.

Further, we can use the above generalized Poincar\'e-Lyapunov method
to study the bifurcation of limit cycles
from the nilpotent center of the switching system \eqref{Eqn21}.
By B-T bifurcation theory, we add a linear perturbation term
$-\varepsilon^2 y$ to system \eqref{Eqn21} such that the origin
of this switching system becomes a linear-type center.
Then, we compute the generalized Lyapunov constants of the perturbed
system with additional all
$\varepsilon^k$-order perturbation terms.
Following the procedure in \cite{ELT2023}, we can obtain the bifurcation
of limit cycles from the nilpotent center as many as possible.
Finally, we obtain one more limit cycle in the switching system
by considering the pseudo-Hopf bifurcation.

\section{The proof of Theorem \ref{Thm2.2}}

Since the multiplicity of monodromic critical point for smooth
nilpotent systems is $2k+1\geq3$ (see \cite{Dumortier2006}),
we know that the smallest multiplicity for a singular point is $3$
if it is a nilpotent focus or center. For the sake of convenience,
we call the singular point with multiplicity $3$
the 3rd-order critical point. However, the multiplicity of the monodromic
singular point in switching systems can be $2$, i.e., a 2nd-order
critical point. Hence, we consider four cases in the following
two subsections.

\subsection{The $3$rd-order critical point $(1,0)$ of the
upper system in \eqref{Eqn9}}

By using the conditions given in \eqref{Eqn13} and \eqref{Eqn14},
we obtain that the singular point $(1,0)$ of the upper system
of \eqref{Eqn9} is a 3rd-order nilpotent focus or center if and only if
\begin{equation}\label{Eqn26}
\mathfrak{f}_2^+=0,\quad \mathfrak{f}_3^+<0,\quad
\Delta^+=(\mathfrak{g}_1^+)^2+8\mathfrak{f}_3^+<0,
\end{equation}
namely,
\begin{equation}\label{Eqn27}
\begin{aligned}
&a_{02}=-a_{12},\quad 2 a_{03} + (b_{02}+b_{12}-a_{21})^2<0.
\end{aligned}
\end{equation}
Then, we obtain
\begin{equation*}
\begin{aligned}
&\mathfrak{f}_2^-=2a_{12},\quad \mathfrak{f}_3^-
=a_{03}+2a_{21}(b_{02}-b_{12}),\quad
\Delta^-=8a_{03}+4(a_{21}+b_{02}-b_{12})^2.
\end{aligned}
\end{equation*}
When $\mathfrak{f}_2^-\neq0$, i.e., $a_{12}\neq0$, the singular
point $(1,0)$ of the lower system in \eqref{Eqn9} is a cusp.
On the other hand, the singular point $(1,0)$ of the lower system
in \eqref{Eqn9} is a center or a focus with multiplicity $3$
if $\mathfrak{f}_2^- \!= 0$ (i.e., $a_{12}=0$), $\mathfrak{f}_3^- \!<\! 0$
and $\Delta^- \!<\! 0$. In order to have a monodromic singular point at $(1,0)$
of \eqref{Eqn9}, it requires that there do not exist seperatrices
which connect this singular point in the upper system or the lower system.
With the aid of Maple,
we present the following example to demonstrate a
phase portrait of system \eqref{Eqn9}, indicating
that the lower system of \eqref{Eqn9} has two seperatrices
connecting $(1,0)$ in the lower half Poincar\'e disc when $a_{12}<0$.
Thus, we only need to consider $a_{12}\geq0$.

\begin{figure}[!htb]
\centering
\includegraphics[width=0.28\textwidth,height=0.2\textheight]{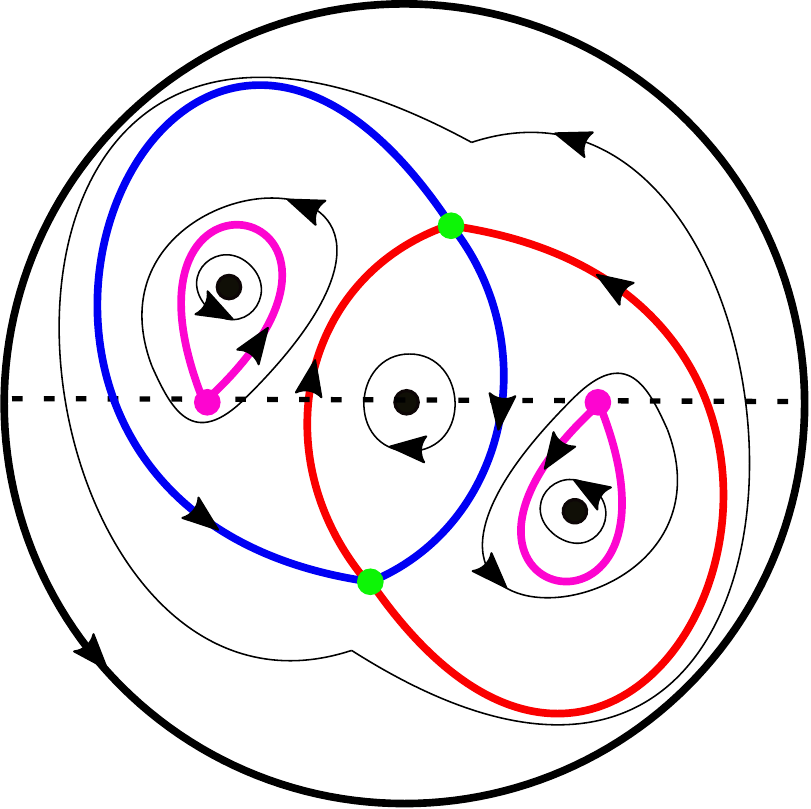}
\vspace*{0.10in}
\caption{The phase portrait of system \eqref{Eqn9} with $a_{02} \!=\!
b_{12} \!=\! 1$, $a_{21} \!=\! a_{12} \!=\! -1$,$a_{03} \!=\! -4$,
$b_{02} \!=\! b_{21}=0$ and $b_{03} \!=\! \frac{1}{3}$,
showing that the singular points $(\pm1,0)$ are two cusps.}
\label{Fig2}
\end{figure}

\begin{example}\label{Exam4.1}
The phase portrait of system \eqref{Eqn9} with
$a_{02}=b_{12}=1$, $a_{21}=a_{12}=-1$, $a_{03}=-4$, $b_{02}=b_{21}=0$
and $b_{03}=\frac{1}{3}$, as depicted in Figure~\ref{Fig2},
shows that the singular points $(\pm1,0)$ are two cusps.
\end{example}

To discuss the bi-center conditions for $(\pm1,0)$ of system \eqref{Eqn9},
we assume that $a_{02}=-a_{12}$, as given in \eqref{Eqn27}.
Then, we show how to apply the Poincar\'e-Lyapunov method to
derive the bi-center conditions for $(\pm1,0)$ of system \eqref{Eqn9}.
By the symmetry of system \eqref{Eqn9}, we only need to study the
center conditions at the singular point $(1,0)$.

With the transformation $x\rightarrow x+1$ applied and perturbations added
to system \eqref{Eqn9} {\color{red}with \eqref{Eqn27}}, we obtain the following perturbed system, {\color{blue} can check from \eqref{Eqn10}}
\begin{equation}\label{Eqn28}
\left(\begin{array}{cc}\dot{x}\\
\dot{y}
\end{array}\right)
=\left\{\begin{aligned}&\left(\begin{aligned}
&-\varepsilon^2 y+2a_{21}xy+a_{21}x^2y+a_{12}xy^2+a_{03} y^3\\
&+\sum_{k=1}\varepsilon^k G_{1k}^+(x,y),\\
&x+\frac{3x^2}{2}+2b_{21}xy+
(b_{02}+b_{12})y^2+\frac{x^3}{2}\\
&+b_{21}x^2y+b_{12}xy^2+b_{03}y^3+\sum_{k=1}\varepsilon^k G_{2k}^+(x,y),
\end{aligned}\right), &\text{if} \ \ y>0,\\[1.0ex]
&\left(\begin{aligned}&-\varepsilon^2 y+2a_{21}xy+2a_{12}y^2
+a_{21}x^2y+a_{12}xy^2 +a_{03}y^3
\\
&+\sum_{k=1}\varepsilon^k G_{1k}^-(x,y),\\
&x+\frac{3x^2}{2}+2b_{21}xy
+(b_{12}-b_{02})y^2+\frac{x^3}{2}\\
&+b_{21}x^2y+b_{12}xy^2+b_{03}y^3+\sum_{k=1}\varepsilon^k G_{2k}^-(x,y),
\end{aligned}\right), &\text{if} \ \ y<0,
\end{aligned}
\right.
\end{equation}
where
\begin{equation*}
G_{1k}^{\pm}(x,y)=\sum_{i+j=2}^{3}p_{ijk}^{\pm}\, x^iy^j,\quad
G_{2k}^{\pm}(x,y)=\sum_{i+j=2}^{3}q_{ijk}^{\pm}\, x^iy^j,
\end{equation*}
in which $p_{ijk}^{\pm}$, $q_{ijk}^{\pm}$ are real parameters.

Note that system \eqref{Eqn28} should be invariant under the
transformation $ x \rightarrow x-1$
(since $(1,0)$ and $(-1,0)$ are symmetric
singular points of \eqref{Eqn9}), which yields
\begin{equation}\label{Eqn29}
\begin{array}{llll}
l_{20k}^{\pm}=l_{30k}^{\pm}=0, & l_{11k}^-=l_{11k}^+
=2l_{21k}^+, & l_{02k}^-=-l_{02k}^++2l_{12k}^+,\\[1.0ex]
l_{21k}^-=l_{21k}^+, & l_{12k}^-=l_{12k}^+, & l_{03k}^-=l_{03k}^+,
\end{array}
\end{equation}
where $l_{ijk}^{\pm}$ represent $p_{ijk}^{\pm}$ or $q_{ijk}^{\pm}$.

Since system \eqref{Eqn28} has a large number of perturbation
parameters, which makes the calculation of the generalized Lyapuov
constants very difficult, we need to reduce the number of parameters.
Without loss generality, we assume that
$q_{21k}^+=0$ because they are redundant for proving the necessity
of the conditions I-VI.
Further, we consider only the $\varepsilon^2$-order perturbed terms,
i.e., setting $l_{ij1}^+=l_{ijk}^+=0$ ($k>2$).
The reason for why only choosing $\varepsilon^2$-order
perturbed terms will be given in Section 5.

For convenience, we let $l_{ij2}^+=l_{ij}$.
Further, introducing the transformation
$(x,y,t)\rightarrow (\varepsilon^3 x, \varepsilon^2 y, \frac{t}{\varepsilon})$
into system \eqref{Eqn28}, we obtain
\begin{equation}\label{Eqn30}
\left( \begin{array}{cc}\dot{x}\\
\dot{y}
\end{array} \right)
=\left\{\begin{aligned}&\left(\begin{aligned}&-y
+2(a_{21}+p_{21}\varepsilon^2)\varepsilon xy
+p_{02}\varepsilon^2 y^2+(a_{21}+p_{21}\varepsilon^2 )\varepsilon^4 x^2y\\
&+(a_{12}+p_{12}\varepsilon^2)\varepsilon^3 xy^2
+(a_{03}+p_{03}\varepsilon^2)\varepsilon^2 y^3,\\[1.0ex]
&x+\tfrac{3}{2}\varepsilon^3 x^2+2b_{21}\varepsilon^2xy
+(b_{02}+b_{12}+q_{02}\varepsilon^2)\varepsilon y^2
+\tfrac{1}{2}\varepsilon^6 x^3\\
&+b_{21}\varepsilon^5x^2y +(b_{12}+q_{12}\varepsilon^2)\varepsilon^4 xy^2
+(b_{03}+q_{03}\varepsilon^2)\varepsilon^3 y^3,
\end{aligned}\right), &\text{if} \ \ y>0,\\[1.5ex]
&\left(\begin{aligned}&-y+2(a_{21}+p_{21}\varepsilon^2)\varepsilon xy
+(2a_{12}-p_{02}\varepsilon^2+2p_{12}\varepsilon^2)y^2\\
&+(a_{21}+p_{21}\varepsilon^2)\varepsilon^4 x^2y
+(a_{12}+p_{12}\varepsilon^2)\varepsilon^3 xy^2\\
&+(a_{03}+p_{03}\varepsilon^2)\varepsilon^2 y^3,\\[1.0ex]
&x+\tfrac{3}{2}\varepsilon^3 x^2
+2b_{21}\varepsilon^2xy
+(b_{12}-b_{02}-q_{02}\varepsilon^2+2q_{12}\varepsilon^2)\varepsilon y^2 \\
&+\tfrac{1}{2}\varepsilon^6 x^3+b_{21}\varepsilon^5x^2y
+(b_{12}+q_{12}\varepsilon^2)\varepsilon^4 xy^2\\
&+(b_{03}+q_{03}\varepsilon^2)\varepsilon^3 y^3,
\end{aligned}\right), &\text{if} \ \ y<0.
\end{aligned}
\right.
\end{equation}

To give a clear view of the proof, we first present a table below to
show the flow of the proof.

\begin{table}[!h]\label{T}
\caption{Outline of the proof for the center conditions I-IV.}
{\scriptsize\begin{center}
{\begin{tabular}{|l|l|l|l|l|c|}
\hline
\multicolumn{5}{|c|}{Cases}  &Conditions  \\
\hline
\multicolumn{5}{|c|}{(i) $a_{02}>0$}  &I \\[1mm]
\hline
\multirow{10}*{(ii) $a_{02}=0$}
&\multicolumn{4}{c|}{(ii-1) $p_{12}=-\tfrac{b_{21}}{2}$}  &II\\[1.0mm]
\cline{2-6}
&\multicolumn{4}{c|}{(ii-2) $b_{12}=-a_{21}$}  &II\\[1.0mm]
\cline{2-6}
&\multirow{8}*{(ii-3) $p_{02}=p_{12}$} &\multirow{6}*{(ii-3-1) $p_{12}=-1$}
&\multirow{2}*{(ii-3-1-1) $b_{12}=-a_{21}$}
& (ii-3-1-1-1)  & II,\ III\\[1mm] \cline{5-6}
& & & & (ii-3-1-1-2)  & --- \\[1mm]
\cline{4-6}
& & &  \multirow{4}*{(ii-3-1-2) $b_{12}\ne-a_{21}$} & (ii-3-1-2-1) & II\\[1mm]
\cline{5-6}
& & & & (ii-3-1-2-2) & II \\[1mm]
\cline{5-6}
& & & & (ii-3-1-2-3) & --- \\[1mm]
\cline{5-6}
& & & & (ii-3-1-2-4) & II,\ IV\\[1mm]
\cline{3-6}
& &  \multirow{2}*{(ii-3-2) $q_{02}=0$, $1+p_{21}\neq0$,}
&\multicolumn{2}{c|}{(ii-3-2-1) $b_{12}=-a_{21}$} & --- \\[1mm]
\cline{4-6}
& &\hspace*{0.36in} $b_{21}=\tfrac{1}{3}p_{12}p_{21}(2p_{21} \!-\! 1)$
& \multicolumn{2}{c|}{(ii-3-2-2) $b_{12}\ne-a_{21}$} & --- \\[1mm]
\hline
\end{tabular}}
\end{center}}
\end{table}

The basic idea of the proof is setting the generalized Lyapunov constants
zero to get a number of ``necessary'' center conditions,
and excluding those which
make certain order Lyapunov constant non-zero, and then later we prove
that the remaining necessary conditions are also sufficient.

With the aid of a computer algebra system, we use the
Poincar\'{e}-Lyapunov method to compute the generalized
Lyapunov constants associated with the origin of system \eqref{Eqn30}.
The first two generalized Lyapunov constants are $V_1(\varepsilon)=0$, and
\begin{equation*}
\begin{aligned}
V_2(\varepsilon)=&\frac{8}{3}\varepsilon\, \big[
b_{02}+(q_{02}-q_{12})\varepsilon^2 \big].
\end{aligned}
\end{equation*}
Setting the $\varepsilon$-order and $\varepsilon^3$-order terms in
$V_2(\varepsilon)$ zero yield the necessary center conditions,
$$
b_{02}=q_{02}-q_{12}=0.
$$
Then, the $3$rd generalized Lyapunov constant is given by
\begin{equation*}
\begin{aligned}
V_3(\varepsilon)=&\frac{\pi}{4}\varepsilon
\big\{a_{12}(a_{21}+b_{12})+ \big[
2(a_{21}+b_{12})p_{12}
+ 3b_{03} +a_{12}(1 +2 p_{21}+ 2 q_{02}) -2b_{12}b_{21} \big]
\varepsilon^2\\
&\qquad +\big[ 3q_{03} - 2b_{21}(1 + q_{02}) +p_{12} (1+2p_{21} +2 q_{02})
\big]\varepsilon^4 \big\}.
\end{aligned}
\end{equation*}

From the previous analysis for the multiplicity, we only consider
two cases (i) $a_{12}>0$ and (ii) $a_{12}=0$.

\begin{enumerate}
\item[{\bf (i)}]
Assume that $a_{12}>0$. Letting the $\varepsilon$-order,
$\varepsilon^{3}$-order and $\varepsilon^5$-order terms
in $V_3(\varepsilon)$ equal zero we obtain the conditions,
\begin{equation*}
\begin{aligned}
b_{12}&=-a_{21},\quad
b_{03}=-\frac{1}{3} \big[ a_{12} (1 +2p_{21}+2q_{02}) +2a_{21}b_{21} \big],
\\
q_{03}&=\frac{1}{3} \big[ 2 b_{21} (1 + q_{02})
-p_{12} (1 + 2p_{21} +2q_{02}) \big].
\end{aligned}
\end{equation*}
Then, we have the 4th generalized Lyapunov constant, given by
\begin{equation*}
\begin{aligned}
V_4(\varepsilon)=&-\frac{16}{45}\varepsilon^3\,
\big[a_{12}-(p_{02}-p_{12})\varepsilon^2\big] \big\{ 4a_{12}(p_{21}+q_{02})
+ \big[3b_{21}+2(b_{21}+2p_{12})(p_{21}+q_{02})\big]\varepsilon^2\big\}.
\end{aligned}
\end{equation*}
Thus, we obtain the conditions $p_{21}+q_{02}=b_{21}=0$ by setting
the $\varepsilon^{3}$-order and $\varepsilon^{7}$-order
terms in $V_4(\varepsilon)$ zero. Therefore, the necessary center
conditions are $a_{21}+b_{12}=a_{12}+3b_{03}=b_{21}=0$,
giving the condition I.
\end{enumerate}

\begin{enumerate}
\item[{\bf (ii)}]
Assume that $a_{12}=0$. By setting the $\varepsilon^{3}$-order
and $\varepsilon^5$-order terms
in $V_3(\varepsilon)$ zero we obtain the conditions,
\begin{equation*}
\begin{aligned}
b_{03}=\frac{2}{3} \big[b_{12}b_{21}-p_{12}(a_{21}+b_{12}) \big],\quad
q_{03}=\frac{1}{3} \big[ 2b_{21}(1 + q_{02})
-p_{12}(1+2p_{21} +2q_{02}) \big].
\end{aligned}
\end{equation*}

Then, we obtain the $4$th generalized Lyapunov constant,
\begin{equation*}
\begin{aligned}
V_4(\varepsilon)=&\frac{16}{45} \, \varepsilon^5
(p_{02}-p_{12})\big\{2(a_{21}+b_{12})(b_{21}+2p_{12})
+\big[3b_{21}+2(b_{21}+2p_{12})(p_{21}+q_{02})\big]\varepsilon^2 \big\}.
\end{aligned}
\end{equation*}

We have the following three subcases (ii-1) $b_{21}+2p_{12}=0$,
(ii-2) $a_{21}+b_{12}=0$ and (ii-3) $p_{02}-p_{12}=0$.
Note that in the following analysis, when the condition in one case
is satisfied, the conditions in other two cases may hold or not
hold, depending upon the analysis for each case.
In subcase (ii-3-1-1), both conditions $p_{02}-p_{12}=0$ and
$a_{21}+b_{12}=0$ are satisfied.

\vspace{0.10in}
{\bf (ii-1)} Assume that $p_{12}=-\frac{b_{21}}{2}$.
Then, we have $b_{21}=0$ by setting $V_4(\varepsilon)=0$,
yielding the condition II.

\vspace{0.10in}
{\bf (ii-2)} Assume that $b_{12}=-a_{21}$. Since the parameter $q_{02}$
is redundant for the analysis of this subcase, without loss of generality,
we let $q_{02}=0$. If $p_{21}=0$,
we have $b_{21}=0$
by using $V_4(\varepsilon)=0$, which gives a special case of the condition II.
Otherwise, if $p_{21}\ne 0$,
we have $p_{12}=-\frac{b_{21}(3+2p_{21})}{4p_{21}}$
by setting $\varepsilon^7$-order term of $V_4(\varepsilon)$ zero.
Then, the 5th generalized Lyapunov constant becomes
\begin{equation*}
\begin{aligned}
V_5(\varepsilon)=&
\frac{1}{768}b_{21}\pi \, \varepsilon^7
\big\{16p_{21}
\big[20a_{21}b_{21}^2
- 3 (a_{03}+2a_{21}^2-8a_{21}b_{21}^2)p_{21}
+ 2 (3a_{03}-10a_{21}^2) p_{21}^2
\big]
\\
& +\big[105b_{21}^2
+ 4 b_{21} (17b_{21} -90p_{02}) p_{21}
-4 (96a_{21} \!+\! 7 b_{21}^2 \!+\! 60b_{21}p_{02}
\!+\! 60p_{02}^2 \!+\! 12p_{03}) p_{21}^2
\\[0.5ex]
& \quad
+ 32 (17 a_{21} - 3 b_{21}^2) p_{21}^3 -448a_{21}p_{21}^4\big]\varepsilon^2
-32 p_{21}^2(1+p_{21}) \big(9+4p_{21}+4p_{21}^2 \big)\varepsilon^4 \big\}.
\end{aligned}
\end{equation*}
Setting the $\varepsilon^{11}$-order term in $V_5(\varepsilon)$ zero
we obtain $p_{21}=-1$. Further, letting $V_5(\varepsilon)=0$ we have
\begin{equation*}
\begin{aligned}
a_{03}=\frac{2}{9}(7a_{21}^2+2a_{21}b_{21}^2),\quad
a_{21}=-\frac{1}{96}(29b_{21}^2-40b_{21}p_{02}+80p_{02}^2+48p_{03}).
\end{aligned}
\end{equation*}
Letting the $\varepsilon^{12}$-order term in $V_6(\varepsilon)$ zero
we obtain $p_{03}=0$. Then, the 6th generalized Lyapunov constant becomes
\begin{equation*}
\begin{aligned}
V_6(\varepsilon)=&-\frac{1}{18900}\, b_{21}(b_{21}
-4p_{02})\,\varepsilon^9\big[3V_{6a}+5V_{6b}\varepsilon^2\big],
\end{aligned}
\end{equation*}
where
\begin{equation*}
\begin{aligned}
V_{6a}=&\ (7b_{21}^2+40b_{21}p_{02}-80p_{02}^2)(29b_{21}^2
-40b_{21}p_{02}+80p_{02}^2),\\
V_{6b}=&\ 551b_{21}^2+1544b_{21}p_{02}-3088p_{02}^2.
\end{aligned}
\end{equation*}
By computing the resultant of $V_{6a}$ and $V_{6b}$ with respect to $p_{02}$ we have
\begin{equation*}
\begin{aligned}
{\rm Res}[V_{6a},V_{6b},p_{02}]=9011459133227925504\, b_{21}^8\neq 0, \quad
(b_{21} \ne 0),
\end{aligned}
\end{equation*}
which implies that $V_{6a}$ and $V_{6b}$ have no common roots.
Hence, we obtain $p_{02}=\frac{b_{21}}{4}$ from $V_6(\varepsilon)=0$.
Then, we obtain the 7th generalized Lyapunov constant,
\begin{equation*}
\begin{aligned}
V_7(\varepsilon)=&-\frac{7}{82944}b_{21}^5\pi\,\varepsilon^9
(2b_{21}^2+9\varepsilon^2)(8b_{21}+75\varepsilon^2),
\end{aligned}
\end{equation*}
which indicates that the $\varepsilon^9$-order term
in $V_7(\varepsilon)$ is non-zero when $b_{21}\ne0$.

\vspace{0.10in}
{\bf (ii-3)} Assume that $p_{02}=p_{12}$. Then, we have $V_4(\varepsilon)=0$
and obtain the 5th generalized Lyapunov constant,
\begin{equation*}
\begin{aligned}
V_5(\varepsilon)=\frac{\pi}{36}\varepsilon^5[3(a_{21}+b_{12})V_{5a}
+V_{5b}\varepsilon^2-V_{5c}\varepsilon^4-6V_{5d}\varepsilon^6],\\
\end{aligned}
\end{equation*}
where
\begin{equation*}
\begin{aligned}
V_{5a}=& \ 4a_{03}b_{21}
+(5a_{03}+4a_{21}^2-6a_{21}b_{12}) p_{12},\\
V_{5b}=&\ b_{21} (9a_{03}-18a_{21}b_{12}-8b_{12}b_{21}^2)
+3(a_{21}+ b_{12}) ( 4 b_{21} + 5 p_{12}) p_{03}\\
&+ 2 [3 (a_{21}+ b_{12})(a_{21}- 3 b_{12})
-2 (a_{21}-4b_{12}) b_{21}^2] p_{12}
+10 b_{21} (a_{21} + b_{12}) p_{12}^2\\
&+ 3 \big[ 4 a_{03}b_{21}+ (5a_{03} +12 a_{21}^2- 4 a_{21}b_{12}- 6 b_{12}^2) p_{12} \big] p_{21}
\\
& +\big[6 a_{21} (a_{21}+ 6 b_{12}) p_{12}
- 3 a_{03} (4 b_{21} + 5 p_{12})\big]q_{02},\\
V_{5c}=&\ 2 b_{21} (9a_{21}+9b_{12}+4b_{21}^2)
-3 b_{21} (3 +4 p_{21} ) p_{03}
+6 (a_{21} +b_{12}-2b_{21}^2) p_{12}\\
&+18b_{12}b_{21}p_{21}
-10 (2 b_{21} +b_{21} p_{21} ) p_{12}^2
- 6 (6a_{21} - b_{12}) p_{12} p_{21}^2-( 12a_{21} -12b_{12}\\
&-4b_{21}^2 +15p_{03}) p_{12}p_{21} +18 a_{21} p_{12} q_{02}^2
+\big[8 b_{21}^3 - 16 b_{21}^2 p_{12}
- 2 b_{21} (6 p_{03} + 5 p_{12}^2) \\
& -3 p_{12} (5 p_{03} - 12 b_{12}- 12b_{12}p_{21})
+ 6 a_{21} (3 b_{21} + 2 p_{12} + 2p_{12}p_{21})\big]q_{02} ,\\
V_{5d}=&\ (1+p_{21})\big[ 3b_{21}+p_{12}p_{21} (1 -2p_{21})
+(3 b_{21} + p_{12} + p_{12} p_{21}) q_{02} + 3 p_{12} q_{02}^2\big].
\end{aligned}
\end{equation*}

Two subcases follow $V_{5d}=0$: (ii-3-1) $p_{21}=-1$ and  (ii-3-2)
$$
b_{21}=\frac{1}{3}\big[ p_{12}p_{21}(2p_{21}-1)+(3 b_{21} + p_{12}
+ p_{12} p_{21}) q_{02} + 3 p_{12} q_{02}^2\big].
$$

\vspace{0.10in}
{\bf (ii-3-1)} When $p_{21}=-1$ we have $V_{5d}=0$.
In order to have  that the $\varepsilon^5$-order term in
$V_5(\varepsilon)$ equals zero, it needs $b_{12}=-a_{21}$ or $V_{5a}=0$.

\vspace{0.10in}
{\bf (ii-3-1-1)} Assume that $b_{12}=-a_{21}$, we have $a_{03}+2 a_{21}^2<0$.
Then, we obtain the simplified $V_{5b}$,
\begin{equation*}
\begin{aligned}
V_{5b}=&-b_{21} \big[2 a_{21} (9 a_{21} + 4 b_{21}^2)
+ 3 a_{03} (4 q_{02}-1)\big]\\
&+5\big[(3 a_{03} + 6 a_{21}^2 + 4 a_{21} b_{21}^2) - 3(a_{03}
+ 2 a_{21}^2) q_{02}\big]p_{12}.
\end{aligned}
\end{equation*}

{\bf (ii-3-1-1-1)} If
$$
q_{02}=\frac{3 a_{03} + 6 a_{21}^2 + 4 a_{21} b_{21}^2}{3 (a_{03}
+ 2 a_{21}^2)},
$$
from $V_{5b}=0$ we have $b_{21}=0$, which gives a special case of
the condition II, or
\begin{equation*}
\begin{aligned}
M_1=9 (a_{03}+2 a_{21}^2)^2 + 8 a_{21}b_{21}^2 (3 a_{03} + 2 a_{21}^2)=0,
\end{aligned}
\end{equation*}
which does not contain the perturbed parameters, as expected,
giving the condition III. Further there exists the free
parameter $p_{03}$ such that $V_{5c}=0$.

{\bf (ii-3-1-1-2)} If
$$
q_{02}\ne\frac{3 a_{03} + 6 a_{21}^2
+ 4 a_{21} b_{21}^2}{3 (a_{03} + 2 a_{21}^2)},
$$
we have
\begin{equation*}
\begin{aligned}
p_{12}=\frac{b_{21} \big[2 a_{21} (9 a_{21} + 4 b_{21}^2)
+ 3 a_{03} (4 q_{02}-1)\big]}{5 (3 a_{03} + 6 a_{21}^2
+ 4 a_{21} b_{21}^2) - 15 (a_{03} + 2 a_{21}^2) q_{02}}
\end{aligned}
\end{equation*}from $V_{5b}=0$. Letting
$$
M_2= -2a_{21}(1 - q_{02})^2 - b_{21}^2 (1- 2 q_{02} ) =0
$$
from $V_{5c}=0$, which yields $V_5(\varepsilon)=V_6(\varepsilon)=0$.
Then, from the $\varepsilon^9$-order term of $V_7(\varepsilon)$,
\begin{equation*}
\begin{aligned}
V_{7a}=&\, b_{21}M_1=0\\
 \end{aligned}
\end{equation*}which also giving the special case of the condition II and the condition III. In fact, from the  $\varepsilon^{11}$-order term $V_{7b}$ and $\varepsilon^{13}$-order term $V_{7c}$ of $V_7(\varepsilon)$, we have the resultants of $M_2$, $V_{7b}$,  and $V_{7c}$
with respect to $q_{02}$, respectively,
$${\rm Res}[M_2,V_{7b},V_{7c},q_{02}]=[0,0],$$
i.e., there exist the free parameters such $V_{7b}=V_{7c}=0$. Then we have $V_7(\varepsilon)=0$.

We remark that since the parameter $q_{02}$ is redundant for the following subcases, without loss of generality, we let $q_{02}=0$.

\vspace{0.10in}
{\bf (ii-3-1-2)} Assume that $b_{12}\ne-a_{21}$.   We solve the
polynomial equations: $V_{5a}=V_{5b}=V_{5c}=0$ to
find other real solutions for parameters $a_{21}$, $a_{03}$,
$b_{12}$ and $b_{21}$.

\vspace{0.10in}
{\bf (ii-3-1-2-1)} If $a_{03} \!=\! -\,\frac{2}{5}a_{21}(2a_{21} \!-\! 3b_{12})
\!<\! 0$, we have $b_{21} \!=\! 0$, and either $p_{12} \!=\! 0$
or $6a_{21}-5p_{03} \!=\! 0$
by using the equations $V_{5a}=V_{5b}=V_{5c}=0$.
The first solution is included in the condition II.
For the second solution, if $p_{03}=\frac{6a_{21}}{5}$ we obtain
$V_6(\varepsilon)=0$ and the 7th generalized Lyapunov constant, given by
\begin{equation*}
\begin{aligned}
V_7(\varepsilon)=&-\frac{\pi}{180}\varepsilon^9a_{21}p_{12}
(6a_{21}-35p_{12}^2)\big[(a_{21}+b_{12})^2-\varepsilon^4\big].
\end{aligned}
\end{equation*}
Setting $V_7(\varepsilon)=0$ yields $a_{21}=\frac{35}{6}p_{12}^2\neq0$,
leading to $V_8(\varepsilon)=0$ and
\begin{equation*}
\begin{aligned}
V_9(\varepsilon)=-\frac{98}{729}p_{12}^7 \pi \varepsilon^{13}
(6b_{12}+35p_{12}^2-6\varepsilon^2)(6b_{12}+35p_{12}^2+6 \varepsilon^2).
\end{aligned}
\end{equation*}
Hence, we obtain that the $\varepsilon^{17}$-order term
in $V_9(\varepsilon)$ is not vanished.

\vspace{0.10in}
{\bf (ii-3-1-2-2)} If $a_{21}=0$, by setting
$V_{5a}=V_{5b}=V_{5c}=V_7(\varepsilon)=0$
we obtain $b_{21}=0$, which is included in the condition II.

\vspace{0.10in}
{\bf (ii-3-1-2-3)} If $2a_{21}-3b_{12}=0$, we obtain
\begin{equation*}
\begin{aligned}
&b_{12}=-\frac{32}{115}b_{21}^2,\quad a_{21}=-\frac{48}{115}b_{21}^2,
\quad a_{03}=-\frac{3072}{13225}b_{21}^4,\quad
p_{12}=-\frac{4}{5}b_{21}, \quad p_{03}=\frac{56}{575}b_{21}^2,
\end{aligned}
\end{equation*}
by solving $V_{5a}=V_{5b}=V_{5c}=V_7(\varepsilon)=0$.
Further, we have $V_8(\varepsilon)=0$ and
\begin{equation*}
\begin{aligned}
V_9(\varepsilon)= \frac{322224\pi}{502839296875}b_{21}^7\varepsilon^{13}
(24576 b_{21}^4-13984 b_{21}^2\varepsilon^2-10051 \varepsilon^4)
\ne 0 \quad \textrm{if} \ b_{21} \ne 0.
\end{aligned}
\end{equation*}
But if
$ b_{21}^2 = \frac{23(19+\sqrt{2185})\varepsilon^2}{1536}$, then $V_9(\varepsilon)=0$.

\vspace{0.10in}
{\bf (ii-3-1-2-4)} Assume that
$a_{21}(2a_{21}-3b_{12})(5a_{03}+4a_{21}^2-6a_{21}b_{12})\ne0$.
If $p_{12}=0$, we have $b_{21}=0$ from $V_{5a}=0$,
which is a special case of the condition II.
Suppose that $p_{12}\ne0$. Then, by using $V_{5a}=V_{5b}=0$,
we have that
\begin{equation*}
\begin{aligned}
&b_{21}=-\frac{1}{a_{03}}\,p_{12}(5a_{03}+4a_{21}^2-6a_{21}b_{12})\neq0,\\
&p_{03}=\frac{1}{48a_{03}^3a_{21}(a_{21}+b_{12})(2a_{21}-3b_{12})}\,
M_3,
\end{aligned}
\end{equation*}
where
\begin{equation*}
\begin{aligned}
M_{3}=&\ 18 a_{03}^2 \big[5a_{03}^2 +4a_{03}a_{21}
(3a_{21}-2 b_{12})-4 a_{21}^2 b_{12} (2 a_{21} - 3 b_{12})\big] \\
& + p_{12}^2 (2 a_{21}-3b_{12}) (a_{03}-2a_{21}b_{12})
(5 a_{03}+4 a_{21}^2-6a_{21}b_{12})
(15a_{03}+4a_{21}^2-6a_{21}b_{12}).
\end{aligned}
\end{equation*}
Then, the polynomial $V_{5c}$ becomes
\begin{equation*}
\begin{aligned}
V_{5c}=\frac{(5a_{03}+4a_{21}^2-6a_{21}b_{12})p_{12}}
{64 a_{03}^3a_{21}(a_{21}+b_{12})(2a_{21}-3b_{12})}
M_{4}M_{5},
\end{aligned}
\end{equation*}
where
\begin{equation*}
\begin{aligned}
M_{4}=&\ 15 a_{03}^2 +  4  a_{21}^2
\big[ 5 a_{03} + (2 a_{21}- 3 b_{12})(4 a_{21} +5 b_{12}) \big],\\[0.5ex]
M_{5}=&\ 18 a_{03}^2+\big[15 a_{03}(2a_{21}-3b_{12})
+2a_{21}(2a_{21}-3b_{12})^2\big]p_{12}^2.
\end{aligned}
\end{equation*}

To have $V_{5c} = 0$, it requires $M_{4}M_{5}=0$ under which
$V_6(\varepsilon)=0$, and the 7th generalized Lyapunov constant becomes
\begin{equation*}
\begin{aligned}
V_7(\varepsilon)=&\ \frac{\pi\varepsilon^9}
{3538944 a_{03}^5a_{21}^2(2a_{21}-3b_{12})^2(a_{21}+b_{12})^2}
\\
& \times \big[1536 a_{03}^2 a_{21}^2(2 a_{21}-3 b_{12})^2 (a_{21}
+ b_{12})^2(a_{03}-2 a_{21} b_{12})\widetilde{V}_{7a}
\\
&\quad \ +16 a_{21}(2a_{21}-3b_{12})(a_{21}+b_{12})(5a_{03}+4a_{21}^2
-6a_{21}b_{12})\widetilde{V}_{7b}\varepsilon^2
-\widetilde{V}_{7c}\varepsilon^4\big],
\end{aligned}
\end{equation*}
where
\begin{equation*}
\begin{aligned}
\widetilde{V}_{7a}=& \ 630a_{03}^4+15a_{03}^3\big[8 a_{21}
(8 a_{21} - 13 b_{12}) + (46 a_{21} - 129 b_{12}) p_{12}^2\big]\\
& +2a_{03}^2a_{21}(2 a_{21} - 3 b_{12})\big[36 a_{21}
(2 a_{21} - 5 b_{12}) + 25 (22 a_{21} - 27 b_{12}) p_{12}^2\big]\\
& +140a_{03}a_{21}^2(2a_{21}-3 b_{12})^3p_{12}^2+8a_{21}^3
(2 a_{21}-5b_{12})(2a_{21}-3b_{12})^3p_{12}^2,\\
\widetilde{V}_{7b}=& \ 5670 a_{03}^7+135 a_{03}^6\big[2368a_{21}^2+6a_{21}
(986b_{12}-19p_{12}^2)+b_{12}(2792b_{12}+411p_{12}^2)\big]\\
&+18a_{03}^5 \big[4 a_{21}^2 (6396 a_{21}^2+16316 a_{21} b_{12}+3305 b_{12}^2)
+10 (2614 a_{21}^3 + 4917 a_{21}^2 b_{12}\\
&-9407a_{21}b_{12}^2-5235 b_{12}^3) p_{12}^2-25(46 a_{21}-129 b_{12})
(2a_{21}-3b_{12})p_{12}^4\big]\\
&+12a_{03}^4 a_{21} (2 a_{21} - 3 b_{12})\big[12 a_{21} (2696 a_{21}^3
+ 8966 a_{21}^2 b_{12} + 9067 a_{21} b_{12}^2 + 3490 b_{12}^3)\\
&+30(1122 a_{21}^3+3451a_{21}^2 b_{12}-540a_{21}b_{12}^2-349 b_{12}^3)p_{12}^2
-5 (1376 a_{21}^2 - 4878 a_{21} b_{12}\\
&+ 5121 b_{12}^2) p_{12}^4\big]+16a_{03}^3a_{21}^2(2 a_{21}-3b_{12})^2p_{12}^2
\big[22882 a_{21}^3 + 78167 a_{21}^2 b_{12}\\
&+ 65395 a_{21} b_{12}^2 + 26175 b_{12}^3 - 15 (180 a_{21}^2
- 854 a_{21} b_{12} + 891 b_{12}^2) p_{12}^2\big] \\
& +16a_{03}^2a_{21}^3(2a_{21}-3b_{12})^3p_{12}^2\big[2696 a_{21}^3
+ 9182 a_{21}^2 b_{12} + 8527 a_{21} b_{12}^2 + 3490 b_{12}^3 \\
& - 30 (16 a_{21}^2 - 140 a_{21} b_{12}+ 159 b_{12}^2) p_{12}^2\big]
+192a_{21}^5(2 a_{21}-5b_{12})(2a_{21}-3b_{12})^5 b_{12}p_{12}^4 \\
& -96 a_{03}a_{21}^4 (2 a_{21} -3 b_{12})^4
(4 a_{21}^2 -96a_{21}b_{12}+145b_{12}^2)p_{12}^4,
\end{aligned}
\end{equation*}and
\begin{equation*}
\begin{aligned}
\widetilde{V}_{7c}=&\ 4252500 a_{03}^9+56700a_{03}^8\big[6 a_{21} (46 a_{21}
- 19 b_{12}) + 125 (2 a_{21} - 3 b_{12}) p_{12}^2\big]\\
&+135 a_{03}^7 \big[32 a_{21} (11934 a_{21}^3 - 682 a_{21}^2 b_{12}
- 17551 a_{21} b_{12}^2 - 6510 b_{12}^3) + 440 a_{21} (337 a_{21}\\
&- 363 b_{12}) (2 a_{21} - 3 b_{12}) p_{12}^2 +
 21875 (2 a_{21} - 3 b_{12})^2 p_{12}^4\big]
\\
& +90 a_{03}^6 a_{21}
(2 a_{21} - 3 b_{12}) \big[32 a_{21} (14258 a_{21}^3 + 18094 a_{21}^2 b_{12}
-10465 a_{21} b_{12}^2 - 3906 b_{12}^3)
\\
& +8 (71196 a_{21}^3-40193 a_{21}^2 b_{12} - 94064 a_{21} b_{12}^2 -
32550 b_{12}^3) p_{12}^2 \\
&+ 25 (2584 a_{21} -  5991 b_{12}) (2 a_{21} - 3 b_{12}) p_{12}^4\big]
+12 a_{03}^5 a_{21}^2 (2 a_{21} - 3 b_{12})^2
\\
& \times \big[48 a_{21} (42952 a_{21}^3
+ 110952 a_{21}^2 b_{12} + 67965 a_{21} b_{12}^2
+ 21700 b_{12}^3)\\
&+80 (38078 a_{21}^3 + 31157 a_{21}^2 b_{12} - 70341 a_{21} b_{12}^2
-13020 b_{12}^3) p_{12}^2 +25 (15556 a_{21}^2\\
&-77328 a_{21} b_{12} - 4509 b_{12}^2) p_{12}^4\big]
+8a_{03}^4a_{21}^3(2a_{21}-3 b_{12})^3
\big[144 a_{21} (5248 a_{21}^3 \\
& + 15748 a_{21}^2 b_{12} +13517 a_{21} b_{12}^2 +4340 b_{12}^3)
+40 (67916 a_{21}^3 + 164530 a_{21}^2 b_{12}
\\
& + 39263 a_{21} b_{12}^2 +28644 b_{12}^3) p_{12}^2
-15 (19608 a_{21}^2 + 110966 a_{21} b_{12} + 242733 b_{12}^2) p_{12}^4\big]
\\
& +16 a_{03}^3a_{21}^4(2 a_{21}-3b_{12})^4p_{12}^2\big[8 (53754 a_{21}^3
+165299 a_{21}^2 b_{12} \!+\! 113365 a_{21} b_{12}^2 \!+\!
 43400 b_{12}^3)
\\
&- 15 (12172 a_{21}^2 + 25588 a_{21} b_{12}
+ 72391 b_{12}^2) p_{12}^2\big]
+32 a_{03}^2 a_{21}^5 (2 a_{21}-3 b_{12})^5p_{12}^2
\\
& \times \big[8 (2156 a_{21}^3
+ 7415 a_{21}^2 b_{12} + 6106 a_{21} b_{12}^2 + 2170 b_{12}^3)
-5 (6944 a_{21}^2 + 9626 a_{21} b_{12}
\\
& + 26097 b_{12}^2) p_{12}^2\big]
-64 a_{03} a_{21}^6 (2 a_{21}-3b_{12})^6p_{12}^4(3276 a_{21}^2
+3656 a_{21} b_{12} +7205 b_{12}^2)
\\
& -128 a_{21}^7 p_{12}^4(2 a_{21}
-3 b_{12})^7 (104 a_{21}^2+ 102 a_{21} b_{12} + 145 b_{12}^2).
\end{aligned}
\end{equation*}

Since $a_{03}<-\frac{1}{2}(b_{12}-a_{21})^2$,
we have $a_{03}-2a_{21}b_{12}<0$. Then, to solve $V_7(\varepsilon)=0$,
we only need to find the solutions to the equations:
$M_{4}M_{5}=\widetilde{V}_{7a}=\widetilde{V}_{7b}=\widetilde{V}_{7c}=0$.
Thus, we compute the resultants of
$M_{4}$, $M_{5}$, $\widetilde{V}_{7a}$,
$\widetilde{V}_{7b}$ and $\widetilde{V}_{7c}$
with respect to $a_{03}$, respectively, and obtain
\begin{equation*}
\begin{aligned}
&\ {\rm Res}[M_{4},\widetilde{V}_{7a},\widetilde{V}_{7b},
\widetilde{V}_{7c},a_{03}]\\
=&\ a_{21}^6(2 a_{21}-3 b_{12})^3(a_{21} + b_{12})N_4\big[C_1,\
C_2 a_{21}^{5}(2 a_{21}-3 b_{12})^3N_3,\\
&\hspace*{2.12in}
C_3 a_{21}^{9}(2 a_{21}-3 b_{12})^6(a_{21} + b_{12})N_3\big],\\[1.0ex]
&\ {\rm Res}[M_{5},\widetilde{V}_{7a},\widetilde{V}_{7b},
\widetilde{V}_{7c},a_{03}]\\
=&\ a_{21}^3(2 a_{21}-3 b_{12})^6(a_{21}
+ b_{12})^2N_1\big[C_4,\ C_5 a_{21}^{3}(2 a_{21}-3 b_{12})^{4}N_2,\\
&\hspace*{2.20in}
C_6 a_{21}^{6}(a_{21}-1)(25 +2a_{21})(2 a_{21}-3 b_{12})^{8}(a_{21}
+ b_{12})^2\big],
\end{aligned}
\end{equation*}
where $C_i$, $i=1,2,\cdots,6$, are real constants, and
\begin{equation*}
\begin{aligned}
N_1=&\ (2a_{21}-3p_{12}^2)(2 a_{21}+3p_{12}^2),\\
N_2=&\ 2 a_{21} (2 a_{21} + 3 b_{12})^2 - (8 a_{21} + 3 b_{12})
(a_{21} + 6 b_{12}) p_{12}^2,\\
N_3=&\ 8 a_{21}^2 (4 a_{21} + 5 b_{12})^2 +
 5 a_{21} (13 a_{21}^2 + 78 a_{21} b_{12} + 85 b_{12}^2) p_{12}^2 +
 50 (a_{21} + b_{12})^2 p_{12}^4,\\
N_4=&\ 12 a_{21}^2 (4 a_{21} + 5 b_{12})^2 (107 a_{21}^2
- 40 a_{21} b_{12} - 175 b_{12}^2)
\\
& - 10 a_{21} (4682 a_{21}^4+ 4193 a_{21}^3 b_{12}
- 16315 a_{21}^2 b_{12}^2 - 23925 a_{21} b_{12}^3 - 7875 b_{12}^4) p_{12}^2
\\
& + 5 (5954 a_{21}^4- 7331 a_{21}^3 b_{12} - 21660 a_{21}^2 b_{12}^2
+ 33975 a_{21} b_{12}^3 + 47250 b_{12}^4) p_{12}^4.
\end{aligned}
\end{equation*}
To find the common factors of the above resultants,
we only need to consider the following three subcases
due to $a_{21}(2 a_{21}-3 b_{12})(a_{21} + b_{12})\neq0$.

\vspace{0.10in}
{\bf (ii-3-1-2-4-1)} If $a_{21}=\frac{3}{2}p_{12}^2$,
we obtain the simplified polynomials,
\begin{equation*}
\begin{aligned}
M_{5}=&\ 9 (a_{03} - b_{12} p_{12}^2 + p_{12}^4)
(2 a_{03} - 3 b_{12} p_{12}^2 + 3 p_{12}^4),\\
\widetilde{V}_{7a}=& \ 9 (a_{03} - b_{12} p_{12}^2 + p_{12}^4)
\big[70 a_{03}^3 - 81 p_{12}^6 (b_{12}- p_{12}^2)^2 (5 b_{12}-3 p_{12}^2)\\
&+ 54 a_{03} p_{12}^4 (b_{12} - p_{12}^2) (10 b_{12} - 13 p_{12}^2)
- 15 a_{03}^2 p_{12}^2 (27 b_{12} - 19 p_{12}^2)\big].
\end{aligned}
\end{equation*}
If $a_{03}=(b_{12} -p_{12}^2) p_{12}^2$, we have $M_{5}=\widetilde{V}_{7a}=0$,
which yields $\Delta^+=\frac{1}{4}(2 b_{12} - p_{12}^2)^2>0$
violating the condition $\Delta^+<0$.
If $a_{03}=\frac{3}{2}(b_{12}- p_{12}^2)p_{12}^2$,
we obtain $b_{12}=-\frac{3}{2}p_{12}^2$ from $\widetilde{V}_{7a}=0$,
which contradicts the condition $b_{12}\neq-a_{21}$.

\vspace{0.10in}
{\bf (ii-3-1-2-4-2)} Assume that $a_{21}=-\frac{3}{2}p_{12}^2$.
Similar to the above process, we use
the equation $M_{5}=\widetilde{V}_{7a}=0$ with the constraint
$\Delta^+<0$ to obtain one set of parameter values:
\begin{equation*}
\begin{aligned}
\{ 2a_{21}+3p_{12}^2=a_{03}-3(b_{12}+p_{12}^2)p_{12}^2
=b_{21}+2p_{12}=p_{03}-3p_{12}^2=0\}.
\end{aligned}
\end{equation*}
Eliminating the perturbed parameter $p_{12}$ we obtain the condition IV.

\vspace{0.10in}
{\bf (ii-3-1-2-4-3)} Now assume $N_4=0$ and $p_{12}\ne0$.
We solve the equations:
$M_4=\widetilde{V}_{7a}=\widetilde{V}_{7b}=\widetilde{V}_{7c}=0$
to obtain the solutions under which
$V_8(\varepsilon)=0$ and the $9$th generalized
Lyapunov constant is simplified to
\begin{equation*}
\begin{aligned}
V_{9}(\varepsilon)=&\ \frac{-\,\pi\, p_{12}\, \varepsilon^{13} }
{122305904640 a_{03}^7 a_{21}^3
(2 a_{21} - 3 b_{12})^3 (a_{21} + b_{12})^3}
\big[ 256 a_{03}^2 a_{21}^2 (2 a_{21}- 3 b_{12})^2
\\
& \times (a_{21} + b_{12})^2 V_{9a}
+ 32 a_{21} (2 a_{21} - 3 b_{12}) (a_{21} + b_{12}) V_{9b}\varepsilon^{2}
-3V_{9c}\varepsilon^{4}\big],
\end{aligned}
\end{equation*}
where $V_{9a}$, $V_{9b}$ and $V_{9c}$ are polynomials in $a_{03}$,
$a_{21}$, $b_{12}$ and $p_{12}$, with $221$, $317$ and $299$ terms,
respectively. We compute the resultants of $M_4$, $V_{9a}$, $V_{9b}$
and $M_{9c}$ with respect to $a_{03}$, respectively, and obtain
\begin{equation*}
\begin{aligned}
{\rm Res}&[M_4, V_{9a}, V_{9b}, V_{9c}, a_{03}]
=\big[C_7 a_{21}^{13}(2 a_{21}-3 b_{12})^7(a_{21} + b_{12})^3N_5,\\
&\ C_8 a_{21}^{19}(2 a_{21}-3 b_{12})^{11}(a_{21}
+ b_{12})^3N_6,\ C_9 a_{21}^{21}(2 a_{21}-3 b_{12})^{13}(a_{21}
+ b_{12})^3N_3N_7\big],\\
\end{aligned}
\end{equation*}
where $C_i$, $i=7,8,9$, are real constants, and $N_5$, $N_6$ and $N_7$
are polynomials in $a_{21}$, $b_{12}$ and $p_{12}$,
having 45 terms, 81 terms and 47 terms, respectively.

We again calculate the resultants of
$N_4$, $N_5$, $N_6$, $N_3N_7$
with respect to $a_{21}$, respectively, to have
\begin{equation*}
\begin{array}{rl}
{\rm Res}[N_4, N_5, N_6, N_3N_7, a_{21}]= \!\!\!\! &b_{12}^{36}
(2 b_{12}-25p_{12}^2)^2 (10 b_{12}+11p_{12}^2)p_{12}^{32} N_8
\\
\!\!\!\! &
\times\big[C_{10}
N_9,\ C_{11} (2 b_{12}-25) (3150 b_{12}-143) p_{12}^{4}N_{10}N_{11},\\
&
\quad \, C_{12} (2 b_{12}-25) (3150 b_{12}-143)p_{12}^{4}
N_{10}N_{12}N_{13}\big],\\
\end{array}
\end{equation*}
where $C_i$, $i=10,11,12$, are real constants, and
\begin{equation*}
\begin{aligned}
N_8=&\, -86943853334p_{12}^{12}-234718996369b_{12}p_{12}^{10}
+107204353618b_{12}^2p_{12}^{8}\\
&+414472065816 b_{12}^3p_{12}^{6}-88238153808b_{12}^4p_{12}^{4}
-2186887248 b_{12}^5p_{12}^{2}+792148896b_{12}^6,
\end{aligned}
\end{equation*}
while the polynomials $N_{i}(b_{12})$, $i=9,12,\cdots,13$, do not
have common nonzero roots. So, we only have possible solutions from
$ (2 b_{12}-25p_{12}^2)^2 (10 b_{12}+11p_{12}^2)N_8$.
If $b_{12}=\frac{25}{2}p_{12}^2$, we solve $N_4=N_5=0$ to obtain
$a_{21}=-\frac{25}{2}p_{12}^2$, which yields $ a_{21}+b_{12}=0$,
contradicting the condition $a_{21}+b_{12}\ne0$.
If $b_{12}=-\frac{11}{10}p_{12}^2$, we solve $M_{4}=N_4=N_5=0$ to
have $a_{21}=\frac{3}{2}p_{12}^2$, and either $a_{03}=-\frac{9}{10}p_{12}^2$ or
$a_{03}=-\frac{21}{10}p_{12}^2$, which violates $\Delta^+<0$.

Now, assume that $N_8=0$. We have $V_{9}(\varepsilon)
=V_{10}(\varepsilon)=0$, and obtain the $\varepsilon^{17}$-order
term $V_{11a}$ in the $11$th generalized Lyapunov constant,
which is a polynomial in $a_{03}$, $a_{21}$ and $b_{12}$
having $237$ terms.
We compute the resultant of $M_4$ and $V_{11a}$ with respect to
$a_{03}$ to obtain
\begin{equation*}
\begin{aligned}
&{\rm Res}[M_4, V_{11a}, a_{03}]=C_{13} a_{21}^{19}
(2 a_{21}-3 b_{12})^{11}(a_{21} + b_{12})^4N_{14},\\
\end{aligned}
\end{equation*}
where $C_{13}$ is a real constant, and $N_{14}$ is a polynomial
in $a_{21}$ and $b_{12}$ having $179$ terms.
We again calculate the resultant of
$N_4$ and $N_{14}$ with respect to $a_{21}$, and have
\begin{equation*}
\begin{aligned}
{\rm Res}[N_4, N_{14}, a_{21}]&=C_{14} b_{12}^{48}p_{12}^{36}
(2 b_{12}-25 p_{12}^2)^3 (10 b_{12}+11p_{12}^2)N_{15},
\end{aligned}
\end{equation*}
where $C_{14}$ is a real constant, and $N_{15}$ is a $76$th-degree
polynomial in $b_{12}$ and $p_{12}$.
Since the resultant of
$N_8$ and $N_{15}$ with respect to $b_{12}$ equals $C_{15}p_{12}^{456}\ne0$, where $C_{15}$ is a real constant,
we conclude that there do not exist parameter values
such that all $\varepsilon^k$-order terms in
$V_{9}(\varepsilon)$, $V_{10}(\varepsilon)$ and
$V_{11}(\varepsilon)$ vanish.

\vspace{0.10in}
{\bf (ii-3-2)} Since $q_{02}=0$ we assume that
$b_{21}=\frac{1}{3}p_{12}p_{21}(2p_{21}-1)$ and
$1+p_{21}\neq0$ based on $V_{5d}=0$.
If $p_{12}=0$, we have $V_5(\varepsilon)=0$,
which is included in the condition II.  Hence, we assume that
$p_{12}\neq0$. Similar to the analysis for the case (ii-3-1),
we consider the two subcases $b_{12}=-a_{21}$ and $M_3=0$.

\vspace{0.10in}
{\bf (ii-3-2-1)} Under the condition $b_{12}=-a_{21}$, we assume that
$V_{5b}=V_{5c}=0$, which yields $V_5(\varepsilon)=V_6(\varepsilon)=0$.
Then, setting the $\varepsilon^{17}$-order term in $V_7(\varepsilon)$
zero we obtain $5-16p_{21}+4p_{21}^2=0$, under which
the polynomials $V_{5b}$ and $V_{5c}$ are simplified to
\begin{equation*}
\begin{aligned}
\widehat{V}_{5b}=\ &-27 a_{03} - 324 a_{21}^2 + 940 a_{21} p_{12}^2
- 3 (153 a_{03} + 54 a_{21}^2 + 916 a_{21} p_{12}^2) p_{21},\\[0.5ex]
\widehat{V}_{5c}=\ &\
243 a_{21} + 54 p_{03} + 2495 p_{12}^2
+ 6 (378 a_{21} + 153 p_{03} - 1249 p_{12}^2) p_{21}.
\end{aligned}
\end{equation*}
It is easy to find from $\widehat{V}_{5b}=\widehat{V}_{5c}=0$ that
\begin{equation*}
\begin{aligned}
a_{03}=&-\frac{2 a_{21}[162 a_{21}- 470p_{12}^2
+ 3 (27 a_{21} + 458 p_{12}^2) p_{21}]}{27 (1 + 17 p_{21})} <0,
\ \ (\textrm{due to \eqref{Eqn27}}) \\
p_{03}=&- \frac{243 a_{21} + 2495 p_{12}^2
+ 6 (378 a_{21} - 1249 p_{12}^2) p_{21}}{54 (1 + 17 p_{21})}.
\end{aligned}
\end{equation*}
Further, we have the following simplified polynomials:
\begin{equation*}
\begin{aligned}
\widehat{V}_{7a}=&\ 1458(2034p_{21}-695)a_{21}^2-15(153111a_{21}
-25780708p_{12}^2)p_{12}^2\\
&+2(3360789 a_{21} - 565883302 p_{12}^2)p_{12}^2 p_{21},\\
\widehat{V}_{7b}=&\ 2916(16169p_{21}-5525)a_{21}^2+ 5 (598986 a_{21}
- 278284327 p_{12}^2 ) p_{12}^2
\\&-18(487243 a_{21}-226233559p_{12}^2)p_{12}^2p_{21}.
\end{aligned}
\end{equation*}
Then, we compute the resultant of $\widehat{V}_{7a}$ and $\widehat{V}_{7b}$
with respect to $p_{21}$ to obtain
\begin{equation*}
\begin{aligned}
{\rm Res}&[\widehat{V}_{7a},\widehat{V}_{7b},p_{12}]=\ 871696100250000 a_{21}^8
(1979853717341274735024025\\
&-23176897174962961258321510 p_{21}+101743987046015667217950816 p_{21}^2\\
&-198508767310489076670643400 p_{21}^3
+145238548628125887939760784 p_{21}^4)^2\neq0,
\end{aligned}
\end{equation*}
under the conditions: $5-16p_{21}+4p_{21}^2=0$
and $a_{21}\ne 0$. This shows that there are no parameter values
satisfying $\widehat{V}_{7a}=\widehat{V}_{7b}=0$.

\vspace{0.10in}
{\bf (ii-3-2-2)} Assume $b_{12}\neq-a_{21}$.
Letting the $\varepsilon^{17}$-order
term in $V_7(\varepsilon)$ be zero we have $5-16p_{21}+4p_{21}^2=0$,
yielding the following simplified polynomials:
\begin{equation*}
\begin{aligned}
\overline{V}_{5a}=&\ (5+28p_{21})a_{03}
+ 6 (2 a_{21}-3 b_{12}) a_{21},
\\[0.5ex]
\overline{V}_{5c}=&-1944a_{21}+3726b_{12}-2295p_{03}+18735p_{12}^2
\\
&+(7452p_{03}+6156a_{21}-12474b_{12}-54962p_{12}^2)p_{21}.
\end{aligned}
\end{equation*}
Setting $\overline{V}_{5a}=\overline{V}_{5c}=0$ gives
\begin{equation*}
\begin{aligned}
a_{03}=&-\frac{6 a_{21} (2a_{21}-3b_{12})}{5+28p_{21}},\\[0.5ex]
p_{03}=&\frac{1944a_{21}-3726b_{12}-18735p_{12}^2
+(12474b_{12}-6156a_{21}) p_{21}+54962p_{12}^2}{27(276p_{21}-85)}.
\end{aligned}
\end{equation*}
We use the above solutions to simplify $\overline{V}_{5b}$
and $V_7(\varepsilon)$ to obtain the following polynomials:
\begin{equation*}
\begin{aligned}
\overline{V}_{5b}=& \
324 a_{21}^2 (6222 p_{21} \!-\! 2185) \!+\!
 2 a_{21} \big[p_{12}^2 (86250505 \!-\! 252436996 p_{21}) \!+\!
    810 b_{12} (12636 p_{21} \\
    &-\! 4345) \big]- 80 b_{12} \big[p_{12}^2 (6373795 - 18653564 p_{21}) +
    810 b_{12} (1204 p_{21}-411)\big],
\\[0.5ex]
\overline{V}_{7a}=& 1296 a_{21}^4 (52553294 p_{21}-18022295 ) -
 27920 b_{12}^3 \big[p_{12}^2 (684743715 - 2004000728 p_{21})\\
&+4050 b_{12} (25872 p_{21}-8839)\big]
+4 a_{21}^3 \big[p_{12}^2 (1599113005145 - 4680077698534 p_{21})\\
&+81 b_{12} (2484285362 p_{21} \!-\! 850029185)\big]
+4a_{21}^2 b_{12}\big[81b_{12}(1225559595 \!-\! 3592556264 p_{21}) \\
&+p_{12}^2 (5277460809917 p_{21}-1803276434885)\big]
+2a_{21} b_{12}^2 \big[810 b_{12} (918918755 \\
& - 2690204032 p_{21}) - p_{12}^2 (14708042417435
- 43045101045752 p_{21})\big],\\
\overline{V}_{7b}=& \
729 a_{21}^3 (24810684334 p_{21}-8478232095)
-80 b_{12} \big[9 b_{12} p_{12}^2 (8334257322035\\
&-24391452032792p_{21})+10p_{12}^4(1107374710345-3240898045749 p_{21})\\
&+3645 b_{12}^2 (3097706668 p_{21}-1058445047 )\big]
+9 a_{21}^2 \big[p_{12}^2 (240122211943205\\
&-702753773748111 p_{21}) +486 b_{12}(61045680279 p_{21}-20858985470)\big]\\
&+a_{21}\big[5 p_{12}^4 (111885081254445 - 327448472749544 p_{21})
+7290 b_{12}^2 (30073269035\\
&-88014820004 p_{21})-9b_{12} p_{12}^2 (503643340635395
- 1473987501858484 p_{21})\big],\\
\end{aligned}
\end{equation*}
\begin{equation*}
\begin{aligned}
\overline{V}_{7c}=& \ 9234 a_{21}^2(333894386 p_{21}-114095455)
+3 a_{21}\big[p_{12}^2 (414050132594835\\
&-1211779848284032 p_{21})-810 b_{12}(23508039475-68799585032 p_{21})\big]\\
&+40\big[350 p_{12}^4(15370001625-44982613958 p_{21})
+2430 b_{12}^2 (1831770089\\
&-5360953016 p_{21})-3b_{12}p_{12}^2(31862998691465
-93251840991288p_{21})\big],\\
\overline{V}_{7d}=& 10 p_{12}^2 (43443 - 127132 p_{21})
+ 648 b_{12} (261 - 764 p_{21}) - 81 a_{21} (707 - 2068 p_{21}). \\
\end{aligned}
\end{equation*}
By computing the Groebner basis for $5-16p_{21}+4p_{21}^2$,
$\overline{V}_{5b}$,
$\overline{V}_{7a}$, $\overline{V}_{7b}$, $\overline{V}_{7c}$
and $\overline{V}_{7d}$, we obtain two polynomial equations:
\begin{equation*}
\begin{aligned}
p_{12}^2 (1812510 b_{12} + 13187705 p_{12}^2 - 11105856 p_{12}^2 p_{21})&=0,\\
8667 a_{21} + 12312 b_{12} - 4030 p_{12}^2 - 10368 b_{12} p_{21}
- 16880 p_{12}^2 p_{21}&=0,
\end{aligned}
\end{equation*}
which yields
\begin{equation*}
\begin{aligned}
&a_{21}=\frac{2}{290907855}\, p_{12}^2
(1066162176 p_{21}^2- 2248798664 p_{21} +1571031845),\\
&b_{12}=\frac{1}{1812510}\, p_{12}^2 (11105856 p_{21}-13187705),
\end{aligned}
\end{equation*}
under which $V_8(\varepsilon)=0$, and the $9$th generalized Lyapunov
constant is simplified as
\begin{equation*}
\begin{aligned}
V_{9}(\varepsilon)=&\
\frac{\pi p_{12}\varepsilon^{13}}{C_{15}(5 + 28 p_{21})^3 (276 p_{21}-85)^3}\\
&\times \big[ 8 p_{12}^{10}
(479401837151164960733153380166419352018673132763343471 p_{21} \\
&\qquad \qquad -163805665483591789907309292300795181773767007636854130) \\
&\quad\ + C_{16}\, p_{12}^8 (
35989441277122180465488301421620623564108062636 p_{21} \\
&\qquad \qquad -12297145988872641623638918735129812334369100955 )\varepsilon^2\\
&\quad\ + C_{17} p_{12}^6(
28047246207962034726462002167497642752864 p_{21} \\
&\qquad \qquad -9583396378659769296319246478020611770295) \varepsilon^4\\
&\quad\ + C_{18} p_{12}^4 (
317522879623822832097089379861158 p_{21} \\
&\qquad \qquad -108493632214946857667054466000365 )\varepsilon^6\\
&\quad \ + C_{19} p_{12}^2 (1459616804250815596262396 p_{21}
-498732969803156525142505)\varepsilon^8 \\
&\quad \ +(2716325631250436 p_{21} -928134798864155 )
\varepsilon^{10}\big],
\end{aligned}
\end{equation*}where $C_i$, $i=15,16,\cdots,19$, are {\color{red}integers.}
Hence, any $\varepsilon^k$-order term in
$V_{9}(\varepsilon)$ is non-zero
under the condition $5-16p_{21}+4p_{21}^2=0$ and $p_{12}\ne 0$.
\end{enumerate}

We have shown that the four
conditions I, II, III and IV in Theorem \ref{Thm2.2} are necessary for the
singular points $(\pm1,0)$ of system \eqref{Eqn9} to be nilpotent
bi-center. Now, we prove that these
four conditions are also sufficient.

If the condition I in Theorem \ref{Thm2.2} holds, system \eqref{Eqn9}
is reduced to
\begin{equation}\label{Eqn31}
\left(\begin{array}{cc}\dot{x}\\
\dot{y}
\end{array}\right)
=\left\{\begin{aligned}&\left(\begin{aligned}&2a_{21}xy+a_{21}x^2y
-3b_{03}xy^2+a_{03}y^3, \\
&x+\tfrac{3}{2}x^2-a_{21}y^2+\tfrac{1}{2}x^3-a_{21}xy^2+b_{03}y^3,
\end{aligned}\right),&\text{if} \ \ y>0, \\[1.0ex]
&\left(\begin{aligned}&2a_{21}xy+a_{21}x^2y-6b_{03}y^2-3b_{03}xy^2
+a_{03}y^3, \\
&x+\tfrac{3}{2}x^2-a_{21}y^2+\tfrac{1}{2}x^3-a_{21}xy^2+b_{03}y^3,
\end{aligned}\right), &\text{if} \ \ y<0,
\end{aligned}\right.
\end{equation}
via $x\rightarrow x+1$. The upper and the lower systems in \eqref{Eqn31}
are Hamiltonian systems, having respectively the Hamiltonian functions,
\begin{equation}\label{Eqn32}
\begin{aligned}
H^+(x,y)&=-\frac{1}{2}x^2-\frac{1}{2}x^3-\frac{1}{8}x^4
+\frac{a_{21}}{2}x^2y^2-b_{03}xy^3+a_{21}xy^2+\frac{a_{03}}{4}y^4,\\[1.0ex]
H^-(x,y)&=-\frac{1}{2}x^2-\frac{1}{2}x^3-\frac{1}{8}x^4+\frac{a_{21}}{2}x^2y^2
-b_{03}xy^3+a_{21}xy^2-2 b_{03}y^3+\frac{a_{03}}{4}y^4,
\end{aligned}
\end{equation}
which shows that the condition, $H^+(x,0)\equiv H^-(x,0)$,
in Proposition 2.1 of \cite{X.W3} is satisfied, and so the origin
of system \eqref{Eqn31} is a center. Hence, system \eqref{Eqn9}
has a nilpotent bi-center at ($\pm1$,0), which consists of
a monodromic singular point and a cusp.

If the condition II in Theorem \ref{Thm2.2}
holds, the system \eqref{Eqn9} is simplified into a smooth one,
\begin{equation}\label{Eqn33}
\left(\begin{array}{cc}\dot{x}\\
\dot{y}
\end{array}\right)
=\left(\begin{aligned}&-a_{21}y+a_{21}x^2y+a_{03} y^3,\\[0.0ex]
&-\tfrac{x}{2}+\tfrac{x^3}{2}+b_{12}xy^2,
\end{aligned}\right),
\end{equation}
which is symmetric with respect to the $x$-axis,
so we know that
the singular points {\rm($\pm1$,0)} of system
\eqref{Eqn9} are bi-center.

If the condition III in Theorem \ref{Thm2.2} holds, system
\eqref{Eqn9} again becomes a smooth one:
\begin{equation}\label{Eqn34}
\left(\begin{array}{cc}\dot{x}\\
\dot{y}
\end{array}\right)
=\left(\begin{aligned}&-a_{21}y+a_{21}x^2y+a_{03}y^3,\\
&-\tfrac{x}{2}+\tfrac{x^3}{2}-b_{21}y+b_{21}x^2y-a_{21}xy^2
-\tfrac{2}{3}a_{21}b_{21}y^3,
\end{aligned}\right),
\end{equation}
which has the algebraic integral curve,
\begin{equation*}
\begin{aligned}
I(x,y)=\,& 3(9a_{03}+2a_{21}^2)
(x^2 + 2 b_{21}xy - 2 a_{21}y^2) \\
& -(3a_{03}-2a_{21}^2)
(3 x^4 + 6 b_{21}x^3y - 12 a_{21} x^2y^2
-4 a_{21} b_{21}xy^3 - 6a_{03}y^4),
\end{aligned}
\end{equation*}
giving an inverse integrating factor $I(x,y)$ for this system.
Thus, system \eqref{Eqn9}
has a nilpotent bi-center at ($\pm 1$,0).

If the condition IV in Theorem \ref{Thm2.2} holds, system
\eqref{Eqn9} is a smooth one:
\begin{equation}\label{Eqn35}
\left(\begin{array}{cc}\dot{x}\\
\dot{y}
\end{array}\right)
=\left(\begin{aligned}
&\tfrac{3}{8}b_{21}^2y-\tfrac{3}{8}b_{21}^2x^2y+a_{03}y^3, \\
&-\tfrac{x}{2}-b_{21}y+\tfrac{x^3}{2}+b_{21}x^2y
+\tfrac{16a_{03}-3b_{21}^4}{12b_{21}^2}xy^2
-\tfrac{32a_{03}-9b_{21}^4}{24b_{21}}y^3,
\end{aligned}\right).
\end{equation}
Actually, for an arbitrary nonzero constant $r$,
by the transformation,
\begin{equation*}
x=X,\qquad y=rY,\qquad t=rT,
\end{equation*}system \eqref{Eqn35} can be changed to
\begin{equation}\label{Eqn36}
\left(\! \begin{array}{cc}\dot{X}\\
\dot{Y}
\end{array} \! \right)
\!=\! \left(\begin{aligned}
&\tfrac{3}{8}(b_{21}r)^2Y-\tfrac{3}{8}(b_{21}r)^2X^2Y+a_{03}r^4Y^3 \\[0.5ex]
&-\tfrac{X}{2}-b_{21}rY+\tfrac{X^3}{2}+b_{21}rX^2Y
+\tfrac{16a_{03}r^4-3(b_{21}r)^4}{12(b_{21}r)^2}XY^2
-\tfrac{32a_{03}r^4-9(b_{21}r)^4}{24(b_{21}r)}Y^3
\end{aligned}\right)\!,
\end{equation}
which implies that system \eqref{Eqn35} is invariant under the following
transformation $(a_{03},b_{21})\rightarrow(a_{03}r^4,b_{21}r)$.
Hence, we can always choose proper $r$ to satisfy
$\mathfrak{f}_3^+=2a_{03}-\frac{3}{16}b_{21}^4=-2$,
yielding $a_{03}=\frac{3}{32}b_{21}^2-1$. Then, the origin of
system \eqref{Eqn35} is a center according to the condition
$C_3$ of Theorem 4.2 in \cite{Li2}.

\begin{figure}[!htb]
\vspace{0.00in}
\centering
\hspace*{-0.0in}
\subfigure[]{
\label{Fig3a}
\includegraphics[width=0.26\textwidth,height=0.18\textheight]{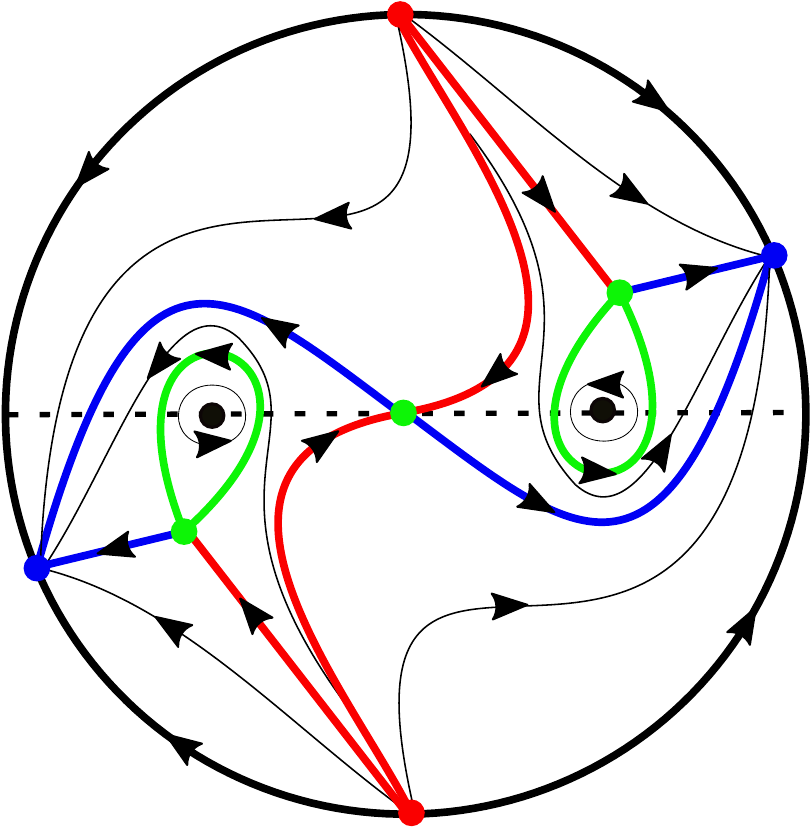}}
\hspace*{0.50in}
\subfigure[]{
\label{Fig3b}
\includegraphics[width=0.26\textwidth,height=0.18\textheight]{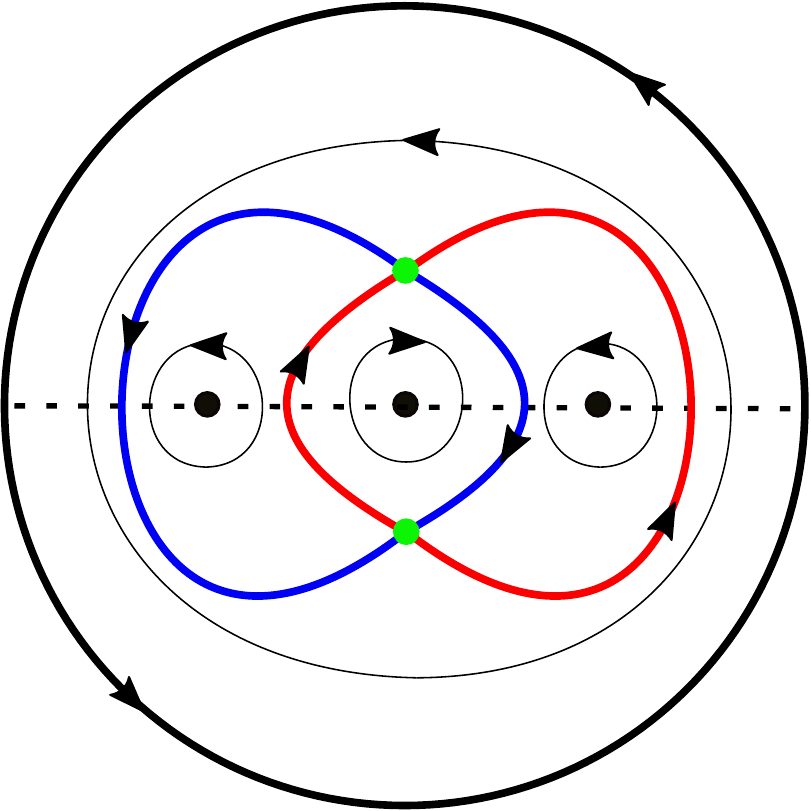}}
\hspace*{0.25in}

\hspace*{0.3in}
\subfigure[]{
\label{Fig3c}
\includegraphics[width=0.26\textwidth,height=0.18\textheight]{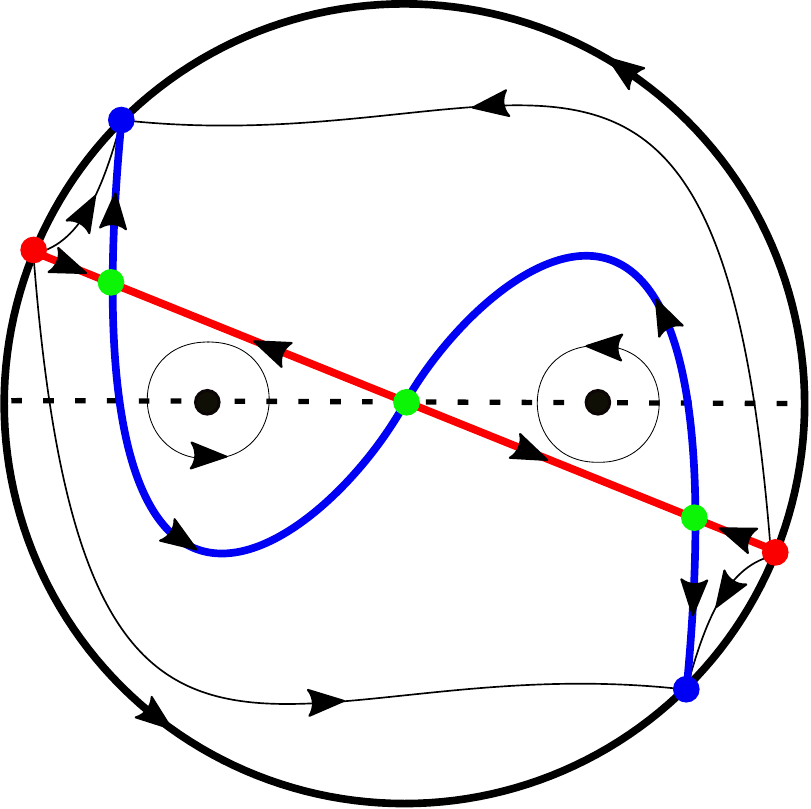}}
\hspace*{0.50in}
\subfigure[]{
\label{Fig3d}
\includegraphics[width=0.26\textwidth,height=0.18\textheight]{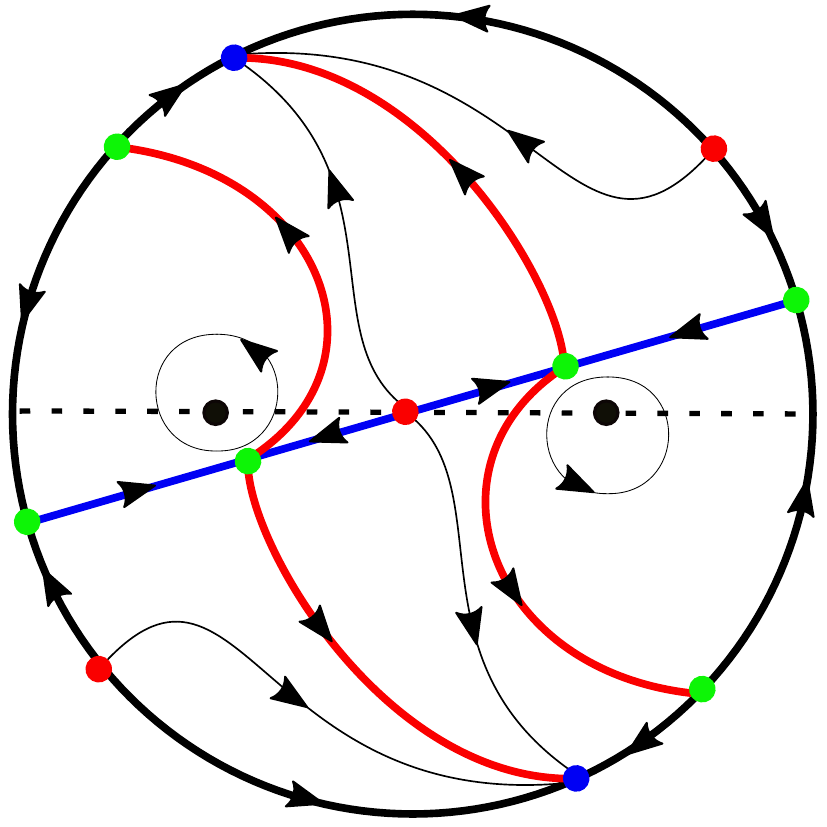}}
\hspace*{0.55in}
\caption{The phase portraits of system \eqref{Eqn9}
showing bi-center at $(\pm 1, 0)$
for (a) Condition I: $a_{21}=1,\, a_{12}=3,\,a_{03}=-3,\, b_{02}=b_{21}=0,\,
a_{02}=b_{12}=b_{03}=-1$; (b) Condition II:
$a_{21}=a_{03}=-1,\, a_{02}=a_{12}=b_{02}=b_{21}=b_{12}=b_{03}=0$;
(c) Condition III:
$a_{21}=1,\,a_{03}=-4,\, a_{02}=a_{12}=b_{02}=0,\,
b_{21}=\frac{3\sqrt{5}}{10},\,
b_{12}=-1, \,b_{03}=- \frac{\sqrt{5}}{5}$;
(d) Condition IV:
$a_{21}=-\frac{3}{2},\,a_{03}=-3,\, a_{02}=a_{12}=b_{02}=0,\, b_{21}=-2,\,
b_{12}=-2, \,b_{03}=5$.}
\label{Fig3}
\end{figure}

\begin{example}\label{Exam4.2}
The global phase portraits of system \eqref{Eqn9}
corresponding to the four bi-center conditions I, II, III and IV,
with respectively three sets of parameter values, show the bi-center at
($\pm 1$,0), as illustrated in Figure \ref{Fig3}.
\end{example}

\subsection{The $2$nd-order critical point $(1,0)$ of the upper system
in \eqref{Eqn9}}

In this subsection,
we consider the center conditions associated with $(1,0)$ of the
first system of \eqref{Eqn9} with multiplicity two.
Thus, assume that $\mathfrak{f}_2^+<0$,
i.e., $a_{02}+a_{12}<0$. Then, the singular point $(1,0)$ in the
upper smooth system of \eqref{Eqn9} is a cusp.
If $\mathfrak{f}_2^-=a_{12}-a_{02}=0$, then $(1,0)$ in the lower system
of \eqref{Eqn9} is a 3rd-order critical point. If
$\mathfrak{f}_2^-=a_{12}-a_{02}\neq0$, $(1,0)$
in the lower smooth system of \eqref{Eqn9} is also a cusp.
Similarly, recall that the singular point ($1,0$) cannot be
a monodromic singular point when $a_{02}-a_{12}>0$, and
hence we only consider the case when $a_{02}-a_{12} \le 0$.
Therefore, combining $a_{02}+a_{12}<0$ and
$a_{02}-a_{12} \le 0$ we have either $a_{02}=a_{12}<0$
or $a_{02}+ |a_{12}|<0$.

\begin{figure}[!htb]
\vspace{0.00in}
\centering
\includegraphics[width=0.26\textwidth,height=0.18\textheight]{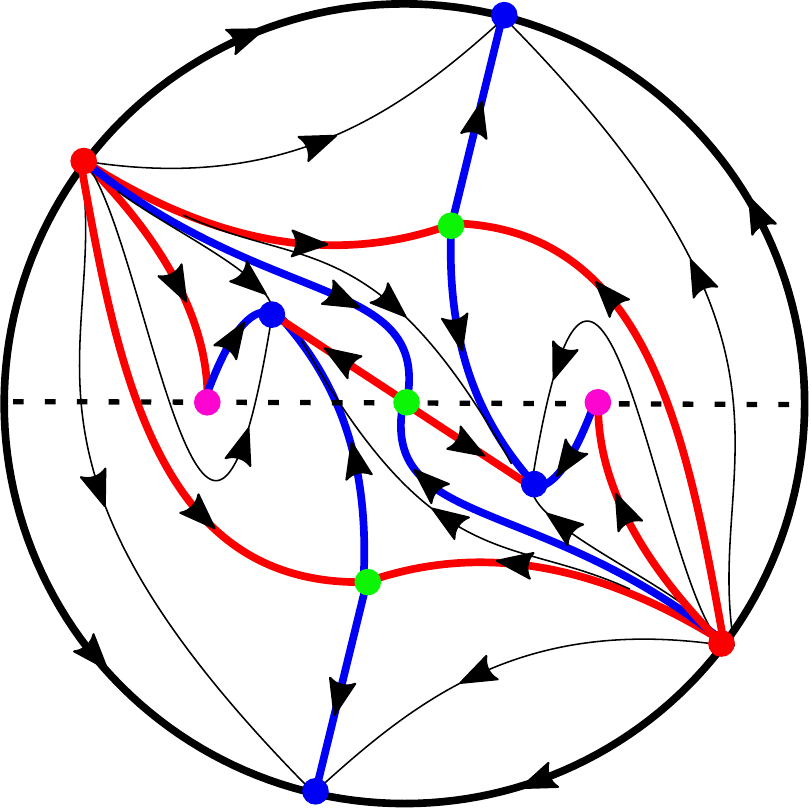}
\caption{The phase portrait of system \eqref{Eqn9}
with $a_{02}=2$, $a_{12}=-3$
and $a_{03}=a_{21}=b_{02}=b_{12}=b_{21}=b_{03}=1$, showing
two cusps at $(\pm 1, 0)$.}
\label{Fig4}
\end{figure}

\begin{example}\label{Exam4.3}
The global phase portrait of system
\eqref{Eqn9} with $a_{02}=2$, $a_{12}=-3$ and $a_{03}=a_{21}=b_{02}=b_{12}
=b_{21}=b_{03}=1$ shows that the nilpotent singular points
$(\pm1,0)$ are two cusps, see Figure \ref{Fig4}.
\end{example}

We apply the generalized Poincar\'e-Lyapunov method again to
system \eqref{Eqn9} with either $a_{02}=a_{12}<0$
or $a_{02}+ |a_{12}|<0$. Introducing the transformation
$(x,y,t)\rightarrow (\varepsilon^3(x+1),\varepsilon^2 y,
\frac{t}{\varepsilon})$ into
system \eqref{Eqn9}, we obtain the perturbed system,
\begin{equation}\label{Eqn37}
\left(\begin{array}{cc}\dot{x}\\
\dot{y}
\end{array}\right)
=\left\{\begin{aligned}&\left(\begin{aligned}
&-y+2 (a_{21}+p_{21}\varepsilon^2)\varepsilon xy
+(a_{02}+a_{12}+p_{02}\varepsilon^2)y^2\\
&+(a_{21}+p_{21}\varepsilon^2)\varepsilon^4 x^2y
+(a_{12}+p_{12}\varepsilon^2)\varepsilon^3 xy^2\\
&+(a_{03}+p_{03}\varepsilon^2)\varepsilon^2 y^3,\\[1.0ex]
&x+\tfrac{3}{2}\varepsilon^3 x^2+\tfrac{1}{2}\varepsilon^6 x^3
+(b_{02}+b_{12}+q_{02}\varepsilon^2)\varepsilon y^2\\
&+2(b_{21}+q_{21}\varepsilon^2)\varepsilon^2 xy
+(b_{21}+\varepsilon^2)\varepsilon^5 x^2y\\
&+(b_{12}+q_{12}\varepsilon^2)\varepsilon^4 xy^2
+(b_{03}+q_{03}\varepsilon^2)\varepsilon^3 y^3,
\end{aligned}\right), &\text{if} \ \ y>0,\\[1.0ex]
&\left(\begin{aligned}&-y+2(a_{21}+p_{21}\varepsilon^2)\varepsilon xy
+(a_{21}+p_{21}\varepsilon^2)\varepsilon^4 x^2y\\
&-(a_{02}-a_{12}+p_{02}\varepsilon^2-2p_{12}\varepsilon^2)y^2\\
&+(a_{12}+p_{12}\varepsilon^2)\varepsilon^3 xy^2
+ (a_{03}+p_{03}\varepsilon^2)\varepsilon^2 y^3,\\[1.0ex]
&x+\tfrac{3}{2}\varepsilon x^2+\tfrac{1}{2}\varepsilon^6 x^3
-(b_{02}-b_{12}+q_{02}\varepsilon^2-2q_{12}\varepsilon^2)\varepsilon y^2\\
&+2(b_{21}+q_{21}\varepsilon^2)\varepsilon^2xy
+(b_{21}+q_{21}\varepsilon^2)\varepsilon^5 x^2y\\
&+(b_{12}+q_{12}\varepsilon^2)\varepsilon^4xy^2
+(b_{03}+q_{03}\varepsilon^2)\varepsilon^3y^3,
\end{aligned}\right), &\text{if} \ \ y<0.
\end{aligned}
\right.
\end{equation}
which satisfies the condition \eqref{Eqn29} for the perturbed parameters.

Similarly, we first present a table to outline the proof for
the center conditions V and VI.

\begin{table}[!h]\label{T2}
\caption{Outline of the proof for the center conditions V-VI.}
{\scriptsize\begin{center}
\begin{tabular}{|l|l|l|l|c|}
\hline
\multicolumn{4}{|c|}{Cases}  &Conditions  \\
\hline
\multirow{5}*{(i) $a_{02}+|a_{02}|<0$} &\multirow{4}*{(i-1) $a_{12}=0$}
&\multirow{2}*{(i-1-1) $q_{21}=0$}  &(i-1-1-1) $p_{12}=\tfrac{b_{21}}{2}$ &  V \\[1mm]
\cline{4-5}
& & &(i-1-1-2) $p_{12}\ne\tfrac{b_{21}}{2}$ &  V \\[1mm]
\cline{3-5}
& & \multirow{2}*{(i-1-2) $q_{12}=-1$} &(i-1-2-1) $b_{12}=-a_{21}$ & V\\[1mm]
\cline{4-5}
& & &(i-1-2-2) $p_{12}=\tfrac{b_{21}}{2}$ &  --- \\[1mm]
\cline{2-5}
&\multirow{2}*{(i-2) $a_{12}\ne0$} &\multicolumn{2}{c|}{(i-2-1)
$q_{21}=0$} &VI\\[1mm]
\cline{3-5}
& & \multicolumn{2}{c|}{(i-2-2) $q_{12}=-1$}  & VI\\[1mm]
\hline
\multicolumn{4}{|c|}{(ii) $a_{02}=a_{12}<0$}  &VI \\[1mm]
\hline
\end{tabular}
\end{center}}
\end{table}

The first two generalized Lyapunov constants at the origin
of system \eqref{Eqn37} are
\begin{equation*}
V_1(\varepsilon)=0 \quad \textrm{and} \quad
V_2(\varepsilon)=\frac{8}{3}\varepsilon
\big[ b_{02}+(q_{02}-q_{12})\varepsilon^2 \big].
\end{equation*}
Setting $V_2(\varepsilon)=0$ we get the necessary center
condition: $b_{02}=0$, $q_{02}=q_{12}$, for system \eqref{Eqn37}.
Then, it follows from the discussion given at the beginning
of this subsection that we consider two cases:
(i) $a_{02}+ |a_{12}|<0$ and (ii) $a_{02}=a_{12}<0$.

\begin{enumerate}
\item[{\bf (i)}] By the $\varepsilon$-order term of $V_3(\varepsilon)$,
we have two subcases (i-1) $a_{12}=0$ and (i-2) $a_{12}\neq0$.

\vspace{0.10in}
{\bf (i-1)} If $a_{12}=0$, we have $a_{02}<0$.
Then, the $3$rd generalized Lyapunov constant is given by
\begin{equation*}
\begin{aligned}
V_3(\varepsilon)=&
-\frac{\pi}{4}\varepsilon^3
\big\{ 3 b_{03}-2b_{12}b_{21}
+2(a_{21} + b_{12})p_{12}
-2(1+q_{12})q_{21}\varepsilon^4
\\
&\qquad \quad  +\big[ 3q_{03} - 2b_{21}(1 + q_{12})
+ (1 +2p_{21} + 2 p_{12} )p_{12} -2b_{12} q_{21} \big] \varepsilon^2
\big\}.\\
\end{aligned}
\end{equation*}
By setting the $\varepsilon^{7}$-order term in $V_3(\varepsilon)$ zero
we have the following two subcases: (i-1-1) $q_{21}=0$ and (i-1-2) $q_{12}=-1$.

\vspace{0.10in}
{\bf (i-1-1)} Assume that $q_{21}=0$. Letting the $\varepsilon^3$-order
and $\varepsilon^5$-order terms in $V_3(\varepsilon)=0$ zero we obtain
\begin{equation*}
\begin{aligned}
b_{03}=\frac{2}{3}\big[ b_{12}b_{21}- (a_{21}+b_{12})p_{12} \big],\quad
q_{03}=\frac{1}{3}\big[ 2 b_{21}- (1+2p_{21}) p_{12}
+2 (b_{21}- p_{12})q_{12} \big].
\end{aligned}
\end{equation*}
Then, we have the $4$th generalized Lyapunov constant,
\begin{equation*}
\begin{aligned}
V_4(\varepsilon)=&\frac{16}{45}\varepsilon^3\big[a_{02}
+(p_{02}-p_{12})\varepsilon^2\big]\big\{2(a_{21}+b_{12})(b_{21}
+2p_{12})
\\
& +\big[3b_{21}
+2p_{21}(b_{21}+2p_{12})+2q_{12}(b_{21}+2p_{12})\big]\varepsilon^2 \big\}.\\
\end{aligned}
\end{equation*}
Two subcases follow the $\varepsilon^{5}$-order term in $V_4(\varepsilon)$
to be zero: (i-1-1-1) $p_{12}=-\frac{b_{21}}{2}$,
and (i-1-1-2) $p_{12}\ne-\frac{b_{21}}{2}$.

\vspace{0.10in}
{\bf (i-1-1-1)} If $p_{12}=-\frac{b_{21}}{2}$, we obtain $b_{21}=0$
by $V_4(\varepsilon)=0$, giving the condition V.

\vspace{0.10in}
{\bf (i-1-1-2)} If $p_{12}\ne-\frac{b_{21}}{2}$,
we obtain $b_{12}=-a_{21}$ and
$$
p_{21}=\frac{-3b_{21}
-2 (b_{21}+2p_{12})q_{12}}{2(b_{21}+2p_{12})}$$
from $V_4(\varepsilon)=0$. Further, we have
\begin{equation*}
\begin{aligned}
V_5(\varepsilon)=&-\frac{\pi}{144(b_{21}+2p_{12})^3}\,
b_{21}\varepsilon^5\big[45a_{02}^2(b_{21}+2p_{12})^3
+2(b_{21}+2p_{12})^2 V_{5a}\varepsilon^2\\
&\hspace*{1.65in}
+(b_{21}+2p_{12})V_{5b}\varepsilon^4-18V_{5c}\varepsilon^6 \big],\\
\end{aligned}
\end{equation*}
where
\begin{equation*}
\begin{aligned}
V_{5a}=&\ 45 a_{02} (p_{02}-p_{12}) (b_{21}+2 p_{12})+9 a_{03} (2 b_{21}+p_{12})
\\
& +2 a_{21} (p_{12}-2 b_{21}) (9 a_{21}+4 b_{21}^2)+80 a_{21} b_{21} p_{12}^2,
\\[0.5ex]
V_{5b}=& \
b_{21}^2 (108 a_{21}+32 b_{21}^2+45 p_{02}^2)
+18 b_{21} (2 b_{21}+5 p_{12}) p_{03}
\\
& -2 b_{21} (9 a_{21}-28 b_{21}^2+45 b_{21} p_{02}-90 p_{02}^2) p_{12}
\\
& +3 \big[96 a_{21}-29 b_{21}^2-120 (b_{21}+p_{12}) p_{02}
+60 p_{02}^2+12 p_{03} \big] p_{12}^2
-212 b_{21} p_{12}^3
\\
& -140 p_{12}^4
+8 (2 p_{12}+b_{21}) \big[(2 b_{21}-p_{12}) (9 a_{21}+2 b_{21}^2)
-20 b_{21} p_{12}^2 \big] q_{12},
\\[0.5ex]
V_{5c}=& \ \big[ b_{21}-4p_{12}+2 (b_{21}+2p_{12}) q_{12} \big]\big[ 2 b_{21} (b_{21}+2 p_{12})
+2 (b_{21}^2-p_{12}^2) q_{12}\\
&+3 (2 p_{12}+b_{21} q_{12}) p_{12} \big].
\end{aligned}
\end{equation*}It follows from $V_5(\varepsilon)=0$ that $b_{21}=0$,
which is included in the condition V. Otherwise,
the $\varepsilon^5$-order term in $V_5(\varepsilon)$ is non-zero.

\vspace{0.10in}
{\bf (i-1-2)} Assume $q_{21}\ne 0$ and $q_{12}=-1$.
Letting the $\varepsilon^{3}$-order and $\varepsilon^5$-order terms
in $V_3(\varepsilon)$ zero, we obtain
\begin{equation*}
\begin{aligned}
b_{03}=\frac{2}{3} \big[b_{12}b_{21}- (a_{21}+b_{12})p_{12} \big],\quad
q_{03}=\frac{1}{3}(p_{12}-2 p_{12}p_{21}+2b_{12}q_{21}).
\end{aligned}
\end{equation*}
Then, we have the $4$th generalized Lyapunov constant,
\begin{equation*}
\begin{aligned}
V_4(\varepsilon)=&\frac{16}{45}\varepsilon^3\big[a_{02}
+(p_{02}-p_{12})\varepsilon^2\big]\big\{2(a_{21}+b_{12})
(b_{21}+2p_{12})
\\
& +\big[ b_{21}-4p_{12}+ 2 (b_{21}+2p_{12}) p_{21}
+2 (a_{21}+b_{12}) q_{21} \big] \varepsilon^2
+(1+2p_{21})q_{21}\varepsilon^4\big\}.\\
\end{aligned}
\end{equation*}
Letting the $\varepsilon^{9}$-order terms in $V_4(\varepsilon)$ zero
we have $p_{21}=-\frac{1}{2}$.

\vspace{0.10in}
{\bf (i-1-2-1)} If $b_{12}=-a_{21}$, the $\varepsilon^3$-order term
in $V_4(\varepsilon)$ becomes zero. By using $V_4(\varepsilon)=0$
we have $p_{12}=0$. Then, the $5$th generalized Lyapunov constant
is simplified as
\begin{equation*}
\begin{aligned}
V_5(\varepsilon)=&\ - \frac{\pi}{144}\varepsilon^5(b_{21}
+q_{21}\varepsilon^2)\big[45a_{02}^2+2(18a_{03}-36a_{21}^2-16a_{21}b_{21}^2
+45a_{02}p_{02})\varepsilon^2\\
&-(36a_{21}-45p_{02}^2-36p_{03}+64a_{21}b_{21}q_{21})\varepsilon^4
-32a_{21}q_{21}^2\varepsilon^6 \big].\\
\end{aligned}
\end{equation*}
If $b_{21}=0$, we obtain the condition included in the
condition V. Otherwise, the $\varepsilon^5$-order term in
$V_5(\varepsilon)$ does not equal zero.

\vspace{0.10in}
{\bf (i-1-2-2)} If $p_{12}=-\frac{b_{21}}{2}$, the $\varepsilon^3$-order term
in $V_4(\varepsilon)$ again becomes zero. From $V_4(\varepsilon)=0$
we obtain $b_{12}=\frac{-3b_{21}-2a_{21}q_{21}}{2q_{21}}$.
Then, we have the $5$th generalized Lyapunov constant,
\begin{equation*}
\begin{aligned}
V_5(\varepsilon)=&\frac{\pi}{576q_{21}^2}\varepsilon^5(b_{21}
+q_{21}\varepsilon^2)\big[\widetilde{V}_{5a}+18\widetilde{V}_{5b}\varepsilon^2
+q_{21}\widetilde{V}_{5c}\varepsilon^4+16q_{21}^3(13b_{21}
+8a_{21}q_{21})\varepsilon^6 \big],\\
\end{aligned}
\end{equation*}
where
\begin{equation*}
\begin{aligned}
\widetilde{V}_{5a}=&\ 36 \big[ 9a_{21}b_{21}^2
- b_{21} (3a_{03}-10a_{21}^2)q_{21}
-5a_{02}^2q_{21}^2 \big],\\[0.5ex]
\widetilde{V}_{5b}=&\ 9 b_{21}^2+2 b_{21} (9 b_{21}^2+19 a_{21}-3 p_{03}) q_{21}
\\
& -2 \big[4 a_{03}+5 (b_{21}+2 p_{02}) a_{02}
-8 a_{21} (a_{21}+b_{21}^2) \big] q_{21}^2 ,
\\[0.5ex]
\widetilde{V}_{5c}=&\
252 b_{21}
+9 \big[ 16 a_{21}+63 b_{21}^2-20 (b_{21}+p_{02}) p_{02}-16 p_{03} \big] q_{21}
+416 a_{21} b_{21} q_{21}^2.
\end{aligned}
\end{equation*}If $b_{21}=0$ we have $\widetilde{V}_{5a}=-180a_{02}^2q_{21}^2\ne0$.
Hence, by using $V_5(\varepsilon)=0$ we have
$a_{21}=-\frac{13}{8q_{21}}b_{21}\ne0$, $\widetilde{V}_{5b}=0$ and
\begin{equation*}
\begin{aligned}
a_{03}=&\ \frac{1}{96b_{21}q_{21}}(377b_{21}^3-160a_{02}^2q_{21}^2),\\
p_{03}=&\ \frac{1}{144q_{21}}\big[
b_{21} (18-109b_{21} q_{21})
-180 (b_{21} + p_{02}) p_{02} q_{21} \big].
\end{aligned}
\end{equation*}
Setting the $\varepsilon^{17}$-order term in $V_6(\varepsilon)$ zero
we have $p_{02}=-\frac{1}{2}b_{21}$. Further,
we have the simplified Lyapunov constant,
\begin{equation*}
\begin{aligned}
V_6(\varepsilon)=&\frac{-1}{9450b_{21}^2q_{21}^2}\varepsilon^7
(b_{21}+q_{21}\varepsilon^2)\big[96 a_{02} q_{21}
(7567b_{21}^3+296b_{21}^4q_{21}-2360a_{02}^2 q_{21}^3)\\
&-22050\pi b_{21}(16b_{21}^3+2b_{21}^4q_{21}-5a_{02}^2q_{21}^3)\varepsilon
+1152a_{02}b_{21}^2q_{21}^2(35-78b_{21}q_{21})\varepsilon^2\\
&-11025\pi q_{21}(16b_{21}^3+2b_{21}^4q_{21}-5a_{02}^2q_{21}^3)\varepsilon^3
-18432a_{02}b_{21}q_{21}^3(1+3b_{21}q_{21})\varepsilon^4\\
&-16384a_{02}b_{21}q_{21}^5\varepsilon^6 \big],
\end{aligned}
\end{equation*}
which has a non-zero $\varepsilon^{15}$-order term.

\vspace{0.10in}
{\bf (i-2)} If $a_{12}\neq0$, we obtain the $3$rd generalized Lyapunov constant,
\begin{equation*}
\begin{aligned}
V_3(\varepsilon)=\frac{\pi}{4}\varepsilon\big\{2a_{12}(a_{21}+b_{12})
+\big[ a_{12}(1 +2p_{21}) +3b_{03}-2b_{12}b_{21}+2 (a_{21}+b_{12}) p_{12}
\big] \varepsilon^2\\
+\big[3q_{03} - 2b_{21} (1 +q_{12} ) + (1+2p_{21} +2 q_{12}) p_{12}
-2b_{12}q_{21} \big] \varepsilon^4
-2(1+q_{12})q_{21}\varepsilon^6 \big\}.\\
\end{aligned}
\end{equation*}
Similar to the previous case, we have the following two subcases:
(i-2-1) $q_{21}=0$ and (i-2-2) $q_{12}=-1$.

\vspace{0.10in}
{\bf (i-2-1)} If $q_{21}=0$, setting the $\varepsilon$-order,
$\varepsilon^{3}$-order and $\varepsilon^5$-order terms in
$V_3(\varepsilon)$ to equal zero yields
\begin{equation*}
\begin{aligned}
b_{12}=&-a_{21},\\
b_{03}=&- \frac{1}{3}\big[
a_{12}(1 + 2 p_{21} + 2 q_{12}) +2a_{21}b_{21} \big],\\
q_{03}=&\frac{1}{3} \big[ 2b_{21} (1 + q_{12} )
- (1+2p_{21}+2 q_{12} )p_{12} \big].
\end{aligned}
\end{equation*}
Further, we have
\begin{equation*}
\begin{aligned}V_4(\varepsilon)=&\frac{16}{45}\varepsilon^3
\big[a_{02}+(p_{02}-p_{12})\varepsilon^2\big]
\big\{4a_{12}(p_{21}+q_{12})+ \big[3b_{21}
+2 (b_{21}+2p_{12}) (p_{21} +q_{12}) \big] \varepsilon^2\big\}.
\end{aligned}
\end{equation*}
By using $V_4(\varepsilon)=0$ we have $b_{21}=0$ and $q_{12}=-p_{21}$,
giving the first condition in VI.

\vspace{0.10in}
{\bf (i-2-2)} If $q_{12}=-1$, vanishing each $\varepsilon$ order term
in $V_3(\varepsilon)$ yields
\begin{equation*}
\begin{aligned}
b_{12}=-a_{21},\quad
b_{03}=\frac{1}{3}(a_{12}-2a_{21}b_{21}-2a_{12}p_{21}),\quad
q_{03}=\frac{1}{3}(p_{12}-2p_{12}p_{21}-2a_{21}q_{21}).
\end{aligned}
\end{equation*}
Then, we obtain the $4$th generalized Lyapunov constant,
\begin{equation*}
\begin{aligned}V_4(\varepsilon)=&\frac{16}{45}\varepsilon^3
\big[a_{02}+(p_{02}-p_{12})\varepsilon^2\big]
\\
& \times \big\{4a_{12}(p_{21}-1)
+ \big[ b_{21} (1+2p_{21})
+4(p_{21} -1)p_{12} \big] \varepsilon^2
-q_{21}(1+2p_{21})\varepsilon^4 \big\}.\\
\end{aligned}
\end{equation*}
It follows from $V_4(\varepsilon)=0$ that $b_{21}=0$, $p_{21}=1$
and $q_{21}=0$, which also leads to the first condition in VI.

\item[{\bf (ii)}]
If $a_{02}=a_{12}<0$, similar to the analysis for the subcase (i-2),
replacing $a_{02}$ by $a_{12}$ in $V_4(\varepsilon)$, then
we obtain the second condition in VI.
\end{enumerate}

The above discussions (or see Table 2) show that the conditions V and VI are necessary
for the origin of the perturbed system \eqref{Eqn37} to be a center,
and so they are necessary conditions for ($\pm1$,0) of the unperturbed
system \eqref{Eqn9} to be nilpotent centers.
We can also prove that these conditions are sufficient.

If the condition V in \eqref{Eqn15} holds, system
\eqref{Eqn9} is reduced to
\begin{equation}\label{Eqn38}
\left(\begin{array}{cc}\dot{x}\\
\dot{y}
\end{array}\right)
=\left\{\begin{aligned}&\left(\begin{aligned}&-a_{21}y+a_{21}x^2y
+a_{02} y^2+a_{03} y^3, \\[0.5ex]
&-\tfrac{x}{2}+\tfrac{x^3}{2}+b_{12}xy^2,
\end{aligned}\right), &\text{if} \ \ y>0,\\[1.0ex]
&\left(\begin{aligned}&-a_{21}y+a_{21}x^2y-a_{02} y^2+a_{03} y^3, \\[0.5ex]
&-\tfrac{x}{2}+\tfrac{x^3}{2}+b_{12}xy^2,
\end{aligned}\right), &\text{if} \ \ y<0.
\end{aligned}\right.
\end{equation}
It is easy to see that system \eqref{Eqn38} is symmetric with respect
to the $x$-axis, and so by the symmetry of switching systems redefined
in Theorem 2.1 of \cite{LY2015}, it implies that the singular points
{\rm($\pm1$,0)} of the system \eqref{Eqn9} are bi-center,
which consists of two 2nd-order nilpotent cusps.

If the condition VI in \eqref{Eqn15} holds, via $x\rightarrow x+1$,
system \eqref{Eqn9} becomes
\begin{equation}\label{Eqn39}
\left(\begin{array}{cc}\dot{x}\\
\dot{y}
\end{array}\right)
=\left\{\begin{aligned}&\left(\begin{aligned}&2a_{21}xy+(a_{02}+a_{12}) y^2
+a_{21}x^2y-3b_{03}xy^2+a_{03}y^3, \\[0.5ex]
&x+\tfrac{3}{2}x^2-a_{21}y^2+\tfrac{1}{2}x^3-a_{21}xy^2+b_{03}y^3,
\end{aligned}\right),&\text{if}\ \ y>0,
\\[1.0ex]
&\left(\begin{aligned}&2a_{21}xy+a_{21}x^2y-(a_{02}+a_{12}+6b_{03})y^2
-3b_{03}xy^2+a_{03}y^3, \\[0.5ex]
&x+\tfrac{3}{2}x^2-a_{21}y^2+\tfrac{1}{2}x^3-a_{21}xy^2+b_{03}y^3,
\end{aligned}\right),&\text{if} \ \ y<0.
\end{aligned}\right.
\end{equation}
Then it can be shown that the upper and lower systems in \eqref{Eqn39}
are Hamiltonian systems, having respectively the Hamiltonian functions:
\begin{equation}\label{Eqn40}
\begin{aligned}
H^+(x,y)&=-\frac{1}{2}x^2-\frac{1}{2}x^3
+a_{21}xy^2 +\frac{a_{02}+a_{12}}{3}y^3
-\frac{1}{8}x^4 +\frac{a_{21}}{2}x^2y^2-b_{03}xy^3
+\frac{a_{03}}{4}y^4, \\[1.0ex]
H^-(x,y)&=-\frac{1}{2}x^2 \!-\! \frac{1}{2}x^3
\!+\! a_{21}xy^2 \!-\! \frac{a_{02} \!+\! a_{12} \!+\! 6b_{03}}{3}y^3
\!-\! \frac{1}{8}x^4 \!+\! \frac{a_{21}}{2}x^2y^2 \!-\! b_{03}xy^3
\!+\! \frac{a_{03}}{4}y^4.
\end{aligned}
\end{equation}
It is straightforward to verify that the condition $H^+(x,0)\equiv H^-(x,0)$
in Proposition 2.1 of \cite{X.W3} is satisfied, indicating that
the origin of system \eqref{Eqn39} is a center. Hence, system
\eqref{Eqn9} has bi-center at ($\pm1$,0).

\begin{figure}[!htb]
\vspace{0.00in}
\centering
\subfigure[]{
\label{Fig5a}
\includegraphics[width=0.26\textwidth,height=0.18\textheight]{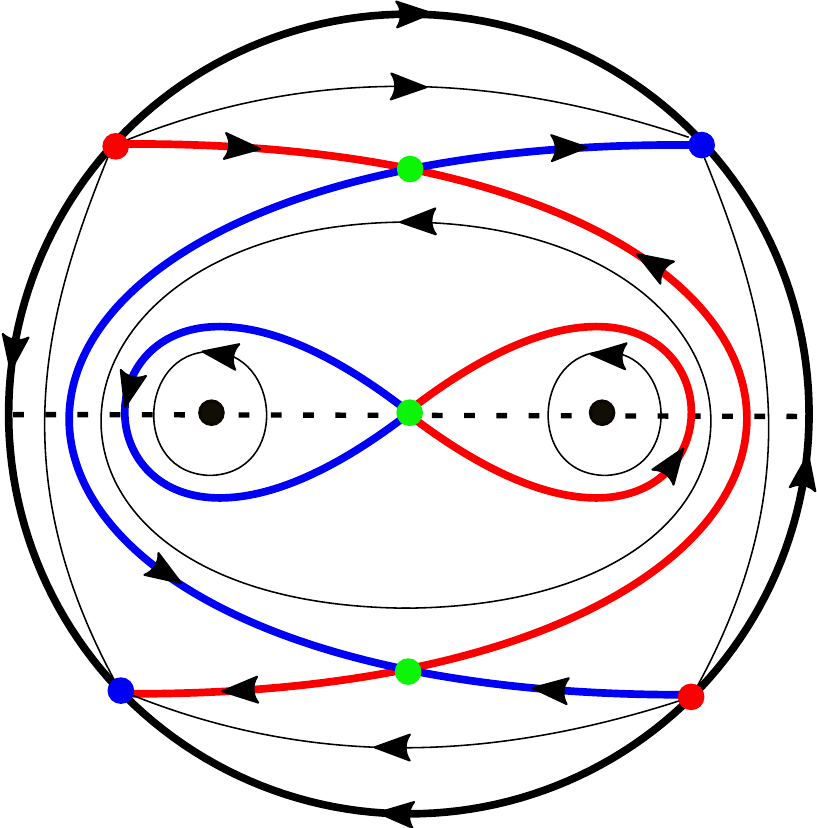}}
\hspace*{0.5in}
\subfigure[]{
\label{Fig5b}
\includegraphics[width=0.26\textwidth,height=0.18\textheight]{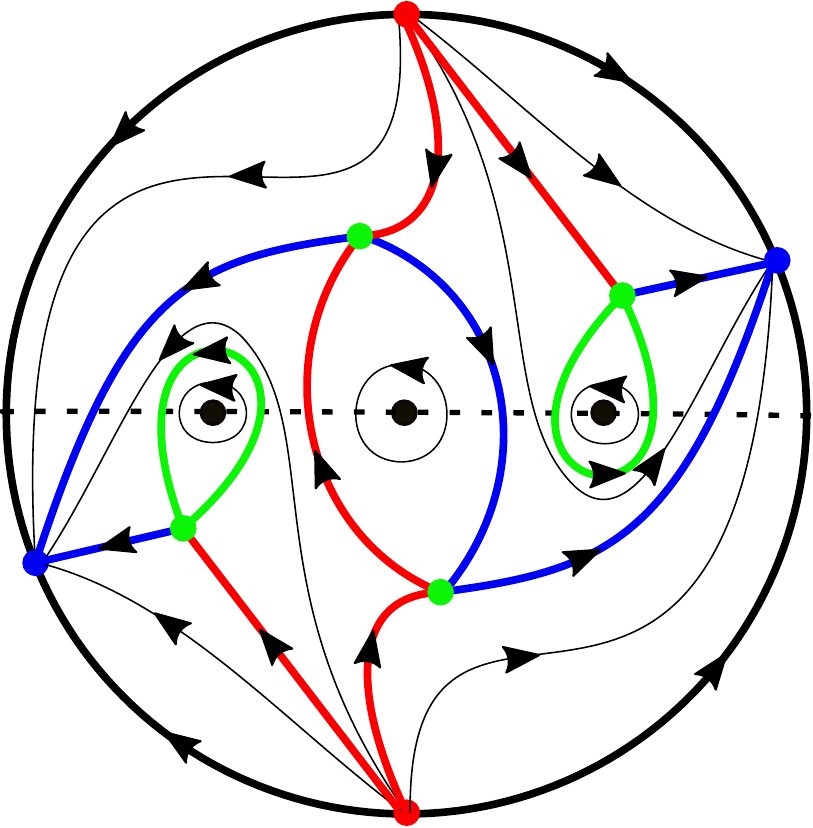}}
\caption{The phase portraits of system
\eqref{Eqn9} showing bi-center at $(\pm 1, 0)$
for (a) Condition V:
$a_{02}=-1,\, a_{12}=b_{03}=b_{02}=b_{21}=0,\, a_{21}=a_{03}=b_{12}=1$;
and (b) Condition VI:
$a_{02}=-4,\,a_{21}=b_{03}=-1,\, a_{12}=3,\, a_{03}=b_{12}=1,\,
b_{02}=b_{21}=0$.}
\label{Fig5}
\end{figure}

\begin{example}\label{Exam4.4}
The global phase portraits of system \eqref{Eqn9} corresponding to
the bi-center conditions V and VI,
with two sets of parameter values, as given in Figure \ref{Fig5},
show the case of bi-center at $(\pm 1, 0)$.
\end{example}

The proof for Theorem \ref{Thm2.2} is completed.

\section{The proof of Theorem \ref{Thm2.3}}

In this section, we will perturb system
\eqref{Eqn9} with the $6$ center conditions
in \eqref{Eqn15} to find the maximal number of small-amplitude
limit cycles bifurcating from the bi-center. Since system
\eqref{Eqn9} is symmetric with the origin, we only need to study
the limit cycles around the singular point $(1,0)$.
Using each of the $6$ center conditions, we can show that system \eqref{Eqn9}
can have $9$ small-amplitude limit cycles around each of the bi-center,
and then plus additional $1$ small limit cycle from the
pseudo-Hopf bifurcation, leading to a total of $20$ small-amplitude limit cycles bifurcating from the bi-center.
Perturbations up to cubic terms are applied to system \eqref{Eqn9},
and it will be shown that there exist maximal $9$ (including
both system parameters and perturbation parameters)
independent perturbation parameters. Therefore, $9$ small-amplitude
limit cycles may bifurcate from each of the bi-center.
Since the proofs for the $6$ center conditions I-VI are similar, we
only prove one case for the center condition VI.
It has been noted that for the first $5$ center conditions I-V,
when the conditions are obtained to satisfy
$V_1 = V_2 = \cdots = V_8 =0$, then $V_9$ and $V_{10}$ equal zero
simultaneously under a same condition, for which $V_{11} \ne 0$,
implying that maximal $9$ small-amplitude limit cycles many exist
around each of the bi-center. For the $6$th center condition VI,
solutions exist such that $V_1= V_2 = \cdots = V_9=0$, $V_{10} \ne 0$,
again indicating the existence of maximal $9$ small-amplitude limit cycles.

In order to include all possible perturbations, we add
the perturbations including all $\varepsilon^k$ ($k\ge 1$) order
terms to system \eqref{Eqn9}
under the center condition VI to obtain
\begin{equation}\label{Eqn41}
{\small\left(\begin{array}{cc}\dot{x}\\
\dot{y}
\end{array}\right)
=\left(\begin{aligned}
&\left(\begin{aligned}&-a_{21}y+a_{02}y^2+a_{21}x^2y-3b_{03}xy^2+a_{03}y^3
-\tfrac{\delta_1}{2}\varepsilon(x-x^3)\\
&+P^+(x,y,\varepsilon),\\[0.5ex]
&-\tfrac{1}{2}x+\tfrac{1}{2}x^3-a_{21}xy^2+b_{03}y^3
+\tfrac{b}{2}\varepsilon^3(x-x^2)+\delta_1\varepsilon y+Q^+(x,y,\varepsilon),
\end{aligned}\right),&\text{if} \ \ y>0, \\[1.0ex]
&\left(\begin{aligned}
&-a_{21}y-a_{02}y^2+a_{21}x^2y-3b_{03}xy^2+a_{03}y^3
-\tfrac{\delta_1}{2}\varepsilon(x-x^3)\\
&+P^-(x,y,\varepsilon),\\[0.5ex]
&-\tfrac{1}{2}x+\tfrac{1}{2}x^3-a_{21}xy^2+b_{03}y^3
+\tfrac{b}{2}\varepsilon^3(x+x^2)+\delta_1\varepsilon y+Q^-(x,y,\varepsilon),
\end{aligned}\right),&\text{if} \ \ y<0,
\end{aligned}\right.}
\end{equation}
where
\begin{equation}\label{Eqn42}
\begin{aligned}
&P^+(x,y,\varepsilon)=\sum_{k=1}\sum_{i+j=0}^3\varepsilon^k p_{ijk}x^iy^j,\qquad
P^-(x,y,\varepsilon)=\sum_{k=1}\sum_{i+j=0}^3\varepsilon^k(-1)^{i+j+1}p_{ijk}x^iy^j,\\
&Q^+(x,y,\varepsilon)=\sum_{k=1}\sum_{i+j=0}^3\varepsilon^kq_{ijk}x^iy^j,\qquad
Q^-(x,y,\varepsilon)=\sum_{k=1}\sum_{i+j=0}^3\varepsilon^k(-1)^{i+j+1}q_{ijk}x^iy^j.
\end{aligned}
\end{equation}
With the transformation
$ x \rightarrow x+1$, the first two polynomials in \eqref{Eqn42} become
\begin{equation}\label{Eqn43}
\begin{aligned}
p^+(x,y,\varepsilon)=&\,\sum_{k=1}\big[ p_{00k}+p_{10k}+p_{20k}+p_{30k}
+(p_{10k}+2p_{20k}+3p_{30k})x \\[-2.0ex]
&\qquad +(p_{20k}+3p_{30k})x^2+p_{30k}x^3+(p_{01k}+p_{11k}+p_{21k})y\\
&\qquad +(p_{11k}+2p_{21k})xy+p_{21k}x^2y+(p_{02k}+p_{12k})y^2
+p_{12k}xy^2+p_{03k}y^3\big]\varepsilon^k,\\[1.0ex]
p^-(x,y,\varepsilon)=&\sum_{k=1}\big[-p_{00k}+p_{10k}-p_{20k}+p_{30k}
+(p_{10k}-2p_{20k}+3p_{30k})x \\[-2.0ex]
&\qquad -(p_{20k}-3p_{30k})x^2 +p_{30k}x^3+(p_{01k}-p_{11k}+p_{21k})y\\
&\qquad -(p_{11k}-2p_{21k})xy+p_{21k}x^2y-(p_{02k}-p_{12k})y^2
+p_{12k}xy^2+p_{03k}y^3\big]\varepsilon^k.\\
\end{aligned}
\end{equation}
In order to keep the origin of the corresponding shifted system
to be a monodromic point, we eliminate the constants and the $x$
and $y$ terms in each $\varepsilon^k$-order term of $p^+(x,y,\varepsilon)$
and  $p^-(x,y,\varepsilon)$, and obtain
\begin{equation*}
\begin{aligned}
p_{00k}=p_{10k}=p_{20k}=p_{11k}=p_{30k}=0, \quad p_{01k}=-p_{21k}, \quad
\forall k \ge 1.
\end{aligned}
\end{equation*}
Similarly, simplifying $Q^{\pm}(x,y,\varepsilon)$
in \eqref{Eqn42} yields
\begin{equation*}
\begin{aligned}
q_{00k}=q_{10k}=q_{20k}=q_{11k}=q_{30k}=0, \quad q_{01k}=-q_{21k}, \quad
\forall k \ge 1.
\end{aligned}
\end{equation*}
Next, with a simple parametrization:
$p_{02k}+p_{12k}\rightarrow p_{02k}$ and
$q_{02k}+q_{12k}\rightarrow q_{02k}$, we have the following
four translational perturbations:
\begin{equation}\label{Eqn44}
\begin{aligned}
p^+(x,y,\varepsilon)=&\ \sum_{k=1}(p_{30k}x^3+2p_{21k}xy+p_{21k}x^2y
+p_{02k}y^2+p_{12k}xy^2+p_{03k}y^3)\varepsilon^k,\\
p^-(x,y,\varepsilon)=&\ \sum_{k=1}\big[p_{30k}x^3+2p_{21k}xy+p_{21k}x^2y
+(2p_{12k}-p_{02k})y^2+p_{12k}xy^2+p_{03k}y^3\big]\varepsilon^k,\\
q^+(x,y,\varepsilon)=&\ \sum_{k=1}(q_{30k}x^3+2q_{21k}xy+q_{21k}x^2y
+q_{02k}y^2+q_{12k}xy^2+q_{03k}y^3)\varepsilon^k,\\
q^-(x,y,\varepsilon)=&\ \sum_{k=1}\big[q_{30k}x^3+2q_{21k}xy+q_{21k}x^2y
+(2q_{12k}-q_{02k})y^2+q_{12k}xy^2+q_{03k}y^3]\varepsilon^k.\\
\end{aligned}
\end{equation}
Therefore, under the transformation $x \rightarrow x+1$,
system \eqref{Eqn41} becomes the following perturbed system,
\begin{equation}\label{Eqn45}
{\small\left(\begin{array}{cc}\dot{x}\\
\dot{y}
\end{array}\right)
=\left\{\begin{aligned}&\left(\begin{aligned}&2a_{21}xy
+(a_{02}-3b_{03})y^2+a_{21}x^2y-3b_{03}xy^2+a_{03}y^3\\
&+\delta_1 \varepsilon(x+\tfrac{3}{2}x^2+\tfrac{1}{2}x^3)-\varepsilon^2y
+p^+(x,y,\varepsilon),\\[0.5ex]
&x+\tfrac{3}{2}x^2+\tfrac{1}{2}x^3-a_{21}y^2-a_{21}xy^2+b_{03}y^3\\
&+\delta_1\varepsilon y-\tfrac{b}{2}\varepsilon^{3}(x+x^2)
+ q^+(x,y,\varepsilon),
\end{aligned}\right),&\text{if} \ \ y>0, \\[1.0ex]
&\left(\begin{aligned}&2a_{21}xy-(a_{02}+3b_{03})y^2
+a_{21}x^2y-3b_{03}xy^2+a_{03}y^3\\
&+\delta_1 \varepsilon(x+\tfrac{3}{2}x^2+\tfrac{1}{2}x^3)-\varepsilon^2 y
+p^-(x,y,\varepsilon),\\[0.5ex]
&x+\tfrac{3}{2}x^2+\tfrac{1}{2}x^3-a_{21}y^2-a_{21}xy^2+b_{03}y^3\\
&+\delta_1\varepsilon y+\tfrac{b}{2}\varepsilon^{3}(2+3x+x^2)
+q^-(x,y,\varepsilon),
\end{aligned}\right),&\text{if} \ \ y<0.
\end{aligned}\right.}
\end{equation}
By the approach described in
\cite{TY2018} we introduce the following near-identity state transformation
and time rescaling:
\begin{equation}\label{Eqn46}
\begin{array}{ll}
x\rightarrow x+d_1(\varepsilon)x+d_2(\varepsilon)y, \quad
y\rightarrow y+d_3(\varepsilon)x+d_4(\varepsilon) y, \quad
t\rightarrow t+d_5(\varepsilon)t, \\
\end{array}
\end{equation}
for the upper system in \eqref{Eqn45}, where
\begin{equation}\label{Eqn47}
d_i (\varepsilon)=d_{i1} \varepsilon+d_{i2} \varepsilon^2
+ \cdots + d_{in} \varepsilon^n, \quad i = 1,2, \cdots, 5.
\end{equation}
A similar transformation can be used to simplify the lower system
in \eqref{Eqn45}. Note that the unperturbed system of \eqref{Eqn45} is
unchanged by the identity map \eqref{Eqn46}$|_{\varepsilon=0}$.
Therefore, we may obtain proper $d_i(\varepsilon)$'s to simplify the
perturbations without loss of generality.
By substituting \eqref{Eqn46} into system \eqref{Eqn45}
and taking the $\varepsilon$-order terms, we obtain
\begin{equation}\label{Eqn48}
\begin{array}{ll}
\dot{x}=\!\!\!&-d_{21}x -\frac{1}{2}(3d_{21}-4a_{21}d_{31})x^2
+2 (a_{02} d_{31} - 3 b_{03} d_{31} + a_{21} d_{41} + a_{21} d_{51}
+ p_{211}) xy\\[0.5ex]
&-(a_{02} d_{11} - 3 b_{03} d_{11} - 3 a_{21} d_{21} - 2 a_{02} d_{41}
+ 6 b_{03} d_{41} - a_{02} d_{51} + 3 b_{03} d_{51} - p_{021})y^2\\[0.5ex]
&-\frac{1}{2}(d_{21}-2a_{21}d_{31})x^3+(a_{21} d_{11} - 6 b_{03} d_{31}
+ a_{21} d_{41} + a_{21} d_{51} + p_{211})x^2y\\[0.5ex]
&+(3 a_{21} d_{21} + 3 a_{03} d_{31} - 6 b_{03} d_{41} - 3 b_{03} d_{51}
+ p_{121})xy^2\\[0.5ex]
&-(a_{03} d_{11} + 4 b_{03} d_{21} - 3 a_{03} d_{41} - a_{03} d_{51}
- p_{031}) y^3,
\end{array}
\end{equation}
and
\begin{equation}\label{Eqn49}
\begin{array}{ll}
\dot{y}=\!\!\!&(d_{11}-d_{41}+d_{51})x+d_{21}\,y
+\frac{3}{2}(2d_{11}-d_{41}+d_{51})x^2+(3 d_{21} - 4 a_{21} d_{31} + 2 q_{211})xy\\[0.5ex]
&-(a_{02} d_{31} - 3 b_{03} d_{31} + a_{21} d_{41} + a_{21} d_{51} - q_{021})\, y^2+\frac{1}{2}(3 d_{21} - 6 a_{21} d_{31} + 2 q_{211})x^2y\\[0.5ex]
&-(a_{21} d_{11} - 6 b_{03} d_{31} + a_{21} d_{41} + a_{21} d_{51} - q_{121})xy^2\\[0.5ex]
&-(a_{21} d_{21} + a_{03} d_{31}
 - 2 b_{03} d_{41} - b_{03} d_{51} - q_{031}) y^3
+\frac{1}{2}(3d_{11}-d_{41}+d_{51})x^3.
\end{array}
\end{equation}
Now, simplify setting
$$
\begin{array}{ll}
d_{11}=&\tfrac{1}{2 (a_{02} - 3 b_{03}) (3 b_{03}^2-a_{03} a_{21})}
\big\{a_{02} a_{21} p_{031} - b_{03}\big[3 a_{21} p_{031}
- 6 b_{03} (p_{021} - p_{121}) - 2 a_{02} p_{121}\big]\\[0.5ex]
&-a_{03}\big[2 a_{21} p_{021} + (a_{02} - 3 b_{03}) q_{121}\big]\big\},
\\[0.0ex]
d_{21}=&0, \\
d_{31}=& \tfrac{1}{6 a_{03}
( 3 b_{03}^2 -a_{03} a_{21} )}
(3 a_{21} b_{03} p_{031} + 2 a_{03} a_{21} p_{121}
- 3 a_{03} b_{03} q_{121}),
\\[1.0ex]
d_{41}=& \tfrac{1}{a_{03} (3 b_{03}-a_{02})}(a_{03} p_{021}
- a_{02} p_{031} + 3 b_{03} p_{031}),\\[1.0ex]
d_{51}=& \tfrac{1}{2 a_{03} (a_{02} - 3 b_{03})
(3 b_{03}^2 -a_{03} a_{21})}\big\{12 b_{03}^2 (3 b_{03}-a_{02}) p_{031}
+3 a_{02}a_{03} a_{21} p_{031}\\[0.5ex]
&-a_{03} b_{03} \big[9 a_{21} p_{031} + 2 a_{02} p_{121}
- 6 b_{03} (2 p_{021} + p_{121})\big]
- a_{03}^2 \big[4 a_{21} p_{021}+ (3 b_{03}-a_{02}) q_{121}\big]\big\},
\end{array}
$$
eliminates the terms $y^2$, $xy^2$, $y^3$
in the $\dot{x}$ equation, and the terms $y$, $xy^2$ in the $\dot{y}$ equation.  It implies that the perturbations
\begin{equation}\label{Eqn50}
\sum_{k=1}\varepsilon^k(p_{02k}y^2+p_{03k}y^3+p_{12k}xy^2) \qquad \text{and} \qquad \sum_{k=1}\varepsilon^kq_{12k}xy^2,
\end{equation}
in \eqref{Eqn45} are redundant and can be removed. Thus, we assume that
$p_{02k}=p_{03k}=p_{12k}=q_{12k}=0$.  Further, introducing the transformation,
$(x,y,t)\rightarrow (\varepsilon^3 x, \varepsilon^2 y,
\frac{t}{\varepsilon})$ into system \eqref{Eqn45} yields
\begin{equation}\label{Eqn51}
{\small \left( \begin{array}{cc}\dot{x}\\
\dot{y}
\end{array} \right)
=\left\{\begin{aligned}
&\left(\begin{aligned}&\delta_1 x-y
+ \tfrac{1}{2}\delta_1(3 \varepsilon^3 x^2 + \varepsilon^6 x^3)
+ 2(a_{21} +\sum_{k=1} p_{21k}\varepsilon^k)\varepsilon xy+a_{03}\varepsilon^2 y^3 \\
&+ (a_{02} - 3 b_{03}) y^2
+(a_{21}+\sum_{k=1}p_{21k}\varepsilon^k)\varepsilon^4 x^2y -3 b_{03}\varepsilon^3 xy^2, \\[0.5ex]
&(1-\tfrac{b}{2}\varepsilon^{3})x+\delta_1 y
+\tfrac{1}{2}(3\varepsilon^3-b \varepsilon^{6})x^2
+2 \sum_{k=1}q_{21k} \varepsilon^{k+2} xy \\
&-(a_{21}-\sum_{k=1}q_{02k}\varepsilon^k)\varepsilon y^2
+\tfrac{1}{2}\varepsilon^6x^3 +\sum_{k=1}q_{21k} \varepsilon^{k+5} x^2y \\
&-a_{21}\varepsilon^4 xy^2
+(b_{03} + \sum_{k=1}q_{03k}\varepsilon^k)\varepsilon^3 y^3,\\
\end{aligned}\right)\! ,&\text{if} \ y>0, \\[1.0ex]
&\left(\begin{aligned}&\delta_1 x-y
+\tfrac{1}{2} \delta_1(3 \varepsilon^3 x^2+\varepsilon^6 x^3)
+2( a_{21}+\sum_{k=1}p_{21k}\varepsilon^k)\varepsilon xy+a_{03}\varepsilon^2 y^3\\
&-(a_{02}+3b_{03}) y^2
+(a_{21}+\sum_{k=1}p_{21k}\varepsilon^k)\varepsilon^4 x^2y-3 b_{03} \varepsilon^3 xy^2,\\[0.5ex]
&b+(1+\tfrac{3}{2}b\varepsilon^{3})x+\delta_1 y
+\tfrac{1}{2}(3\varepsilon^3+b \varepsilon^{6})x^2
+2\sum_{k=1} q_{21k}\varepsilon^{k+2} xy \\
& -(a_{21}+\sum_{k=1}q_{02k}\varepsilon^2)\varepsilon y^2
+\tfrac{1}{2}\varepsilon^6x^3
+\sum_{k=1}q_{21k}\varepsilon^{k+5} x^2y
\\
& -a_{21}\varepsilon^4 xy^2
+(b_{03} + \sum_{k=1}q_{03k}\varepsilon^k)\varepsilon^3 y^3,
\end{aligned}\right)\! ,&\text{if} \ y<0.
\end{aligned}\right.
}
\end{equation}
To further simplify system \eqref{Eqn51}, we apply the following scaling,
\begin{equation}\label{Eqn52}
\begin{array}{lllll}
a_{21} = \varepsilon^2 A_{21}, &
a_{02} = \varepsilon^3 A_{02}, &
b_{03} = \varepsilon^3 B_{03}, &
a_{03} = \varepsilon^4 A_{03}, \\[1.0ex]
q_{03k} = \varepsilon^{3-k} Q_{03k}, &
p_{21k} = \varepsilon^{2-k} P_{21k}, &
q_{02k} = \varepsilon^{2-k} Q_{02k}, &
q_{21k} = \dfrac{1}{\varepsilon^{k-1}}\, Q_{21k}, &
\end{array}
\end{equation}
to \eqref{Eqn51}, and obtain
\begin{equation}\label{Eqn53}
{\small \left( \begin{array}{cc}\dot{x}\\
\dot{y}
\end{array} \right)
=\left\{\begin{aligned}
&\left(\begin{aligned}&
-y + \delta_1 x
+ \tfrac{1}{2}\delta_1(3 \varepsilon^3 x^2 + \varepsilon^6 x^3)+ \varepsilon^3 \big[ 2(A_{21} + \sum_{k=1}P_{21k}) xy
\\[0.5ex]
&+ (A_{02} - 3 B_{03}) y^2 \big] + \varepsilon^6 \big[ (A_{21}+\sum_{k=1}P_{21k}) x^2y
-3 B_{03} xy^2
+A_{03} y^3 \big],\\[0.5ex]
& x -\tfrac{b}{2}\varepsilon^{3}x + \delta_1 y
-\tfrac{1}{2} b \varepsilon^{6} x^2+ \varepsilon^3 \big[ \tfrac{3}{2} x^2
+2 \sum_{k=1}Q_{21k} xy
-(A_{21}-\sum_{k=1}Q_{02k} ) y^2 \big] \\
&+ \varepsilon^6 \big[
\tfrac{1}{2} x^3 +\sum_{k=1}Q_{21k} x^2y
-A_{21} xy^2
+(B_{03} + \sum_{k=1}Q_{03k}) y^3 \big],\\
\end{aligned}\right)\! ,&\text{if} \ y>0, \\[2.0ex]
&\left(\begin{aligned}&
-y + \delta_1 x +\tfrac{1}{2} \delta_1(3 \varepsilon^{4} x^2+\varepsilon^{7} x^3) + \varepsilon^3 \big[ 2(A_{21}+\sum_{k=1}P_{21k}) xy\\[0.5ex]
&-(A_{02}+3B_{03}) y^2 \big]+ \varepsilon^6 \big[(A_{21}+\sum_{k=1}P_{21k}) x^2y
-3 B_{03} xy^2+A_{03} y^3 \big],\\[0.5ex]
&x+ b+\tfrac{3 b}{2}\varepsilon^{3} x+\delta_1 y
+\tfrac{b}{2} \varepsilon^{6} x^2 + \varepsilon^3 \big[ \tfrac{3}{2} x^2
+2 \sum_{k=1}Q_{21k} xy
-(A_{21} \\
&+\sum_{k=1}Q_{02k}) y^2 \big]+ \varepsilon^6 \big[
\tfrac{1}{2} x^3 +\sum_{k=1}Q_{21k} x^2y
-A_{21} xy^2
+(B_{03} + \sum_{k=1}Q_{03k}) y^3 \big],
\end{aligned}\right)\! ,&\text{if} \ y<0.
\end{aligned}\right.
}
\end{equation}
Therefore, the singular point $(1,0)$ of system \eqref{Eqn9}
corresponds to the origin of system \eqref{Eqn53}.
To compute the generalized Lyapunov constants $V_k(\varepsilon)$ ($k\ge 2$),
we let $ b=\delta_1 = 0$, which yields $V_1(\varepsilon)=0$.
It is seen from \eqref{Eqn53} that there exists an infinite number of
parameters, but most of them are redundant.  Without loss of generality, we
let
\begin{equation}\label{Eqn54}
\begin{array}{lllll}
P_{21k} = Q_{21k} \!\!\!\!&= Q_{02k} \!\!\!\! & = Q_{03k} = 0, \ \
\forall k\ge 1, \ k \ne 2,
\end{array}
\end{equation}
which, together with $b=\delta_1=0$, are substituted into
\eqref{Eqn53} to yield the final system,
\begin{equation}\label{Eqn55}
\left( \begin{array}{cc}\dot{x}\\
\dot{y}
\end{array} \right)
=\left\{\begin{aligned}
&\left(\begin{aligned}&
-y
+ \varepsilon^3 \big[ 2(A_{21} + P_{212}) xy
+ (A_{02} - 3 B_{03}) y^2 \big] \\[0.5ex]
&+ \varepsilon^6 \big[ (A_{21}+P_{212}) x^2y
- 3 B_{03} xy^2
+A_{03} y^3 \big],\\[0.5ex]
& x + \varepsilon^3 \big[ \tfrac{3}{2} x^2
+2 Q_{212} xy -(A_{21}-Q_{022} ) y^2 \big] \\
&+ \varepsilon^6 \big[ \tfrac{1}{2} x^3 +Q_{212} x^2y - A_{21} xy^2
+(B_{03} + Q_{032}) y^3 \big],\\
\end{aligned}\right)\! ,&\text{if} \ y>0, \\[1.0ex]
&\left(\begin{aligned}&
-y + \varepsilon^3 \big[ 2(A_{21}+P_{212}) xy
-(A_{02}+3B_{03}) y^2 \big] \\[0.5ex]
&+ \varepsilon^6 \big[ (A_{21}+P_{212}) x^2y
-3 B_{03} xy^2 + A_{03} y^3 \big],\\[0.5ex]
&x + \varepsilon^3 \big[ \tfrac{3}{2} x^2
+2 Q_{212} xy -(A_{21}+Q_{022} ) y^2 \big] \\
&+ \varepsilon^6 \big[ \tfrac{1}{2} x^3 +Q_{212} x^2y
- A_{21} xy^2 +(B_{03} + Q_{032}) y^3 \big],
\end{aligned}\right)\! ,&\text{if} \ y<0,
\end{aligned}\right.
\end{equation}
for computing the generalized Lyapunov constants $V_k(\varepsilon) \ (k\ge 2)$.
Returning to system \eqref{Eqn51}, \eqref{Eqn55} becomes
\begin{equation}\label{Eqn56}
{\small \left( \begin{array}{cc}\dot{x}\\
\dot{y}
\end{array} \right)
=\left\{\begin{aligned}
&\left(\begin{aligned}&- y + \delta_1
\big( x + \tfrac{3}{2} \varepsilon^{3} x^2
+ \tfrac{1}{2} \varepsilon^{6} x^3 \big)
+ 2\, \varepsilon\, (a_{21} + \varepsilon^2 p_{212}) xy
\\
& + (a_{02} - 3 b_{03} ) y^2
+ \varepsilon^4 (a_{21} + \varepsilon^2 p_{212} ) x^2y
- 3\, \varepsilon^3 b_{03} xy^2 + \varepsilon^2 a_{03} y^3 , \\[0.5ex]
&x -\tfrac{b}{2}\varepsilon^{3} ( x + \varepsilon^{3} x^2)
+\delta_1 y + \tfrac{3}{2} \varepsilon^3 x^2
+2 \varepsilon^4 q_{212} xy - \varepsilon (a_{21}- \varepsilon^2 q_{022}) y^2\\
&+  \tfrac{1}{2} \varepsilon^6 x^3
+ \varepsilon^7 q_{212} x^2y - \varepsilon^4 a_{21} xy^2
+ \varepsilon^3 (b_{03} + \varepsilon^2 q_{032}) y^3 ,\\
\end{aligned}\right)\! ,&\text{if} \ y>0, \\[1.0ex]
&\left(\begin{aligned}&-y
+ \delta_1
\big( x + \tfrac{3}{2} \varepsilon^{3} x^2
+ \tfrac{1}{2} \varepsilon^{6} x^3 \big)
+ 2\, \varepsilon (a_{21}+  \varepsilon^2 p_{212}) xy
\\
& -(a_{02}+3 b_{03}) y^2
+ \varepsilon^4 (a_{21}+ \varepsilon^2 p_{212}) x^2y
-3 \varepsilon^3 b_{03} xy^2 + \varepsilon^2 a_{03} y^3 ,\\[0.5ex]
&x + b \!+\! \tfrac{b}{2} \varepsilon^{3} (3 x \!+\! \varepsilon^{3} x^2)
\!+\! \delta_1 y
\!+\!  \tfrac{3}{2} \varepsilon^3  x^2
\!+\! 2 \varepsilon^4 q_{212} xy - \varepsilon (a_{21} \!+\!
\varepsilon^2 q_{022})y^2 \\
&+ \tfrac{1}{2} \varepsilon^6 x^3 + \varepsilon^7 q_{212} x^2y
- \varepsilon^4 a_{21} xy^2 + \varepsilon^3 (b_{03}+\varepsilon^2 q_{032})y^3,
\end{aligned}\right)\! ,&\text{if} \ y<0.
\end{aligned}\right.
}
\end{equation}
We will show that system \eqref{Eqn56} can have at least
$10$ small-amplitude limit cycles around the origin.

Now, system \eqref{Eqn55} has only $8$
independent parameters: $Q_{022}$, $Q_{032}$, $Q_{212}$,
$A_{03}$, $A_{21}$, $B_{03}$, $A_{02}$ and $P_{212}$, which will be
used to solve the Lyapunov equations $V_k(\varepsilon)=0$.
Note that later we will perturb the constant term $b$
to obtain one more limit cycle by using the pseudo-Hopf bifurcation.
The 1st generalized Lyapunov constant for
system \eqref{Eqn55} is given by
$V_1(\varepsilon)=2 \pi \delta_1=0$ when $\delta_1=0$.
Then, we have the 2nd generalized Lyapunov constant,
\begin{equation*}
\begin{aligned}
V_2(\varepsilon)=&\frac{8}{3}\varepsilon^3\, Q_{022} .
\end{aligned}
\end{equation*}
Setting $Q_{022}= 0$ yields $V_2(\varepsilon)=0$ and
the 3rd generalized Lyapunov constant,
\begin{equation*}
\begin{aligned}
V_3(\varepsilon)=&\frac{1}{4}\pi\varepsilon^6
\big[ 3 Q_{032} - 2(1- A_{21}) Q_{212} - 6 B_{03} P_{212} \big].
\end{aligned}
\end{equation*}
We obtain
\begin{equation*}
\begin{aligned}
Q_{032}=\frac{2}{3} (1 - A_{21})\, Q_{212}
+ 2 B_{03} P_{122}
\end{aligned}
\end{equation*}
by solving $V_3(\varepsilon)=0$.
Further, we have the 4th generalized Lyapunov constant,
\begin{equation*}
\begin{aligned}
V_4(\varepsilon)=&\frac{16}{45} \varepsilon^9 A_{02}
\big[ (3 + P_{212} ) Q_{212} - 12 B_{03} P_{212} \big].
\end{aligned}
\end{equation*}
Since the condition required for the center condition VI
does not allow $a_{02} =0$, we have $A_{02} \ne 0$.
Thus, solving $V_4(\varepsilon)=0$ for $Q_{212}$ we obtain
\begin{equation*}
\begin{aligned}
Q_{212}=& \dfrac{12 B_{03} P_{212}}
{3+2 P_{212}}, \quad  (3+2 P_{212} \ne 0).
\end{aligned}
\end{equation*}
Then, solving the 5th generalized Lyapunov constant equation
$V_5(\varepsilon)=0$ for $A_{03}$ yields
\begin{equation*}
\begin{aligned}
A_{03}=& \dfrac{1}{3 (2 P_{212}-1) (3+2 P_{212})^2}
\big\{(3 + 2 P_{212})^2
\big[15 A_{02}^2+2 (A_{21}+P_{212}+1)
\big( A_{21}(10 P_{212} +3)
\\
& +4P_{212}^2+4 P_{212}+9 \big) \big]
     -72 B_{03}^2 (6 P_{212} +5)
\big(8 A_{21} P_{212} +2 P_{212}^2-P_{212}+6 \big) \big\},
\end{aligned}
\end{equation*}
leading to
\begin{equation*}
\begin{aligned}
V_6(\varepsilon) = &
- \frac{128 \varepsilon^{15} A_{02} B_{03} P_{212}}
{525 (2 P_{212} -1) (3+2 P_{212})^3} \, V_{6a}, \\[0.5ex]
V_7(\varepsilon) = &
- \frac{\varepsilon^{18} B_{03} P_{212}}
{50400 (2 P_{212} -1)^2 (3+2 P_{212})^5} \, V_{7a}, \\[0.5ex]
V_8(\varepsilon) = &
- \frac{\varepsilon^{21} B_{03} P_{212}}
{6350400 (2 P_{212} -1)^2 (3+2 P_{212})^5} \, V_{8a}, \\[0.5ex]
V_9(\varepsilon) = &
- \frac{\varepsilon^{24} B_{03} P_{212}}
{76204800 (2 P_{212} -1)^3 (3+2 P_{212})^7} \, V_{9a}, \\[0.5ex]
V_{10}(\varepsilon) = &
- \frac{\varepsilon^{27} B_{03} P_{212}}
{76204800 (2 P_{212} -1)^3 (3+2 P_{212})^7} \, V_{10a}, \\[0.5ex]
\end{aligned}
\end{equation*}
where $V_{6a}$, $V_{7a}$, $V_{8a}$, $V_{9a}$, and $V_{10a}$ are polynomials in
$ A_{21}$, $B_{03}^2$, $A_{02}$ and $P_{212}$.
In particular,
\begin{equation*}
\begin{aligned}
V_{6a} = & \ 216 B_{03}^2 \big[ 24 A_{21} P_{212}
\big(4 P_{212}^2-20 P_{212} -15 \big)
+16 P_{212}^4-224 P_{212}^3-280 P_{212}-265 \big]
\\
&+(3+2 P_{212})^2
\big[ 405 A_{02}^2+A_{21}
\big(864 A_{21} P_{212} +1148 P_{212}^2+728 P_{212} +645 \big)
\\
&+(P_{212} +1) (284 P_{212}^2-136 P_{212} +645) \big],
\end{aligned}
\end{equation*}
and other lengthy expressions of $V_{7a}$, $V_{8a}$, $V_{9a}$, and $V_{10a}$
are omitted for brevity. Here, we assume that
\begin{equation}\label{Eqn57}
(2 P_{212} -1) (3+2 P_{212}) \ne 0.
\end{equation}
Then, we use Maple with its built-in command {\it eliminate},
$$
{\rm eliminate}(\{V_{6a},V_{7a},V_{8a},V_{9a}\},A_{21}),
$$
to obtain a solution $A_{21} = -\,\frac{A_{21n}}{A_{21d}}$, where
$$
\begin{aligned}
A_{21n}=&\ 2539579392 B_{03}^6 P_{212}^2
( 16 P_{212}^4-224 P_{212}^3-280 P_{212}-265 )
\\[0.5ex]
&+279936 B_{03}^4 P_{212}\big[ 17010 A_{02}^2 P_{212} (3+2 P_{212})^2
   +37408 P_{212}^6+185392 P_{212}^5
\\[0.5ex]
&+265584 P_{212}^4
  +266792 P_{212}^3+424802 P_{212}^2+207315 P_{212}-130248 \big]
\\[0.5ex]
&+54 B_{03}^2 (3+2 P_{212})^2
 \big[ 243 A_{02}^2 P_{212} \big(20932 P_{212}^2-26676 P_{212}
+13977 \big)
\\[0.5ex]
&-322112 P_{212}^6+853408 P_{212}^5+4354704 P_{212}^4
   +4303208 P_{212}^3+1241408  P_{212}^2
\\[0.5ex]
&-18636 P_{212} -9540 \big]+81 A_{02}^2 (3+2 P_{212})^4 ( 46388 P_{212}^2
+89553 P_{212}+45 )
\\[0.5ex]
& -5 (3+2 P_{212})^4 (2 P_{212} -1)
(P_{212} +1) (7 P_{212} +3) (34 P_{212} +3) (58 P_{212} +129)
\end{aligned}
$$
and
$$
\begin{aligned}
A_{21d}=&\ 60949905408 B_{03}^6 P_{212}^3
(4 P_{212}^2-20 P_{212}-15)
\\[0.5ex]
&+5038848 B_{03}^4 P_{212}^2 (3+2 P_{212})^2
(1484 P_{212}^2+3288 P_{212} -1039)
\\[0.5ex]
&+3888 B_{03}^2 P_{212} (3+2 P_{212})^2
( 122472 A_{02}^2 P_{212} -3416 P_{212}^4
     +13700 P_{212}^3
\\[0.5ex]
& +34622 P_{212}^2+18267 P_{212} +1512 )
+118098 A_{02}^2 P_{212} (3+2 P_{212})^4
\\[0.5ex]
&-5 (3+2 P_{212})^4 (2 P_{212} -1) (7 P_{212} +3)
(34 P_{212} +3) (58 P_{212} +129),
\end{aligned}
$$
and three polynomial resultant equations,
\begin{equation*}
R_k = R_k(B_{03}^2, A_{02}, P_{212}) = 0, \quad k=1,2,3.
\end{equation*}
Further, we apply the command {\it eliminate} again to
the three resultants to obtain a solution
$B_{03}^2 = B_{03}^2(A_{02}, P_{212})$, and two lengthy
polynomial resultants,
\begin{equation*}
\begin{array}{ll}
R_{12} = A_{02}\, P_{212} (3+2 P_{212}) R_{12a} R_{12b},
\quad
R_{13} = A_{02}\, P_{212} (3+2 P_{212}) R_{13a} R_{13b},
\end{array}
\end{equation*}
where $R_{12a}$, $ R_{12b}$, $R_{13a}$ and $ R_{13b}$ are
lengthy polynomials in $A_{02}$ and $P_{212}$, with
$245$, $5880$, $3086$ and $19498$ terms, respectively.
There are four combinations: $(R_{12a},R_{13a})$,
$(R_{12a},R_{13b})$, $(R_{12b},R_{13a})$ and $(R_{12b},R_{13b})$
to find the solutions for $(A_{02},P_{212})$.
However, it can be shown that only the combination
$(R_{12a},R_{13b})$ generates maximal number of limit cycles.
Finally, we use the Maple-built command {\it resultant},
$$
{\rm resultant}(R_{12a},R_{13a},A_{02}),
$$
to obtain the following resultant,
\begin{equation}\label{Eqn58}
\begin{aligned}
\!\!\!\! R_{1213} = & \ \tilde{C} \pi^{14} P_{212}^{70}
(2 P_{212} -1)^{20} (3+2 P_{212})^{36} \\
& \times
\big[(P_{212} +2)^6
( 18000 P_{212}^4+111504 P_{212}^3+182360 P_{212}^2
+98444 P_{212} +9965 )
\\
&\quad \
 \times \! ( 33888 P_{212}^6+237600 P_{212}^5
+463840 P_{212}^4 +311472 P_{212}^3
\\
&\qquad \ +\! 75078 P_{212}^2+7338 P_{212} +243 ) \big]
\, R_{1213a} R_{1213b},
\end{aligned}
\end{equation}
where $\tilde{C}$ is a big positive integer,
$R_{1213a}$
and $R_{1213b}$ are respectively
$54$th-degree and $1080$th-degree polynomials in $P_{212}$.
With the aid of Maple using $1000$ digits accuracy, we obtain
$9$, $22$ and $221$ real solutions from
the terms in the script bracket of $R_{1213}$, the polynomials
$R_{1213a}$ and $R_{1213b}$, respectively.

\begin{remark}\label{Rem5.1}
In most research articles, researchers prefer to not use numerical
computation in a rigorous mathematical proof.
Here, we should point out that (1) our numerical computation is
only used at the last step, and all computations in the previous
steps are symbolic and exact. (2) we can, instead of using numerical
computation, apply the {\it interval computation}
(or {\it interval arithmetic}) \cite{Hickey2001} to prove the existence
of the solutions in an interval satisfying a required accuracy. This
does not really change the basics of computation, and
a numerical computation makes the presentation simple and clear.
(3) the most important point is
that the accuracy to be taken (e.g., $1000$ decimal points)
should be good enough to guarantee the reliability of proof.
That is, the accumulation of round errors in the numerical
computation does not affect the conclusion for the existence of solutions.
In other words, it is not a matter of numerical computation, but the
high enough accuracy.
\end{remark}

It can be shown that all the $221$ solutions obtained from
$R_{1213b}$ produces
$V_1=V_2=\cdots=V_8=0$, but $V_9 \ne 0$, giving $8 \times 2 = 16 $
small amplitude limit cycles.
While from the $9+22$ solutions, we obtain $8$ solutions
(with $\delta_1 =Q_{022} =0$):
$$
\begin{array}{rlrlrlrl}
&(P_{212},\,A_{02},\,B_{03},\,A_{21},\,A_{03},\, Q_{212},\,Q_{032})
& & &
\\[0.5ex]
= \!\!\!\!
& (\hspace*{0.145in} -84.39736368 \cdots,\, 246.66126379 \cdots,
\hspace*{0.162in} \pm 2.18680584 \cdots, \hspace*{0.10in}30.37184157 \cdots,\\
&\ -1169.84368584 \cdots , \pm 13.35825220 \cdots, \mp 630.69227526 \cdots),\\
& (\hspace{0.22in} -2.12093549 \cdots,\hspace*{0.194in} 0.68942957 \cdots,
\hspace*{0.152in} \pm 0.05796919 \cdots,
\hspace*{0.175in} 1.33657024 \cdots, \\
&\hspace*{0.40in} 0.13502593 \cdots , \hspace*{0.072in}\pm 1.18803564 \cdots,
\hspace*{0.152in} \mp 0.51246946 \cdots), \\
& (\hspace{0.22in} -1.64608697 \cdots,\hspace*{0.194in} 2.82165606 \cdots,
\hspace*{0.152in} \pm 0.10278968 \cdots,
-20.81805724 \cdots, \\
&\hspace*{0.278in} -0.15805791 \cdots, \hspace*{0.085in}\pm 6.94931604 \cdots,
\mp 100.74198202 \cdots), \\
& (\hspace{0.22in} -0.95226778 \cdots,\hspace*{0.194in} 0.20225860 \cdots,
\hspace*{0.152in} \pm 0.06175799 \cdots,\hspace{0.08in} -0.20115348 \cdots, \\
&\hspace*{0.40in} 0.09330998 \cdots, \hspace*{0.085in}\mp 0.64422154 \cdots,
\hspace*{0.152in} \mp 0.63349293 \cdots), \\
\end{array}
$$
which result in $V_2=V_3=\cdots=V_9=0$, but $V_{10} \ne 0$.
Moreover, all the $8$ solutions satisfy the requirement of the
center condition VI,
$$
a_{02} +|a_{12}|
= a_{02} + |-3 b_{03} |
= A_{02} \varepsilon^3 + | 3 B_{03} \varepsilon^3 | < 0
\ \ \Longleftrightarrow \ \
A_{02} - | 3 B_{03}| >0 \ \ (\varepsilon<0).
$$
We choose one of the $8$ solutions,
$$
\begin{array}{rlrlrlrl}
&(P_{212},\,A_{02},\,B_{03},\,A_{21},\,A_{03},\, Q_{212},\,Q_{032},\,
Q_{022},\,\delta_1)_{\rm C}
& & &
\\[1.0ex]
= \!\!\!\!
&( -2.12093549 \cdots,\ 0.68942957 \cdots, \hspace*{0.15in} 0.05796919 \cdots,
\ 1.33657024 \cdots, \\
&\hspace*{0.18in} 0.13502593 \cdots , \ 1.18803564 \cdots, \,
-0.51246946 \cdots, \ 0.0, \hspace*{0.80in} 0.0), \\
\end{array}
$$
where the subscript C indicates a critical point, under which
$$
V_1= V_2 = V_3 = \cdots = V_9 = 0, \quad V_{10} =
0.05995729 \cdots \varepsilon^{27} < 0 \ \ (\varepsilon<0).
$$
More precisely, with the accuracy of $1000$ decimal points, we have
$$
\begin{array}{rlll}
V_1=\!\!\!\!&0.0\, , \quad
V_2=0.0\, \varepsilon^3, \quad
V_3=0.0\, \varepsilon^6, \\[0.5ex]
V_4=\!\!\!\!& 0.63836071 \cdots \times 10^{-999}\, \varepsilon^9,
\hspace*{0.20in}
V_5=\ \ \ 0.26153383 \cdots \times 10^{-997}\, \varepsilon^{12}, \\[0.5ex]
V_6=\!\!\!\!&0.62071107 \cdots \times 10^{-970}\, \varepsilon^{15}, \quad
V_7=-\,0.32934572 \cdots \times 10^{-969}\, \varepsilon^{18}, \\[0.5ex]
V_8=\!\!\!\!&0.10608263 \cdots \times 10^{-968}\, \varepsilon^{21}, \quad
V_9=-\,0.27441458 \cdots \times 10^{-968}\, \varepsilon^{24},
\end{array}
$$
for which we are definitely confident that a solution exists, satisfying
$V_1=V_2 = \cdots = V_9=0$.
In addition, a direct calculation shows that
$$
\det \left[\frac{\partial(V_1(\varepsilon),V_2(\varepsilon),V_3(\varepsilon),
V_4(\varepsilon),V_5(\varepsilon),V_6(\varepsilon),V_7(\varepsilon),
V_8(\varepsilon),V_9(\varepsilon))}
{\partial(\delta_1,Q_{022},Q_{032},Q_{212},A_{03},A_{21},B_{03},
A_{02},P_{212})}\right]_{\rm C}
\!\!\! = 329.04383261 \cdots \, \varepsilon^{108} > 0.
$$
Hence, by Lemma \ref{Lem3.1} we know that system \eqref{Eqn53} has at least
$9$ small-amplitude limit cycles bifurcating from the origin, leading to
the existence of
$18$ small-amplitude limit cycles bifurcating from
the bi-center of such $Z_2$-equivariant cubic switching systems.

Further, applying the pseudo-Hopf bifurcation,
we obtain one more small-amplitude limit cycle
for system \eqref{Eqn56} by perturbing the coefficient $b$.
For $b$ and $\varepsilon$ sufficiently small, the switching
system \eqref{Eqn56} has a small sliding segment on the switching
manifold $y=0$ with the end points at
$(0,0)$ and $(\varrho_b,0)$,  where
$\varrho_b$ is the root of the equation $g^-(x,0)=0$ in \eqref{Eqn56},
namely,
$$
x+b+ \dfrac{b}{2}\varepsilon^3 (3x + \varepsilon^3 x^2)
+ \dfrac{3}{2}\varepsilon^3 x^2 + \dfrac{1}{2} \varepsilon^6 x^3
= \dfrac{1}{2}(b+x)(1+ \varepsilon^3 x)(2+ \varepsilon^3 x)=0,
$$
which yields the solution $\varrho_b=-b$.
Thus, the sliding segment shrinks to $(0,0)$ when $b$ goes to zero.

We consider the point $(\varrho_\varepsilon,0)$ on the switching
line $y=0$ with sufficiently small $\varrho_b<\varrho_\varepsilon$,
where $\varrho_\varepsilon$ is defined as a bifurcation function, satisfying
\begin{equation}\label{Eqn59}
d(\varrho_\varepsilon,\varepsilon)=\Upsilon^+
(\varrho_\varepsilon,\varepsilon)
-(\Upsilon^-)^{-1}(\varrho_\varepsilon,\varepsilon).
\end{equation}
Then, we have two half-return maps given in the form of
\begin{equation}\label{Eqn60}
\Upsilon^+(\varrho_\varepsilon,\varepsilon)
=V_1^+(\varepsilon) \varrho_\varepsilon+O(\varrho_\varepsilon^2) \quad
\text{and} \quad
(\Upsilon^-)^{-1}(\varrho_\varepsilon,\varepsilon)
=V_0^-(\varepsilon)+V_1^-(\varepsilon) \varrho_\varepsilon
+O(\varrho_\varepsilon^2),
\end{equation}
respectively, where
\begin{equation}\label{Eqn61}
V_0^-(\varepsilon)=b\,[2+O(\varepsilon^3)] \quad \text{and} \quad
V_1^{\pm}(\varepsilon)=\pm\delta_1\pi+O(\varepsilon^3).
\end{equation}
It follows from \eqref{Eqn59} and \eqref{Eqn60} that
$V_{0}(\varepsilon)=-b\,[2+O(\varepsilon^3)]$
and $V_{1}(\varepsilon)|_{b=0}=2\delta_1\pi$,
implying that  $V_{0}(\varepsilon)=V_{1}(\varepsilon)=0$
when $b=\delta_1=0$. Then, to get more small-amplitude limit cycles,
we let $b=\delta_1=0$, and compute higher-order generalized
Lyapunov constants. This means that we return to exactly where we start
and obtain the system \eqref{Eqn55} under the assumption $b=\delta_1=0$.
Therefore, by Lemma \ref{Lem3.1} and combining the
$9$ small-amplitude limit cycles obtained above, we know that
$10$ small-amplitude crossing limit cycles exist near the origin of
system \eqref{Eqn56}, and so
$20$ small-amplitude limit cycles can bifurcate from
the bi-center of the $Z_2$-equivariant cubic switching system
\eqref{Eqn9}.

Note that since we include all $\varepsilon^k$-order
perturbations in constructing system \eqref{Eqn41}, the number $20$ is actually
the maximal number of small-amplitude limit cycles which can be obtained
with the $6$ center conditions I-VI under Hopf and
pseudo-Hopf bifurcations for the $Z_2$-equivariant cubic switching system
considered in this paper.

This finishes the proof for Theorem \ref{Thm2.3}.

\section{Conclusion}

In this paper, we have studied the bi-center problem and
the cyclicity problem for planar switching nilpotent systems.
First, we generalize the Poincar\'e-Lyapunov method for
switching systems with linear-type centers to compute
the generalized Lyapunov constants for switching nilpotent systems.
Then, we derive $6$ bi-center conditions for cubic $Z_2$-equivariant
switching systems with two symmetric nilpotent singular points.
In particular, we find that the bi-center $(\pm 1, 0)$
consists of the combination of two second-order nilpotent cusps.
Further, we construct perturbed systems with the $6$ center conditions
and present one with the center condition VI to show the existence
at least $20$ small-amplitude limit cycles
bifurcating from the nilpotent bi-center,
which provides a great improvement
from $12$, leading to a new lower bound on the number of limit cycles
in $Z_2$-equivariant cubic switching systems.

Motivated from this work,
some questions naturally arise, which may promote future research
in this direction.
\begin{enumerate}
\item
In this work, we use some special perturbation to obtain
$6$ center conditions for $Z_2$-equivariant cubic switching
system with nilpotent singular points. Does this special
perturbation generate all possible center conditions?
If not, what center conditions might be missed,
and what kind of perturbations should be used to find all
possible center conditions?

\item
The computation of the Lyapunov constants heavily depends upon
a computer algebra system and the algorithms to be used.
From this work, we have found that the current techniques for symbolic
computation is good enough for computing the generalized Lyapunov
constants. However, the current well-known approaches such as
Groebner basis, regular chains, and Maple built-in programs
{\it eliminate}, {\it resultant}, etc. seem still not powerful enough
to solve the Lyapunov constant (polynomial) equations as well as
many other such equations related to Hilbert's 16th problem.
So, how to develop a more efficient method (or algorithm, or program)
for solving the system of multivariate polynomial equations
is a very challenging task and needs further research.

\item
In this paper, our generalized Lyapunov constants contain different
$\varepsilon^k$ order terms (before scaling), and we treat them
in a consistent way by scaling.
However, in the literature, many works treat such Lyapunov
constants or focus values
according to the $\varepsilon^k$ orders, similar to
the consideration of the 1st-order, 2nd-order, etc. Melnikov functions
(for example, see \cite{TY2018}).
Then, what are the advantages and disadvantages of these two perturbation
methods, and for what systems one of the two approaches is better to be
applied?
\end{enumerate}

\section*{Acknowledgment}

This work was supported by the National Natural Science
Foundation of China Nos.~12071091 (T. Chen), 12071198 (F. Li)
and 12371175 (Y. Tian), the Guangdong Basic and Applied Basic Research
Foundation No.~2022A1515011827 (T. Chen), the Project funded by China Postdoctoral Science Foundation No.~2023M744329 (T. Chen) and
the Natural Sciences and Engineering Research Council of Canada,
No.~R2686A02 (P. Yu).


\vspace{0.10in}
\begin{thebibliography}{99}

\bibitem{Andronov}
A. Andronov, A. Vitt, S. Khaikin,
Theory of Oscillations,
Pergamon Press, Oxford, 1966.

\bibitem{Algaba}
A. Algaba, C. Garc\'ia, J. Gin\'e, J. Llibre,
The center problem for $Z_2$-symmetric nilpotent vector fields,
J. Math. Anal. Appl. 466 (2018) 183--198.

\bibitem{Banerjee}
S. Banerjee, G. Verghese,
Nonlinear Phenomena in Power Electronics: Attractors, Bifurcations,
Chaos, and Nonlinear Control,
Wiley-IEEE Press, New York, 2001.

\bibitem{Bastos}
J. Bastos, C. Buzzi, J. Llibre, D. Novaes,
Melnikov analysis in nonsmooth differential systems
with nonlinear switching manifold,
J. Differential Equations 267 (2019) 3748--3767.

\bibitem{Bernardo}
M.D. Bernardo, P. Kowalczyk, A.B. Nordmark,
Sliding bifurcations: a novel mechanism for the sudden onset
of chaos in dry friction oscillators,
Int. J. Bifurcation Chaos 13 (2003) 2935--2948.

\bibitem{Braun}
F. Braun, L. P. C. da Cruz, J. Torregrosa,
On the number of limit cycles in piecewise planar quadratic
differential systemss, (2023) preprint.

\bibitem{ADRIANA}
A. Buic$\rm{\breve{a}}$, J. Llibre, O. Makarenkov,
Asymptotic stability of periodic solutions for nonsmooth differential
equations with application to the nonsmooth van der Pol oscillator,
SIAM J. Math. Anal. 40 (2009) 2478--2495.

\bibitem{Castillo}
J. Castillo, J. Llibre, F. Verduzco,
The pseudo-Hopf bifurcation for planar discontinuous piecewise linear
differential systems,
Nonlinear Dynam. 90 (2017) 1829--1840.

\bibitem{X.W1}
X. Chen, Z. Du,
Limit cycles bifurcate from centers of discontinuous quadratic systems,
Comput. Math. Appl. 59 (2010) 3836--3848.

\bibitem{X.W3}
X. Chen, W. Zhang,
Isochronicity of centers in switching Bautin system,
J. Differential Equations 252 (2012) 2877--2899.

\bibitem{Chen4}
T. Chen, L. Huang, P. Yu,
Center condition and bifurcation of limit cycles for quadratic
switching systems with a nilpotent equilibrium point,
J. Differential Equations 303 (2021) 326--368.

\bibitem{Chen3}
T. Chen, L. Huang, P. Yu,  W. Huang,
Bifurcation of limit cycles at infinity in piecewise polynomial systems,
Nonlinear Anal.: Real World Appl. 41 (2018) 82--106.

\bibitem{T.Chen2}
T. Chen, S. Li, J. Llibre,
$Z_2$-equivariant linear type bi-center cubic polynomial Hamiltonian
vector fields,
J. Differential Equations  269 (2020) 832--861.

\bibitem{Colak2}
I. Colak, J. Llibre and C. Valls,
Hamiltonian nilpotent centers of linear plus cubic homogeneous
polynomial vector fields,
Adv. Math. 259 (2014) 655--687.

\bibitem{Cruz}
L. Cruz, D. Novaes, J. Torregrosa,
New lower bound for the Hilbert number in piecewise quadratic
differential systems,
J. Differential Equations 266 (2019) 4170--4203.

\bibitem{Dumortier2006}
F. Dumortier, J. Llibre, J. Art\'es,
Qualitative Theory of Planar Differential Systems,
Universitext, Springer-Verlag, New York, 2006.

\bibitem{ELT2023}
R.D. Euz\'ebio, J. Llibre, D.J. Tonon,
Lower bounds for the number of limit cycles in a generalised
Rayleigh-Li\'enard oscillator,
Nonlinearity 35 (2022) 3883--3906.

\bibitem{Filippov}
A.F. Filippov,
Differential Equation with Discontinuous Right-Hand Sides,
Kluwer Academic, Netherlands, 1988.

\bibitem{Freire3}
E. Freire, E. Ponce, F. Torres,
A general mechanism to generate three limit cycles in planar
Filippov systems with two zones,
Nonlinear Dynam. 78 (2014) 251--263.

\bibitem{Garcia}
I. Garc\'ia,
Cyclicity of some symmetric nilpotent centers,
J. Differential Equations 260 (2016) 5356--5377.

\bibitem{Gasull}
A. Gasull,  J. Torregrosa,
Center-focus problem for discontinuous planar differential equations,
Int. J. Bifurcation Chaos  13 (2003) 1755--1765.

\bibitem{Giacomini}
H. Giacomini, J. Gin\'e, J. Llibre,
The problem of distinguishing between a center and a focus
for nilpotent and degenerate analytic systems,
J. Differential Equations 227 (2006) 406--426.

\bibitem{Gouveia}
L. Gouveia, J. Torregrosa,
Local cyclicity in low degree planar piecewise polynomial vector fields,
Nonlinear Anal.: Real World Appl. 60 (2021) 103278 1-19.

\bibitem{L.Guo2018}
L. Guo, P. Yu, Y. Chen,
Bifurcation analysis on a class of $Z_2$-equivariant cubic
switching systems showing eighteen limit cycles,
J. Differential Equations 266 (2019) 1221--1244.

\bibitem{Han2}
M. Han, P. Yu,
Normal Forms, Melnikov Functions and Bifurcations of Limit Cycles,
Springer-Verlag, New York, 2012.

\bibitem{Han1}
M. Han, W. Zhang,
On Hopf bifurcation in non-smooth planar system,
J. Differential Equations 248 (2010) 2399--2416.

\bibitem{Hickey2001}
T. Hickey, Q. Ju, M.~H. Van Emden,
Interval arithmetic: From principles to implementation,
J. ACM 48 (2001) 1038--1068.

\bibitem{Huan}
S. Huan, X. Yang,
Existence of limit cycles in general planar piecewise linear systems
of saddle-saddle dynamics,
Nonlinear Anal. 92 (2013) 82--85.

\bibitem{Kunze}
M. Kunze, T. Kupper,
Qualitative bifurcation analysis of a non-smooth friction oscillator model,
Math. Phys. 48 (1997) 87--101.

\bibitem{Leine}
R. Leine, H. Nijmeijer,
Dynamics and Bifurcations of Non-Smooth Mechanical Systems,
Lect. Notes Appl. Comput. Mech., vol.18, Springer-Verlag, Berlin, 2004.

\bibitem{Li}
F. Li, H. Li, Y. Liu,
New double bifurcation of nilpotent focus,
Int. J. Bifurcation Chaos 31 (2021) 2150053.

\bibitem{Li2}
F. Li, Y. Liu, Y. Liu, P. Yu,
Bi-center problem and bifurcation of limit cycles from nilpotent
singular points in $Z_2$-equivariant cubic vector fields,
J. Differential Equations 265 (2018) 4965--4992.

\bibitem{Li3}
F. Li, Y. Liu, P. Yu, J. Wang,
Complex integrability and linearizability of cubic
$Z_2$-equivariant systems with two 1:q resonant singular points,
J. Differential Equations 300 (2021) 786--813.

\bibitem{Li4}
F. Li, Y. Jin, Y. Tian, P. Yu,
Integrability and linearizability of cubic $Z_2$ systems
with non-resonant singular points,
J. Differential Equations 269 (2020) 9026--9049.

\bibitem{LY2015}
F. Li, P. Yu, Y. Tian, Y. Liu,
Center and isochronous center conditions for switching systems
associated with elementary singular points,
Commun. Nonlinear Sci. Numer. Simul. 28 (2015) 81--97.

\bibitem{J.Li2003}
J. Li, Hilbert's 16th problem and bifurcations of planar
polynomial vector fields,
Int. J. Bifurcation Chaos  3 (2013) 47--106.

\bibitem{J.Li2010}
J. Li, Y. Liu,
New results on the study of $Z_{q}$-equivariant planar
polynomial vector fields,
Qual. Theo. Dyna. Syst. 9 (2010) 167--219.

\bibitem{tian2022}
J. Li, Y. Tian, Y. Liu,
Twenty-six crossing limit cycles around one
singular point in a cubic switching system,
Int. J. Bifurcation Chaos  32(10) (2022) 2250158 (7 pages).

\bibitem{Liu2014}
Y. Liu, F. Li,
Double bifurcation of nilpotent focus,
Int. J. Bifurcation Chaos 25(3) (2015) 1550036 (10 pages).

\bibitem{Lv}
Y. Lv, R. Yuan, Y. Pei,
Dynamics in two nonsmooth predator-prey models with threshold harvesting,
Nonlinear Dynam. 74 (2013) 107--132.

\bibitem{Minorsky}
N. Minorsky,
Nonlinear Oscillations,
Van Nostrand, New York, 1962.

\bibitem{TIAN}
Y. Tian,  P. Yu,
Center conditions in a switching Bautin system,
J. Differential Equations 259 (2015) 1203--1226.

\bibitem{TY2018}
Y. Tian, P. Yu,
Bifurcation of small limit cycles in cubic integrable systems using
higher-order analysis,
J. Differential Equations 264(9) (2018) 5950--5976.

\bibitem{Strozyna}
E. Str\'{o}\.{z}zyna, H. {\rm \.{Z}o{\l}\c{a}dek},
The analytic normal for the nilpotent singularity,
J. Differential Equations 179 (2012) 479--537.

\bibitem{Yang}
J. Yang, L. Zhao,
The cyclicity of period annuli for a class of cubic Hamiltonian
systems with nilpotent singular points,
J. Differential Equations 263(9) (2017) 5554--5581.


\bibitem{Yu2021}
P. Yu, M. Han, X. Zhang,
Eighteen limit cycles around two symmetric foci in a cubic planar
switching polynomial system,
J. Differential Equations 275(2) (2021) 939--959.

\end{thebibliography}
\end{document}